\documentclass{article}
\usepackage{amsmath}
\usepackage{amsthm}
\usepackage{graphics}
\usepackage{epsfig}
\usepackage{stmaryrd}

\newtheorem{theorem}{Theorem}[section]
\newtheorem{cor}[theorem]{Corollary}
\newtheorem{lem}[theorem]{Lemma}
\newtheorem{prop}[theorem]{Proposition}

\theoremstyle{definition}

\newtheorem{rem}[theorem]{Remark}

\numberwithin{equation}{section}

\usepackage{amssymb}
\usepackage{amsfonts}
\usepackage{psfrag}
\usepackage{url}
\usepackage{mathtools}

\usepackage[makeroom]{cancel}

\usepackage{stmaryrd}
\DeclareMathOperator*{\bigtimes}{\vartimes}

\newcommand{\K}{\mathbb{K}}
\newcommand{\N}{\mathbb{N}}
\newcommand{\R}{\mathbb{R}}

\newcommand{\Q}{\mathbb{Q}}
\newcommand{\C}{\mathbb{C}}

\newcommand{\cC}{\mathcal{C}}
\newcommand{\cH}{\mathcal{H}}
\newcommand{\cM}{\mathcal{M}}

\begin{document}

\title{H\"{o}lder classifications of finite-dimensional linear flows}
\author{Arno Berger and Anthony Wynne\\Mathematical and Statistical Sciences\\University of Alberta\\Edmonton, Alberta, {\sc Canada}}

\date{\today}
\maketitle

\begin{abstract}
\noindent
Two flows on a finite-dimensional normed space $X$ are equivalent if
some homeomorphism $h$ of $X$ preserves all orbits, i.e., $h$ maps
each orbit onto an orbit. Under the assumption that $h$,
$h^{-1}$ both are $\beta$-H\"{o}lder continuous near the origin for
some (or all) $0<\beta< 1$, a complete classification with respect to
some-H\"{o}lder (or all-H\"{o}lder) equivalence is established for
linear flows on $X$, in terms of basic linear algebra properties of 
their generators. Consistently utilizing equivalence
instead of the more restrictive conjugacy, the classification theorems
extend and unify known results. Though entirely elementary, the analysis is
somewhat intricate and highlights, more clearly than does the existing
literature, the fundamental roles played by linearity and the
finite-dimensionality of $X$. 
\end{abstract}
\hspace*{6.6mm}{\small {\bf Keywords.} Equivalence between flows,
  linear flow, H\"{o}lder equivalence, Lyapunov similarity.}

\noindent
\hspace*{6.6mm}{\small {\bf MSC2020.} 34A30, 34C41, 34D08, 37C15.}

%%%%%%%%%%%%%%%%%%%%%%%%%%%%%%%%%%%%%%%%%%%%%%%%%%%%%%%%%%%%%%%

\section{Introduction}\label{sec1}

%%%%%%%%%%%%%%%%%%%%%%%%%%%%%%%%%%%%%%%%%%%%%%%%%%%%%%%%%%%%%%%

Let $X\ne \{0\}$ be a finite-dimensional normed space over $\R$ and $\varphi$
a flow on $X$, i.e., $\varphi : \R \times X \to X$ is
continuous with $\varphi(t+s , x) = \varphi \bigl( t, \varphi (s,x)
\bigr)$ and $\varphi(0,x) = x$ for all $t,s\in \R$, $x\in X$. A
fundamental question throughout dynamics is that of classification:
When, precisely, are two flows $\varphi$, $\psi$ on $X$ {\em the
  same}, and in what sense? A geometrically motivated 
approach to this question is as follows: Say that $\varphi$, $\psi$ are {\bf
  equivalent}, in symbols $\varphi \thicksim \psi$, if there exists a
homeomorphism $h:X\to X$ that maps every $\varphi$-orbit onto
a $\psi$-orbit, i.e., 
\begin{equation}\label{eq01}
  h \bigl(
  \bigl\{
\varphi (t,x) : t \in \R
  \bigr\}
  \bigr) = \bigl\{ \psi \bigl( t,h(x)\bigr) :t\in \R \bigr\} \qquad
  \forall x \in X \, .
\end{equation}
Imposing additional regularity requirements on $h$ naturally yields
further, narrower forms of equivalence. Specifically, if $h$, $h^{-1}$
both are $\beta$-H\"{older} continuous for
some $0<\beta < 1$ (or all $0<\beta < 1$, or $\beta = 1$) then
$\varphi$, $\psi$ are {\bf some-H\"{o}lder} (or {\bf all-H\"{o}lder},
or {\bf Lipschitz}) {\bf equivalent}, in symbols $\varphi
\stackrel{0^+}{\thicksim} \psi$ (or $\varphi
\stackrel{1^-}{\thicksim} \psi$, or $\varphi
\stackrel{1}{\thicksim} \psi$). More restrictively still, if $h$,
$h^{-1}$ both are differentiable (or linear) then $\varphi$, $\psi$ are
{\bf differentiably} (or {\bf linearly}) {\bf equivalent}, in symbols $\varphi
\stackrel{{\sf diff}}{\thicksim} \psi$ (or $\varphi
\stackrel{{\sf lin}}{\thicksim} \psi$). As discussed in detail in Section
\ref{sec2} below, these equivalences constitute but six familiar
``vertices'' in an infinite ``graph'' of equivalences, no two of
which coincide entirely, that is, for all pairs of flows on $X$. This
in turn leads to an infinitude of natural, genuinely different
classifications of flows (see Figure \ref{fig1} below). 

Building on the classical literature briefly reviewed below, the
present article, together with \cite{BW, BW3}, completely answers the
question of equivalence for {\em linear\/} flows. As it turns out, for
such flows all (infinitely many) equivalences coalesce, rather
amazingly, into a mere {\em four\/} different forms, informally
referred to, respectively, as {\bf topological}, {\bf H\"{o}lder}, {\bf Lipschitz},
and {\bf smooth equivalence} (see Figures \ref{fig0a}
and \ref{fig2} below).
Recall that a flow $\varphi$ on $X$ is {\bf linear} if the time-$t$
map $\varphi_t = \varphi (t,\cdot):X\to X$ is linear, or equivalently
if $\varphi_t = e^{t A^{\varphi}}$,
for every $t\in \R$, with a (unique) linear operator $A^{\varphi}$ on
$X$ called
the {\bf generator} of $\varphi$. Henceforth, upper case Greek letters $\Phi$, $\Psi$ are used
exclusively to denote linear flows. All four equivalences between linear
flows $\Phi$, $\Psi$ just alluded to are fully characterized below, in terms of basic
linear algebra properties of $A^{\Phi}$, $A^{\Psi}$. This yields four
classification theorems, each of which in one way or another
extends, complements, or unifies earlier results in the literature.

The first main result of this article, then, is the following {\bf
  topological classification theorem} which also
shows that, perhaps surprisingly, equivalence between linear flows always
entails some-H\"{o}lder equivalence. To state the result,
recall that every linear flow $\Phi$ on $X$ determines a unique
$\Phi$-invariant decomposition $X= X_{\sf S}^{\Phi} \oplus X_{\sf
  C}^{\Phi} \oplus X_{\sf U}^{\Phi}$ into stable, central, and
unstable subspaces, along with a unique decomposition $\Phi
\stackrel{{\sf lin}}{\cong} \Phi_{\sf S} \times \Phi_{\sf C} \times
\Phi_{\sf U}$; see Sections \ref{sec2} and \ref{sec2a} below for formal details. 

\begin{theorem}\label{thma}
Let $\Phi$, $\Psi$ be linear flows on $X$. Then each of the following
three statements implies the other two:
\begin{enumerate}
\item $\Phi \stackrel{0^+}{\thicksim}\Psi$, i.e., $\Phi$, $\Psi$ are some-H\"{o}lder equivalent;
\item $\Phi \thicksim \Psi$, i.e., $\Phi$, $\Psi$ are equivalent;
  \item $\{\dim X_{\sf S}^{\Phi}, \dim X_{\sf U}^{\Phi}\} = \{\dim
    X_{\sf S}^{\Psi}, \dim X_{\sf U}^{\Psi}\}$, and there exists an $\alpha
    \in \R \setminus \{0\}$ so that $A^{\Phi_{\sf
        C}}$, $\alpha A^{\Psi_{\sf C}}$ are similar.
\end{enumerate}
\end{theorem}

An important insight implicit in Theorem \ref{thma} is that the
validity of (\ref{eq01}) for linear $\varphi$, $\psi$ guarantees not
only that $h$, $h^{-1}$ are (or can be chosen to be) $\beta$-H\"{o}lder
continuous for some $\beta > 0$, but also that, with an appropriate $\alpha \in
\R\setminus \{0\}$,
\begin{equation}\label{eq02}
h \bigl( \varphi (t, x)\bigr) = \psi \bigl( \alpha t , h(x)\bigr)
\qquad \forall t\in \R , x\in X \, .
\end{equation}
Notice how (\ref{eq02}) in general is much more restrictive than
(\ref{eq01}). Virtually all studies on equivalences between (linear) flows in the literature are based on
(\ref{eq02}), often with the additional requirement that $\alpha > 0$,
or indeed $\alpha = 1$. By contrast, the natural, significantly more general
form (\ref{eq01}) is referred to only perfunctorily, if at all
\cite{Kuiper, Ladis, MM, Willems}; see also Section \ref{sec2} and the discussion in
\cite[Sec.\ 5]{BW}.

The second main result of this article is a {\bf H\"{o}lder
  classification theorem} involving the concept of Lyapunov
similarity, introduced rigorously in Section
\ref{sec2a}. For now, simply say that two linear operators are {\bf
  Lyapunov similar} if they (more precisely, the flows they generate) have the same
Lyapunov exponents, with matching multiplicities.
  
\begin{theorem}\label{thmb}
Let $\Phi$, $\Psi$ be linear flows on $X$. Then each of the following statements implies the other:
\begin{enumerate}
\item $\Phi \stackrel{1^-}{\thicksim} \Psi$, i.e., $\Phi$, $\Psi$ are all-H\"{o}lder equivalent;
\item there exists an $\alpha
    \in \R \setminus \{0\}$ so that $A^{\Phi}$, $\alpha A^{\Psi}$ are
    Lyapunov similar and $A^{\Phi_{\sf  C}}$, $\alpha A^{\Psi_{\sf C}}$ are similar.
\end{enumerate}
\end{theorem}

Variants of (ii)$\Leftrightarrow$(iii) in Theorem \ref{thma} utilizing
(\ref{eq02}) were first proved in
\cite{Kuiper, Ladis}, though for {\em hyperbolic\/} flows, i.e., for $X_{\sf
  C}^{\Phi} = X_{\sf C}^{\Psi} = \{0\}$, the result is much older;
see, e.g., \cite{Amann, Irwin, R} for broad context, as well as
\cite{ACK1, AK, DSS, He, LZ, Willems} and references therein for
specific subsequent studies. As far as the authors have been able
to ascertain, neither the full
strength of Theorem \ref{thma} utilizing only (\ref{eq01}) nor
the fact that (i)$\Leftrightarrow$(ii) have yet
been documented in the literature. Similarly, a weaker variant
of Theorem \ref{thmb} may be gleaned from the examples in \cite{MM},
albeit with considerable hand-waving, but again its full strength and proof
appear to be new. Given the simple, definitive nature of
Theorems \ref{thma} and \ref{thmb}, as well as the importance of
linear differential equations throughout science (education), the
present article aims to provide elementary, self-contained proofs of
both results which, together with \cite{BW, BW3}, hopefully will
inform future applications and pedagogy.

To put the results in context, it is instructive to compare them to their Lipschitz and
smooth counterparts; stated here without proof, these
have been proved by the authors elsewhere \cite{BW, BW3}.
Though structurally analogous to Theorem \ref{thmb}, the following {\bf Lipschitz
classification theorem} significantly differs from its H\"{o}lder
counterpart, due to the discrepancy between Lipschitz and
Lyapunov similarities. Motivated by precursors in \cite{KS, MM},
Lipschitz similarity is introduced and discussed in detail in \cite{BW3}. For the purpose
of the present article, it suffices to note that Lipschitz similarity
of two linear operators on $X$ requires (most of) their eigenvalues
and multiplicities to match, whereas Lyapunov similarity only requires
the matching of real parts of eigenvalues (and cumulative
multiplicities). Correspondingly, two linear operators are Lyapunov
similar whenever Lipschitz similar, and they are Lipschitz similar whenever
similar, but neither implication is reversible for $\dim X \ge 2$.

\begin{prop}\label{thmc}
Let $\Phi$, $\Psi$ be linear flows on $X$. Then each of the following
statements implies the other:
\begin{enumerate}
\item $\Phi \stackrel{1}{\thicksim} \Psi$, i.e., $\Phi$, $\Psi$ are Lipschitz equivalent;
\item there exists an $\alpha
    \in \R \setminus \{0\}$ so that $A^{\Phi}$, $\alpha A^{\Psi}$ are
    Lipschitz similar and $A^{\Phi_{\sf  C}}$, $\alpha A^{\Psi_{\sf C}}$ are similar.
\end{enumerate}
\end{prop}

In essence the following {\bf smooth classification theorem} has
been established in \cite[Thm.1.2]{BW}, with weaker versions found in
many textbooks \cite{Amann, CK, R}. Although \cite{BW}
employs a more restrictive notion of equivalence than the present article, the result is readily seen
to carry over verbatim.

\begin{prop}\label{thmd}
Let $\Phi$, $\Psi$ be linear flows on $X$. Then each of the following
three statements implies the other two:
\begin{enumerate}
\item $\Phi \stackrel{{\sf lin}}{\thicksim}\Psi$, i.e., $\Phi$, $\Psi$ are linearly equivalent;
\item $\Phi \stackrel{{\sf diff}}{\thicksim}\Psi$, i.e., $\Phi$, $\Psi$ are differentiably equivalent;
  \item there exists an $\alpha    \in \R \setminus \{0\}$ so that  $A^{\Phi}$, $\alpha A^{\Psi}$ are similar.
\end{enumerate}
\end{prop}

A striking consequence of Theorem \ref{thmb} as well as
Propositions \ref{thmc} and \ref{thmd} is that, in analogy to Theorem
\ref{thma}, assuming $\Phi \stackrel{\bigstar}{\thicksim} \Psi $ with
$\bigstar \in \{1^-, 1, {\sf diff}, {\sf lin}\}$ guarantees that the
(all-H\"{older}, Lipschitz, differentiable, or linear) homeomorphism $h$ can be chosen so as to
satisfy (\ref{eq02}). In other words, for linear $\varphi$, $\psi$,
and for every degree of regularity of $h$ considered herein,
(\ref{eq01}) always entails (\ref{eq02}). This remarkable property is
indicative of the extraordinary coherence between individual orbits of
linear flows. It does not appear to be shared by any wider class of
flows on $X$.  

To illustrate the four theorems above, first
notice that for $\dim X = 1$ trivially all classifications
coincide: Every linear flow on $X=\R^1$ is (smoothly, Lipschitz,
H\"{o}lder, or topologically) equivalent to the flow generated by precisely
one of $[0]$ and $[1]$. Already for $\dim X = 2$,
however, the discrepancies between the four classifications become apparent:
Every linear flow on $X=\R^2$ is {\em smoothly\/} equivalent to the flow
generated by precisely one of either
\begin{equation}\label{eq0p1}
  \left[
    \begin{array}{cc}
0 & 0 \\ 0 & 0 
    \end{array}
  \right] ,
\left[
    \begin{array}{cc}
0 & 1 \\ 0 & 0 
    \end{array}
  \right] ,
  \left[
    \begin{array}{cr}
0 & -1 \\ 1 & 0 
    \end{array}
  \right] ,
\end{equation}
or a (necessarily unique) matrix from
$$
\left[
    \begin{array}{cc}
a & 0 \\ 0 & 1 
    \end{array}
  \right] ,
  \left[
    \begin{array}{cc}
1 & 1 \\ 0 & 1 
    \end{array}
  \right] ,
  \left[
    \begin{array}{cr}
1 & -b \\ b & 1 
    \end{array}
  \right] \qquad a\in [-1,1], b\in \R^+; 
$$
it is {\em Lipschitz\/} equivalent to the flow
generated by precisely one of either (\ref{eq0p1}) or
$$
\left[
    \begin{array}{cc}
a & 0 \\ 0 & 1 
    \end{array}
  \right] ,
  \left[
    \begin{array}{cc}
1 & 1 \\ 0 & 1 
    \end{array}
  \right] \qquad a\in [-1,1];
  $$
it is {\em H\"{o}lder\/} equivalent to the flow
generated by precisely one of either (\ref{eq0p1}) or
$$
\left[
    \begin{array}{cc}
a & 0 \\ 0 & 1 
    \end{array}
  \right] \qquad a\in [-1,1] ;
$$
and it is {\em topologically\/} equivalent  to the flow
generated by precisely one of either (\ref{eq0p1}) or
$$
\left[
    \begin{array}{rc}
-1 & 0 \\ 0 & 1
    \end{array}
  \right] ,
  \left[
    \begin{array}{cc}
0 & 0 \\ 0 & 1
    \end{array}
  \right] ,
  \left[
    \begin{array}{cc}
1 & 0 \\ 0 & 1
    \end{array}
  \right] ;
$$
see also Figures \ref{fig0a} and \ref{fig0}.

\begin{figure}[ht] 
  \psfrag{tl1a}[]{\small $O_2$}
  \psfrag{tl1b}[]{\small $ J_2$}
  \psfrag{tl1c}[]{\small $ J_1(i)$}
  \psfrag{tl1e}[]{\small $ \mbox{\rm diag}\, [a,1]$ with $a\in [-1,1]$}
  \psfrag{tl1d}[]{\small $J_2(1)$}
  \psfrag{tl1f}[]{\small $J_1(1+ib)$ with $b\in \R^+$}
  \psfrag{tl4a}[]{\small $ \mbox{\rm diag}\, [-1,1]$}
  \psfrag{tl4b}[]{\small $ \mbox{\rm diag}\, [0,1]$}
  \psfrag{tl4c}[]{\small $I_2$}
  \psfrag{tsmoo}[r]{{\bf smooth}}
  \psfrag{tlip}[r]{{\bf Lipschitz}}
  \psfrag{thoel}[r]{{\bf H\"{o}lder}}
  \psfrag{ttop}[r]{{\bf topological}}
%
% scale 0.85
%
\begin{center}
\includegraphics{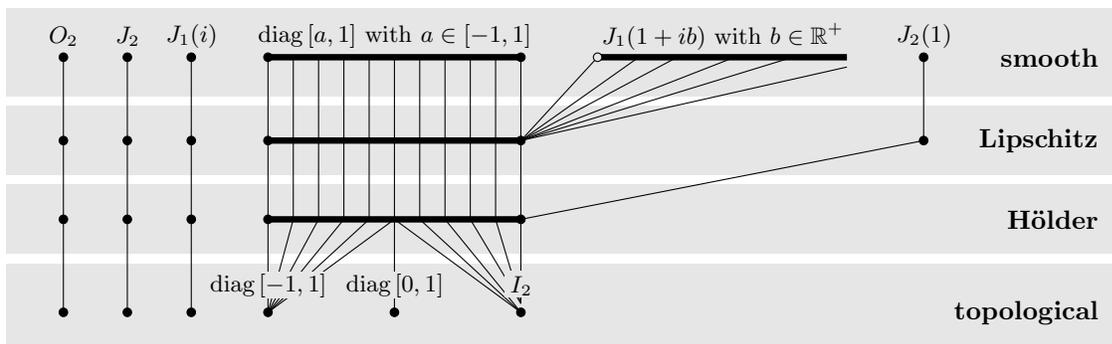}
\end{center}
\vspace*{-4mm}
\caption{No two of the four classifications of all linear flows on $X=\R^2$
  coincide.}\label{fig0a}
\end{figure}

\begin{figure}[t] 
  \psfrag{tl1}[]{$A^{\Phi} = O_2$}
  \psfrag{tl2}[]{$A^{\Phi} = J_2$}
  \psfrag{tl3}[]{$A^{\Phi} = J_1(i)$}
  \psfrag{tl4}[]{$A^{\Phi} = \mbox{\rm diag}\, [a, 1]$}
  \psfrag{tl1a}[]{$\mbox{\rm Fix} \, \Phi = \R^2$}
  \psfrag{tl4a}[]{$a=-1$}
  \psfrag{tl4b}[]{$-1< a<0$}
  \psfrag{tl4c}[]{$a=0$}
  \psfrag{tl4d}[]{$0<a<1$}
  \psfrag{tl4e}[]{$a=1$}
%
% scale 0.85
%
\begin{center}
\includegraphics{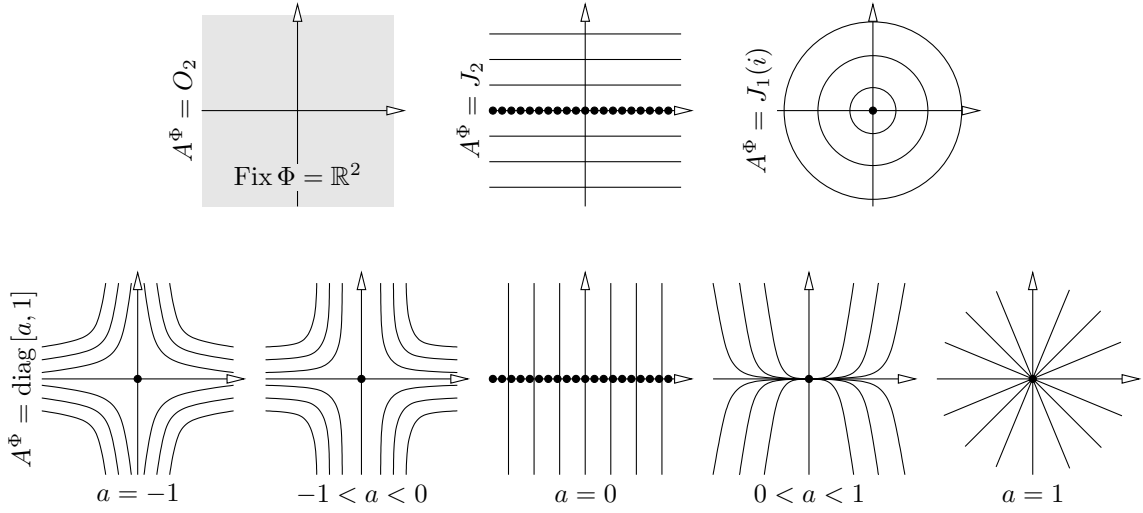}
\end{center}
\vspace*{-2mm}
\caption{Displaying all possible phase portraits (without orientation)
  of a linear flow $\Phi$ on $X=\R^2$, up to H\"{o}lder equivalence
  (Theorem \ref{thmb}). In the bottom half, the two left-most flows are
  (topologically) equivalent, and so are the two right-most flows
  (Theorem \ref{thma}); see also the lower half of Figure \ref{fig0a}.}\label{fig0}
\end{figure}

\medskip

The remainder of this article is organized as follows: Section
\ref{sec2} properly introduces the fundamental notion of equivalence
between flows on $X$, motivated by (\ref{eq01}), as well as natural refinements
thereof. Section \ref{sec2a} briefly reviews a few basic concepts
specific to linear flows, notably irreducibility and Lyapunov
exponents. Sections \ref{sec3} and \ref{sec4} carry out detailed 
analyses of $\beta$-H\"{o}lder relations between linear flows
($0<\beta < 1$) and the behaviour of minimal periods under such
relations, respectively. The observations in both sections are of an
auxiliary nature but may also be of independent interest. Section
\ref{sec5} presents the proof of the main results, Theorems
\ref{thma} and \ref{thmb}, in mildly extended form. A brief concluding
Section \ref{sec7} clarifies how the main results naturally carry over
to {\em complex\/} spaces.

Throughout, the familiar symbols $\N$, $\N_0$, $\Q^+$, $\Q$, $\R^+$,
$\R$, and $\C$ denote the sets of all positive whole, non-negative
whole, positive rational, rational, positive real, real, and
complex numbers respectively, each with their usual arithmetic, order, and
topology. Every $z\in \C$ can be written uniquely as $z=a+ib$ where
$a= {\sf Re}z$, $b={\sf Im}z$ are real numbers, with complex conjugate
$\overline{z} = a-ib$ and modulus $|z|=\sqrt{a^2+b^2}$. Given any
$v,w\in \C$ and $Z\subset \C$, let $v+wZ = \{v+wz: z \in Z\}$.

%%%%%%%%%%%%%%%%%%%%%%%%%%%%%%%%%%%%%%%%%%%%%%%%%%%%%%%%%%%%%%%%

\section{Equivalences between flows}\label{sec2}

%%%%%%%%%%%%%%%%%%%%%%%%%%%%%%%%%%%%%%%%%%%%%%%%%%%%%%%%%%%%%%%%

Throughout, let $X = \R^d$, where the actual value of $d\in \N$ is either clear from the
context or irrelevant. Endow $X$ with the Euclidean norm
$|\cdot|$; this is solely for convenience, as all concepts and
results herein are readily seen to be independent of any particular
norm. Denote by $e_1, \ldots , e_d$ the canonical basis of $X$,
by $O_X=O_d$, $I_X=I_d$ the zero and identity operator (or $d\times
d$-matrix) respectively, and let
$B_r(x) = \{y\in X : |y-x|<r\}$ for every $r\in \R^+$, $x\in X$.
In accordance with a familiar tenet of linear analysis
\cite{Jam}, the case of a (finite-dimensional) normed space over $\C$ does not pose any
additional challenge; it is only considered briefly in Section \ref{sec7} below.

Given a flow $\varphi$ on $X$, the $\varphi$-{\bf orbit} of any
$x\in X$ is $\varphi_{\R}( x) := \{\varphi_t (x) : t\in \R\}$. For any
two flows $\varphi$, $\psi$ on $X$ and any homeomorphism $h:X\to X$, say
that $\varphi$ is $h$-{\bf related} to $\psi$, in symbols $\varphi
\stackrel{h}{\thicksim} \psi$, if (\ref{eq01}) holds, that is, if
$$
  h \bigl(
\varphi_{\R}( x)
\bigr) = \psi_{\R}\bigl(  h(x)
\bigr) \qquad \forall x \in X \, ,
$$
or equivalently if $h$, $h^{-1}$ both map orbits into orbits. An
orbit-wise characterization of $\varphi
\stackrel{h}{\thicksim} \psi$ is readily established.

\begin{prop}\label{prop12}
Let $\varphi$, $\psi$ be flows on $X$. For every homeomorphism $h:X\to
X$ the following are equivalent:
\begin{enumerate}
\item $\varphi \stackrel{h}{\thicksim} \psi$;
\item for every $x\in X$ there exists a continuous bijection
  $\tau_x:\R\to \R$ with $\tau_x (0)=0$ so that
  $$h \bigl( \varphi_t (x)\bigr) = \psi_{\tau_x(t)} \bigl( h(x) \bigr)
  \qquad \forall t \in \R \, .
  $$ 
\end{enumerate}
\end{prop}

In light of Proposition \ref{prop12}, the simplest, most
fundamental equivalence between flows, previewed in the Introduction, is as follows: Say
that $\varphi$, $\psi$ are {\bf equivalent}, in symbols $\varphi
\thicksim \psi$, if $\varphi \stackrel{h}{\thicksim} \psi$ for some
homeomorphism $h$. Clearly, this defines an equivalence relation on the class of all
flows on $X$. Informally put, $\varphi \thicksim \psi$ means that every
$\varphi$-orbit is, up to a change of spatial coordinates (via $h$)
and a (possibly orbit-dependent) re-parametrization of time (via
$\tau_x$), also a $\psi$-orbit and vice versa.

Observe that $\tau_x$ in Proposition \ref{prop12} is uniquely
determined unless $\varphi_{\R}(x) = \{x\}$, i.e., unless $x$ is a
fixed point of $\varphi$, in symbols $x\in \mbox{\rm Fix}\, \varphi$;
in the latter case the continuous bijection $\tau_x$ is arbitrary.
Imposing additional requirements on the family $\tau = (\tau_x)_{x\in
  X}$ naturally yields other, narrower
equivalences. For instance, say that $\varphi$, $\psi$ are {\bf
  strictly equivalent}, in symbols $\varphi
\thickapprox \psi$, if $\varphi \stackrel{h}{\thicksim} \psi$ for some
$h$ so that either $\tau_x$ is increasing for every $x\in X\setminus
\mbox{\rm Fix} \, \varphi$ or else $\tau_x$ is
decreasing for every $x$. A more stringent condition is that
$\tau_x$ be independent of $x$ altogether. In this case, it is readily
seen that, with some $\alpha \in \R
\setminus \{0\}$, simply $\tau_x (t) = \alpha t$ for all $x\in
X\setminus \mbox{\rm Fix}\, \varphi$,
$t\in \R$; in other words, (\ref{eq02}) holds. The latter situation
henceforth is denoted
$\varphi \simeq \psi$; in case $\alpha>0$ it is referred to in
\cite{BW} as $\varphi$, $\psi$ being {\bf flow equivalent}. In summary,
\begin{equation}\label{eq12}
\varphi \simeq \psi \quad \Longrightarrow \quad  \varphi  \thickapprox
\psi \quad \Longrightarrow \quad  \varphi  \thicksim \psi \, ,
\end{equation}
and simple examples show that the left and right
implication in (\ref{eq12}) cannot be reversed in general for $d\ge 2$
and for any $d\in \N$, respectively. (For $d=1$ trivially $ \varphi
\thickapprox  \psi$ implies $\varphi \simeq \psi $.)

Many other equivalences between flows are conceivable beyond
the three forms appearing in (\ref{eq12}). To see 
but one example, define $\varphi \bowtie \psi$ to mean that $\varphi
\stackrel{h}{\thicksim} \psi$ for some $h$ so that $\lim_{|t|\to
  \infty} \tau_x(t)/t$ exists and is nonzero for every $x\in
X\setminus \mbox{\rm Fix} \, \varphi$. Again, this defines a
bona fide equivalence relation, with
$$
\varphi \simeq \psi \quad \Longrightarrow \quad
\varphi \bowtie \psi \quad \Longrightarrow
\quad  \varphi  \thicksim \psi  \, ,
$$
and again neither of these implications can be reversed in general for
$d\ge 2$. Examples like this suggest that $\thicksim$
is the most general equivalence, whereas $ \simeq$ is the most
restrictive, and $\bowtie $, $\thickapprox$ are somehow
intermediate between these two. With the additional requirement that $\alpha = 1$, and thus 
simply $h\bigl( \varphi_t (x)\bigr) = \psi_t \bigl( h(x) \bigr)$ for
some $h$ and all $t$, $x$, the relation $\simeq$ has often been
employed (sometimes implicitly or with different notation) in the literature, with $\varphi$,
$\psi$ referred to as being ({\bf topologically}) {\bf conjugate},
here in symbols $\varphi \cong
\psi$; see, e.g., \cite{ACK1, AK, DSS, He, KS, Kuiper, LZ, Willems}.

\begin{figure}[ht] 
\psfrag{t1lin}[]{$\varphi \stackrel{{\sf lin}}{\simeq} \psi$}
\psfrag{t2lin}[]{$\varphi \stackrel{{\sf lin}}{\thickapprox} \psi$}
\psfrag{t3lin}[]{$\varphi \stackrel{{\sf lin}}{\thicksim} \psi$}
\psfrag{t1diff}[]{$\varphi \stackrel{{\sf diff}}{\simeq} \psi$}
\psfrag{t2diff}[]{$\varphi \stackrel{{\sf diff}}{\thickapprox} \psi$}
\psfrag{t3diff}[]{$\varphi \stackrel{{\sf diff}}{\thicksim} \psi$}
\psfrag{t1lip}[]{$\varphi \stackrel{1}{\simeq} \psi$}
\psfrag{t2lip}[]{$\varphi \stackrel{1}{\thickapprox} \psi$}
\psfrag{t3lip}[]{$\varphi \stackrel{1}{\thicksim} \psi$}
\psfrag{t1all}[]{$\varphi \stackrel{1^-}{\simeq} \psi$}
\psfrag{t2all}[]{$\varphi \stackrel{1^-}{\thickapprox} \psi$}
\psfrag{t3all}[]{$\varphi \stackrel{1^-}{\thicksim} \psi$}
\psfrag{t1some}[]{$\varphi \stackrel{0^+}{\simeq} \psi$}
\psfrag{t2some}[]{$\varphi \stackrel{0^+}{\thickapprox} \psi$}
\psfrag{t3some}[]{$\varphi \stackrel{0^+}{\thicksim} \psi$}
\psfrag{t10}[]{$\varphi \stackrel{0}{\simeq} \psi$}
\psfrag{t20}[]{$\varphi \stackrel{0}{\thickapprox} \psi$}
\psfrag{t30}[]{$\varphi \stackrel{0}{\thicksim} \psi$}
\psfrag{tif}[]{\small ($d=1$)}
%
% scale 0.75
%
%
%\vspace*{-6mm}
%
\begin{center}
\includegraphics{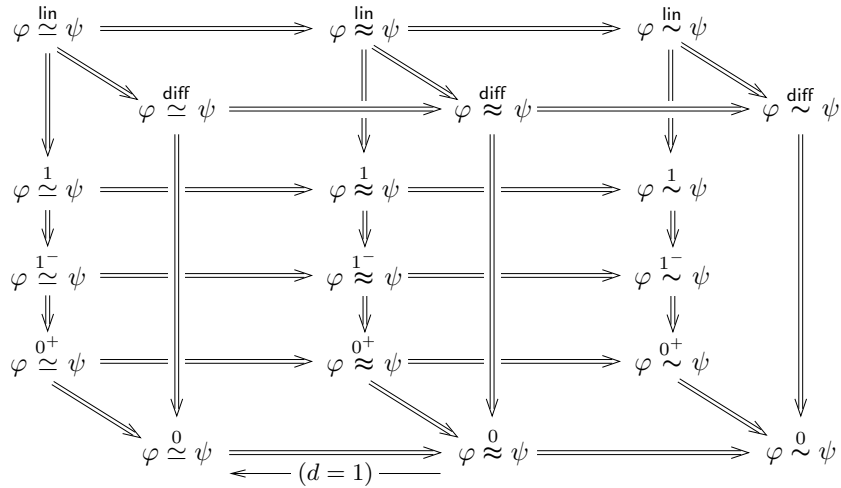}
\end{center}
\vspace*{-4mm}
\caption{Relating {\em eighteen\/} natural equivalences between flows
  $\varphi$, $\psi$ on $X=\R^d$, $d\in \N$. All equivalences are
  genuinely different in that no conceivable implication
not shown in the diagram is valid in general.}\label{fig1}
\end{figure}

Apart from imposing additional requirements on $\tau$, an
important, natural way of refining $ \varphi  \stackrel{h}{\thicksim}
\psi  $, alluded to in the Introduction, is
to require additional regularity of $h$. Note that if $ \varphi
\stackrel{h}{\thicksim} \psi  $ then also $ \varphi
\stackrel{\overline{h}}{\thicksim} \overline{\psi} $, where
$\overline{h} = h - h(0)$ and $\overline{\psi}_t = \psi_t \bigl( \cdot
+ h(0)\bigr) - h(0)$ for all $t\in \R$. Thus, no generality is lost by
assuming that $h(0)=0$. Bearing this in mind, denote by $\cH= \cH (X)$ the set of all
homeomorphisms $h:X\to X$ with $h(0)=0$, and let $\cH_{\beta}=
\cH_{\beta} (X)$ with $0\le \beta \le 1$ be the set of all $h\in \cH$ for which $h$,
$h^{-1}$ both satisfy a $\beta$-H\"{o}lder condition (a.k.a.\
Lipschitz condition in case $\beta=1$) near $0$, i.e.,
$$
\cH_{\beta} = \left\{ h \in \cH : \exists r \in \R^+ \: \mbox{\rm
    s.t.} \: \sup\nolimits_{x,y\in B_r(0), x\ne y}  \frac{ |h(x) -
    h(y)| +  |h^{-1}(x) - h^{-1}(y)| }{|x-y|^{\beta}} < \infty \right\}  \, ;
$$
see, e.g., \cite{F, Hei} for comprehensive accounts on H\"{o}lder and
Lipschitz analysis. Since $\beta \mapsto \cH_{\beta}$ is decreasing, one may
also consider
$$
\cH_{\beta^-}:= \bigcap_{\gamma < \beta} \cH_{\gamma} \quad (\mbox{\rm
  if } \beta > 0) \, , \qquad
\cH_{\beta^+}:= \bigcup_{\gamma > \beta} \cH_{\gamma} \quad (\mbox{\rm
  if } \beta < 1) \, .
$$
Furthermore, let
$$
\cH_{\sf diff} = \bigl\{ h\in \cH : h, h^{-1} \enspace \mbox{\rm are
  differentiable at } 0\, \bigr\} \, , \qquad
\cH_{\sf lin} = \bigl\{h\in \cH : h \enspace \mbox{\rm is
 linear}\,   \bigr\} \, .
$$
This yields a strictly decreasing family of subsets of $\cH_0=\cH$, 
$$
\cH_0 \supset \cH_{0^+} \supset \ldots \supset \cH_{\beta^-} \supset
\cH_{\beta} \supset \cH_{\beta^+} \supset \ldots \supset
\cH_{1^-}\supset \cH_1 \supset \cH_{\sf lin} \qquad \forall 0<\beta < 1 \, ,
$$
and clearly also $\cH_0 \supset \cH_{\sf diff} \supset \cH_{\sf lin}$,
whereas $\cH_{0^+}\not \supset \cH_{\sf diff}$ and $\cH_{\sf diff}\not
\supset \cH_{1}$. Correspondingly, given any $\bigstar \in\{ 0, 0^+, \beta^- ,
\beta , \beta^+, 1^- ,1,  {\sf diff}, {\sf lin}\}$ with $0<\beta < 1$, understand $\varphi
\stackrel{\bigstar }{\thicksim} \psi$ to mean that $\varphi
\stackrel{h}{\thicksim} \psi$ for some $h\in \cH_{\bigstar}$. Though reflexive
and symmetric by definition,  the relation $\stackrel{\bigstar}{\thicksim}$
is not transitive, and hence not an equivalence relation if $\bigstar\in 
\{\beta^- , \beta , \beta^+\}$ and $0<\beta < 1$. Only in the six
cases $\bigstar\in \{0,0^+, 1^- , 1, {\sf diff}, {\sf lin}\}$, therefore, does
$\stackrel{\bigstar}{\thicksim}$ lead to a classification. Say that
$\varphi$, $\psi$ are {\bf topologically}, 
{\bf some-H\"{o}lder}, {\bf all-H\"{older}}, {\bf Lipschitz}, {\bf
  differentiably}, and {\bf linearly equivalent} if $\varphi
\stackrel{0}{\thicksim} \psi$, $\varphi
\stackrel{0^+}{\thicksim} \psi$, $\varphi
\stackrel{1^-}{\thicksim} \psi$, $\varphi
\stackrel{1}{\thicksim} \psi$, $\varphi
\stackrel{{\sf diff}}{\thicksim} \psi$, and
$\varphi \stackrel{{\sf lin}}{\thicksim} \psi$ respectively. Clearly,
\begin{equation}\label{eq13}
\varphi \stackrel{{\sf lin}}{\thicksim} \psi \quad \Longrightarrow
\quad
\varphi
\stackrel{1}{\thicksim} \psi
\quad \Longrightarrow
\quad \varphi
\stackrel{1^-}{\thicksim} \psi
 \quad \Longrightarrow
\quad \varphi
\stackrel{0^+}{\thicksim} \psi
\quad \Longrightarrow
\quad \varphi
\stackrel{0}{\thicksim} \psi \, ,
\end{equation}
as well as 
\begin{equation}\label{eq13a}
\varphi \stackrel{{\sf lin}}{\thicksim} \psi \quad \Longrightarrow
\quad
\varphi
\stackrel{{\sf diff}}{\thicksim} \psi \quad \Longrightarrow
\quad \varphi
\stackrel{0}{\thicksim} \psi \, ,
\end{equation}
and simple examples again show that none of the implications in
(\ref{eq13}), (\ref{eq13a}) can be reversed in general, not even for
$d=1$. Also, $\varphi \stackrel{{\sf diff}}{\thicksim} \psi \not
\Rightarrow \varphi \stackrel{0^+}{\thicksim} \psi$ and
$\varphi \stackrel{1}{\thicksim} \psi\not \Rightarrow \varphi
\stackrel{{\sf diff}}{\thicksim} \psi $ in general. 
In a similar vein, one may consider the equivalence relations
$\stackrel{\bigstar}{\simeq}$ and $\stackrel{\bigstar}{\thickapprox}$
for any $\bigstar \in \{0,0^+, 1^- , 1, {\sf diff}, {\sf
  lin}\}$. Altogether, then, there are (at 
least) eighteen different natural equivalences between flows on
$X$, leading in turn to an equal number of different classifications;
see Figure \ref{fig1}. 

%%%%%%%%%%%%%%%%%%%%%%%%%%%%%%%%%%%%%%%%%%%%%%%%%%%%%%%%%%%%%%%%

\section{Linear flow preliminaries}\label{sec2a}

%%%%%%%%%%%%%%%%%%%%%%%%%%%%%%%%%%%%%%%%%%%%%%%%%%%%%%%%%%%%%%%%

\begin{figure}[ht]
\psfrag{t1lin}[]{$\Phi \stackrel{{\sf lin}}{\simeq} \Psi$}
\psfrag{t2lin}[]{$\Phi \stackrel{{\sf lin}}{\thickapprox} \Psi$}
\psfrag{t3lin}[]{$\Phi \stackrel{{\sf lin}}{\thicksim} \Psi$}
\psfrag{t1diff}[]{$\Phi \stackrel{{\sf diff}}{\simeq} \Psi$}
\psfrag{t2diff}[]{$\Phi \stackrel{{\sf diff}}{\thickapprox} \Psi$}
\psfrag{t3diff}[]{$\Phi \stackrel{{\sf diff}}{\thicksim} \Psi$}
\psfrag{t1lip}[]{$\Phi \stackrel{1}{\simeq} \Psi$}
\psfrag{t2lip}[]{$\Phi \stackrel{1}{\thickapprox} \Psi$}
\psfrag{t3lip}[]{$\Phi \stackrel{1}{\thicksim} \Psi$}
\psfrag{t1all}[]{$\Phi \stackrel{1^-}{\simeq} \Psi$}
\psfrag{t2all}[]{$\Phi \stackrel{1^-}{\thickapprox} \Psi$}
\psfrag{t3all}[]{$\Phi \stackrel{1^-}{\thicksim} \Psi$}
\psfrag{t1some}[]{$\Phi \stackrel{0^+}{\simeq} \Psi$}
\psfrag{t2some}[]{$\Phi \stackrel{0^+}{\thickapprox} \Psi$}
\psfrag{t3some}[]{$\Phi \stackrel{0^+}{\thicksim} \Psi$}
\psfrag{t10}[]{$\Phi \stackrel{0}{\simeq} \Psi$}
\psfrag{t20}[]{$\Phi \stackrel{0}{\thickapprox} \Psi$}
\psfrag{t30}[]{$\Phi \stackrel{0}{\thicksim} \Psi$}
\psfrag{tsmoo}[r]{{\bf smooth}}
\psfrag{tlip}[r]{{\bf Lipschitz}}
\psfrag{thol}[r]{{\bf H\"{o}lder}}
\psfrag{ttop}[r]{{\bf topological}}
\psfrag{tif}[]{\small ($X_{\sf C}=X$ or $d=1$)}
\psfrag{tif2}[]{\small ($d\le 3$)}
%
% scale 0.75
%
\begin{center}
\includegraphics{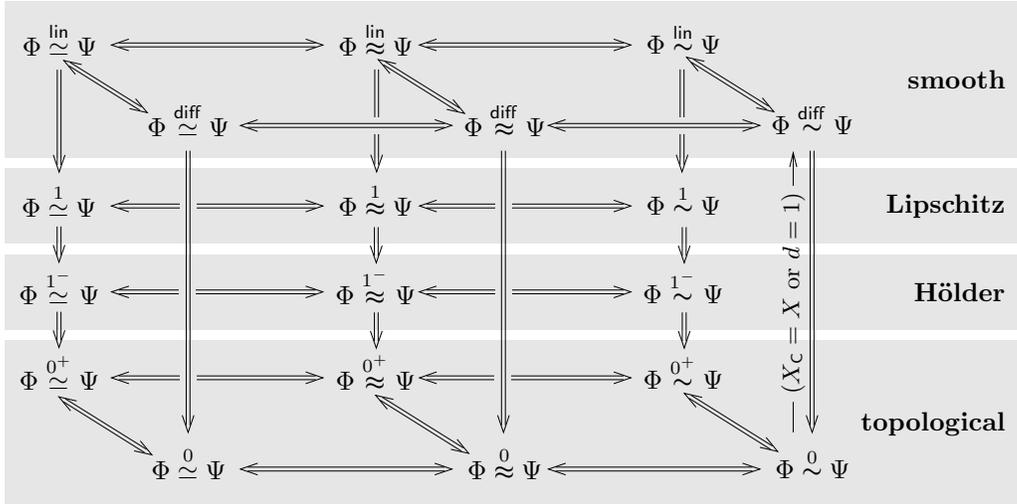}
\end{center}
\vspace*{-4mm}
\caption{As a consequence of Theorems \ref{thma} and \ref{thmb}, as
  well as Propositions \ref{thmc} and \ref{thmd}, all equivalences between
  linear flows $\Phi$, $\Psi$ on $X=\R^d$ coalesce into no more than
  {\em four\/} different forms.}\label{fig2} 
\end{figure}

The present section briefly recalls basic terminology and notation
pertaining to linear flows. As indicated above, a main objective of
this article is to demonstrate how all of the (at least eighteen, and
in fact infinitely many, as alluded to earlier) different equivalences
between arbitrary flows shown in Figure \ref{fig1} coalesce into
a mere {\em four\/} different forms --- provided that all flows
considered are linear; see Figure \ref{fig2}. Another, closely related
objective is to characterize, in combination with \cite{BW, BW3}, each
form of equivalence using basic linear algebra as outlined in the
Introduction. 

Regarding equivalence of arbitrary (not necessarily linear) flows
$\varphi$, $\psi$ on $X$, recall that the right implication in
(\ref{eq12}) cannot be reversed, and correspondingly $\varphi 
\stackrel{\bigstar}{\thicksim}\psi$ does not in general imply 
$\varphi \stackrel{\bigstar}{\thickapprox}\psi$ for any
$\bigstar$. It is a simple but consequential fact that this reverse
implication is valid for {\em linear\/} flows, for all six forms of (strict)
equivalence considered herein. 

\begin{prop}\label{lemH1}
Let $\Phi$, $\Psi$ be linear flows on $X$ and $\bigstar\in \{0, 0^+, 1^-, 1,
\mbox{\sf diff}, \mbox{\sf lin}\}$. Then $\Phi
\stackrel{\bigstar}{\thickapprox} \Psi$ if and only if $\Phi
\stackrel{\bigstar}{\thicksim} \Psi$. 
\end{prop}

Let $\Phi$ be a linear flow on $X$. A set $Y\subset X$ is
$\Phi$-{\bf invariant} if $\Phi_t Y = Y$ for every $t\in \R$, or
equivalently if $\Phi_{\R} y \subset Y$ for every $y\in Y$. A linear flow
$\Phi$ is {\bf irreducible} if $X= Z \oplus \widetilde{Z}$ with
$\Phi$-invariant sub{\em spaces\/} $Z$, $\widetilde{Z}$ implies that $Z=\{0\}$
or $\widetilde{Z} = \{0\}$. Thus, $\Phi$ is irreducible if and only
if, relative to an appropriate basis, the generator $A^{\Phi}$ is a single real
Jordan block. In particular, for an irreducible $\Phi$ the spectrum
$\sigma (\Phi):= \sigma (A^{\Phi})$, i.e., the set of all eigenvalues of
$A^{\Phi}$, is either a real singleton or a non-real complex conjugate
pair, that is, $\sigma(\Phi) = \{z, \overline{z}\}$ for some $z\in
\C$. Let $J_1 =[0]\in \R^{1\times 1}$, and for $m\in \N\setminus
\{1\}$ denote by $J_m$ the standard nilpotent $m\times m$-Jordan
block, 
$$
J_m = \left[ 
\begin{array}{ccccc}
0 & 1 & 0 & \cdots & 0 \\
\vdots & \ddots & \ddots &   &  \vdots \\
&  &  & \ddots & 0 \\
\vdots  & & & \ddots  & 1\\
0 & \cdots & & \cdots & 0
\end{array}
\right]\in \R^{m\times m} .
$$
Moreover, for every $m\in \N$ let
$$
J_m (a) = aI_m + J_m \, , \quad J_m(a+ib) = a I_{2m} +  \left[
\begin{array}{c|c}
J_m & - b I_m \\ \hline
b  I_m & J_m
\end{array}
\right]  \qquad \forall  a\in \R , b\in \R  \setminus \{0\} \, .
$$
For every $z\in \C$, therefore, $J_m(z)$ simply is a {\em real\/} Jordan block with $\sigma \bigl(
J_m(z)\bigr) = \{z, \overline{z}\}$. Note that $J_m(0) = J_m$ and $J_m^m
= O_m$; moreover, $J_m(z)\in \R^{m\times m}$ if $z\in \R$, whereas $J_m(z)\in \R^{2
  m\times 2m}$ if $z\in \C\setminus \R$. Observe that for any $a\in \R$,
$$
e^{tJ_m(a)} = e^{a t} e^{tJ_m} = e^{a t} \sum\nolimits_{j=0}^{m-1}
\frac{t^j}{j!} J_m^j \qquad \forall t \in \R \, ,
$$
whereas for any $a\in \R$, $b\in \R \setminus \{0 \}$,
$$
e^{tJ_m(a+ib)} = e^{a t } e^{t J_m(i b)} = e^{a t}
 \left[
\begin{array}{c|r}
\cos ( b t ) I_m  & - \sin ( b t )   I_m \\ \hline
\sin (b t ) I_m  &  \cos (b t )   I_m
\end{array}
\right]  
\mbox{\rm diag} \left[ e^{t J_m}, e^{t J_m}\right] \qquad \forall t \in \R \, .
$$

In general, given any linear flow $\Phi$ on $X$, recall that the subspaces
\begin{align*}
& X_{\sf S}^{\Phi} := \Bigl\{ x\in X : \lim\nolimits_{t\to \infty}
  \Phi_t x = 0 \Bigr\} \, , \\
& X_{\sf C}^{\Phi} := \Bigl\{ x\in X : \lim\nolimits_{|t|\to \infty}
  e^{-\varepsilon |t|}\Phi_t x = 0 \enspace \forall \varepsilon > 0 \Bigr\} \, , \\
& X_{\sf U}^{\Phi} := \Bigl\{ x\in X : \lim\nolimits_{t\to - \infty}
                                                                                      \Phi_t x = 0 \Bigr\} \, , \\
& X_{\sf H}^{\Phi}:= X_{\sf S}^{\Phi} \oplus X_{\sf U}^{\Phi} \, , 
\end{align*}
referred to as the {\bf stable}, {\bf central}, {\bf unstable}, and
{\bf hyperbolic} subspace of $\Phi$ respectively, are $\Phi$-invariant, and $X = X_{\sf
  S}^{\Phi} \oplus X_{\sf C}^{\Phi} \oplus X_{\sf
  U}^{\Phi} = X_{\sf H}^{\Phi} \oplus X_{\sf C}^{\Phi}$. Moreover,  say that $\Phi$ is {\bf stable}, {\bf central}, {\bf
  unstable}, and {\bf hyperbolic} if $X$ equals $X_{\sf S}^{\Phi}$, $X_{\sf
  C}^{\Phi}$, $X_{\sf U}^{\Phi} $, and $X_{\sf H}^{\Phi}$ respectively.
For convenience throughout, usage of the word {\em flow\/} in conjunction with
any of these adjectives, as well as {\em irreducible\/} or
{\em generated by}, automatically implies that the flow under
consideration is linear. 
Additionally, for $\bullet\in \{{\sf S}, {\sf C}, {\sf U}, {\sf H}\}$,
let $d_{\bullet}^{\Phi} = \dim X_{\bullet}^{\Phi}$,
write $\Phi|_{\R \times X_{\bullet}^{\Phi}}$ simply as
$\Phi_{\bullet}$, and denote by $P_{\bullet}^{\Phi}$ the linear
projection of $X$ onto $X_{\bullet}^{\Phi}$, along $\bigoplus_{\circ
  \in \{{\sf S}, {\sf C}, {\sf U} \} \setminus \{  \bullet \} } X_{\circ}^{\Phi}$ and $X_{\sf C}^{\Phi}$ if $\bullet \in
\{{\sf S}, {\sf C}, {\sf U}\}$ and $\bullet = {\sf H}$ respectively. Clearly, $\Phi \stackrel{{\sf
    lin}}{\cong} \bigtimes_{\bullet\in \{{\sf S}, {\sf C}, {\sf U}\}} \Phi_{\bullet}$ via the linear
isomorphism $h = \bigtimes_{\bullet\in \{{\sf S}, {\sf C}, {\sf U}\}}
P_{\bullet}^{\Phi}$, and $d_{\sf H}^{\Phi} = d_{\sf S}^{\Phi} + d_{\sf
  U}^{\Phi} = d - d_{\sf C}^{\Phi}$.
The {\bf time-reversal} $\Phi^*$ of $\Phi$ is the linear flow
on $X$ with $\Phi^*_t = \Phi_{-t}$ for every $t\in \R$; in other
words, $\Phi^*$ is generated by $-A^{\Phi}$. Obviously, 
$\Phi^* \stackrel{\sf lin}{\simeq} \Phi$, and $X_{\sf S}^{\Phi^*} =
X_{\sf U}^{\Phi}$, $X_{\sf C}^{\Phi^*} =
X_{\sf C}^{\Phi}$, $X_{\sf U}^{\Phi^*} =
X_{\sf S}^{\Phi}$, as well as $X_{\sf H}^{\Phi^*} =
X_{\sf H}^{\Phi}$.

The following is a simple but useful general observation regarding the
H\"{o}lder property of maps  relative to a decomposition of $X$ into
complementary subspaces. 

\begin{prop}\label{prop3add0}
Let $Y$, $Z$ be subspaces of $X$ with $X=Y\oplus Z$ and $0\le \beta
\le 1$.
\begin{enumerate}
\item If $h\in \cH_{\beta}(X)$ and $h(Y) =Y$, then $h|_Y \in
  \cH_{\beta}(Y)$.
\item If $f\in \cH_{\beta}(Y)$ and $g\in \cH_{\beta}(Z)$, then $f\times
  g\in \cH_{\beta} (X)$; here $f\times g (y+z) = f(y) + g(z)$ for every
  $y\in Y$, $z\in Z$.
\end{enumerate}
\end{prop}

For the analysis in subsequent sections, it is helpful to
recall one further classical concept: Given any linear flow
$\Phi$ on $X$, the (forward) {\bf Lyapunov exponent}
$$
\lambda_+^{\Phi} (x) = \lim\nolimits_{t\to \infty} \frac{\log |\Phi_t x|}{t}
$$
exists for every $x\in X\setminus \{0\}$, and the range of $x\mapsto
\lambda_+^{\Phi}(x)$ equals $\{{\sf Re}z : z \in \sigma(\Phi)\}$. With
$\lambda_+^{\Phi}(0):= -\infty$ for convenience, the set $L^{\Phi}
(s):= \{x\in X : \lambda_+^{\Phi}
(x) \le s\}$ is a $\Phi$-invariant subspace for every $s\in \R$,
referred to as the {\bf Lyapunov space} of $\Phi$ at $s$. Writing, for
every linear operator $A$ on $X$ and $z\in \C$,
$$
\mbox{\rm gker}\, (A - z I_X) := \bigcup\nolimits_{n\in \N} \ker \, (A^2 -
2 \, {\sf
  Re} z A + |z|^2 I_X)^n \, ,
$$
it is readily seen that $L^{\Phi}(s) = \sum_{{\sf Re} z \le s}
\mbox{\rm gker}\, (A^{\Phi} - zI_X)$ for every $s\in \R$. (In
\cite{CK,KS} the term {\em Lyapunov space\/} instead refers to any subspace
$ \sum_{{\sf Re} z = s} \mbox{\rm gker}\, (A^{\Phi} -
zI_X)$. Such spaces may behave poorly under equivalence, and with the
(good) behaviour of key objects being crucial for the present
article this terminology is not adopted here.) Letting
$\ell^{\Phi}(s)= \dim L^{\Phi}(s)$, clearly the
integer-valued function $\ell^{\Phi}$ is non-decreasing and
right-continuous, with $\lim_{s\to -\infty}\ell^{\Phi}(s) = 0$ and
$\lim_{s\to \infty} \ell^{\Phi}(s) = d$. Observe that $\lambda_+^{\Phi} (x) = a$ for
some $x\in X$, $a\in \R$ precisely if $\ell^{\Phi} (a) > \ell^{\Phi}
(a^-)$, and refer to the non-negative integer $\ell^{\Phi} (a) - \ell^{\Phi}
(a^-)$ as the {\bf multiplicity} of $a$. Let $\lambda_1^{\Phi}\le \lambda_2^{\Phi} \le \ldots \le
\lambda_d^{\Phi}$ be the (not necessarily different) Lyapunov
exponents of the linear flow $\Phi$, that is, $\bigl\{ \lambda_j^{\Phi} : j\in \{1,
\ldots , d\} \bigr\} = \bigl\{ \lambda_+^{\Phi} (x): x\in X \setminus
\{0\} \bigr\}$, with each exponent repeated according to its 
multiplicity; see, e.g., \cite{BaPe, CK} for authoritative accounts
of the theory and applications of Lyapunov exponents. For convenience, let
$$
\Lambda^{\Phi} := \Lambda^{A^{\Phi}} := \mbox{\rm diag} \, 
[\lambda_1^{\Phi}, \ldots , \lambda_d^{\Phi}] \, .
$$
Note that if
$\Phi$ is irreducible with $\sigma (\Phi) = \{z,\overline{z}\}$ for
some $z\in \C$, then simply $\Lambda^{\Phi} = {\sf Re} z I_d$. Also,
$\Phi$ is stable, unstable, central, and hyperbolic precisely if
$\lambda_j^{\Phi}<0$, $\lambda_j^{\Phi}>0$, $\lambda_j^{\Phi}=0$, and
$\lambda_j^{\Phi}\ne 0$ for every $j\in \{1, \ldots , d\}$,
respectively. Moreover, $\ell^{\Phi^*}(-s) = d - \ell^{\Phi}(s^-)$ for
all $s\in \R$, and consequently $\lambda_j^{\Phi^*} = -
\lambda_{d+1-j}^{\Phi}$ for every $j\in \{1, \ldots , d\}$, i.e.,
$\Lambda^{\Phi^*} = -\mbox{\rm diag}\, [\lambda_d^{\Phi}, \ldots ,
\lambda_1^{\Phi}]$. Say that two linear flows $\Phi, \Psi$, or their generators, are {\bf Lyapunov
  similar} if $\Lambda^{\Phi} = \Lambda^{\Psi}$. (In \cite{ACK1} the
term {\em Lyapunov equivalent\/} is used instead.) Thus
$\Phi$, $\Psi$, or $A^{\Phi}$, $A^{\Psi}$, are Lyapunov similar
precisely if they have the same Lyapunov exponents, with matching
multiplicities, or equivalently if $\ell^{\Phi} = \ell^{\Psi}$. Note
that if $A^{\Phi}$, $A^{\Psi}$ are similar then clearly $\ell^{\Phi} =
\ell^{\Psi}$, whereas the converse is not true in general for $d\ge 2$.

\begin{rem}\label{rem32}
(i) This article is based entirely on (\ref{eq01}),
whereby {\em equivalence\/} between flows on $X=\R^d$ means the
preservation of all orbits, up to a bijection $h:X\to X$ that exhibits 
some additional regularity. Without such regularity
this approach would be too crude to be truly meaningful: For instance, for $d\ge 4$
and any linear flow $\varphi$ on $X$, (\ref{eq01}) holds with $\psi$
generated by precisely one of either 
$$
O_d, \, \mbox{\rm diag}\, [O_{d-1}, 1] , \,  \mbox{\rm diag}\,
[O_{d-2}, J_1(i)], \,  \mbox{\rm diag}\, [O_{d-3}, J_1(i), 1] , \,
\mbox{\rm diag}\, [J_1(i), I_{d-2}], \, I_d \, ,
$$
or, in case $d$ is even, $ \mbox{\rm diag}\, [J_1(i), \ldots ,
J_1(i)]$. However, the bijection $h$ may fail to be
measurable, let alone continuous, $\beta$-H\"{o}lder, etc. 

(ii) Equivalence between flows on
$X$ can of course be defined differently altogether. To see but one such
definition specifically for linear flows, say that $\Phi$, $\Psi$
are {\bf kinematically similar}, in symbols $\Phi  \leftrightarrows  \Psi$, if there
exists an invertible linear operator $Q$ on $X$ so that 
\begin{equation}\label{eq3addxx1}
\sup\nolimits_{t\in \R} (\|\Phi_t Q^{-1} \Psi_{-t}\| + \|\Psi_t Q
\Phi_{-t}\|) < \infty \, ,
\end{equation}
where $\|\cdot\|$ denotes any operator norm; see, e.g., \cite[Sec.\
5]{cop}. To relate this classical concept to the present
article, note on the one hand that if $A^{\Phi}$, $A^{\Psi}$ are
similar, say $QA^{\Phi} = A^{\Psi} Q$, then $\Phi_t Q^{-1} \Psi_{-t} =
Q^{-1}$ for all $t\in \R$, so (\ref{eq3addxx1}) automatically holds. On the
other hand, $\Phi$, $\Psi$ are readily seen to be Lyapunov similar if
and only if
\begin{equation}\label{eq3addxx2}
\sup\nolimits_{t\in \R} (\|\Phi_t Q^{-1} \Psi_{-t}\| + \|\Psi_t Q
\Phi_{-t}\|)e^{-\varepsilon |t|}< \infty \qquad \forall \varepsilon > 0
\,  .
\end{equation}
Clearly, (\ref{eq3addxx1}) implies (\ref{eq3addxx2}), and this
implication is not reversible for $d\ge 2$. Also, it turns out that
$\Phi  \leftrightarrows  \Psi$ can be 
characterized easily in terms of $A^{\Phi}$, $A^{\Psi}$, and
\begin{equation}\label{eq3om}
  \Phi \stackrel{1}{\cong} \Psi \quad \Longrightarrow \quad
  \Phi  \leftrightarrows   \Psi \quad \Longrightarrow \quad
   \Phi \stackrel{1^-}{\cong} \Psi \, ;
 \end{equation}
see, e.g., \cite[Sec.\ 4]{KS}. Thus $\Phi  \leftrightarrows  \Psi$ entails (\ref{eq02}) for
some $h\in \cH_{1^-}$ and $\alpha = 1$. However, the precise regularity
of $h$ is not characterized by $\Phi  \leftrightarrows  \Psi$. For instance, $\Phi
 \leftrightarrows  \Psi$ with $\Phi$, $\Psi$ generated by $\mbox{\rm diag}\, [J_1(1+i),
J_1(1+i)]$,  $I_4$ respectively, but also with $\mbox{\rm diag}\, [J_2(1),
J_2(1)]$, $J_2(1+i)$ instead; in the former case, $ \Phi \stackrel{1}{\cong} \Psi$ whereas
in the latter case $ \Phi \, \cancel{\stackrel{1}{\thicksim}} \, 
\Psi$. Similarly, $\Phi  \leftrightarrows  \Psi$ may or may not hold whenever $\Phi
\stackrel{1^-}{\thicksim} \Psi $ but $\Phi \, \cancel {\stackrel{1}{\thicksim}}\, \Psi
$. These examples also illustrate how neither implication in
(\ref{eq3om}) can be reversed in general for $d\ge 4$. 
\end{rem}

%%%%%%%%%%%%%%%%%%%%%%%%%%%%%%%%%%%%%%%%%%%%%%%%%%%%%%%%%%%%%%%%

\section{$\beta$-H\"{o}lder relations between (un)stable flows}\label{sec3}

%%%%%%%%%%%%%%%%%%%%%%%%%%%%%%%%%%%%%%%%%%%%%%%%%%%%%%%%%%%%%%%%

This section studies $\stackrel{\bigstar}{\thicksim}$ for $\bigstar
\in \{\beta^- , \beta , \beta^+\}$ with $0<\beta < 1$. Although they
are not transitive, a careful analysis of these relations, at least in
the context of stable or unstable flows, nonetheless is
essential for the H\"{o}lder classification(s) to be established in
Section \ref{sec5}. 

To study $\beta$-H\"{o}lder relations between flows, first
consider irreducible flows. The following four lemmas establish
all-H\"{o}lder relations between flows generated by $J_m(z)$
with $m\in \N$ and $z\in \C \setminus i\R$. The case of $z= a \in \R
\setminus \{0\}$ is especially simple.

\begin{lem}\label{lemr1}
Given $m\in \N$ and $a\in \R \setminus \{0\}$, let $\Phi$, $\Psi$ be
the flows on $\R^m$ generated by $J_m(a)$, $aI_m$
respectively. Then $\Phi \stackrel{1^-}{\thicksim} \Psi$, i.e., $\Phi$,
$\Psi$ are all-H\"{o}lder equivalent.
\end{lem}

\begin{proof}
It will be shown that in fact $\Phi \stackrel{1^-}{\cong} \Psi$, i.e.,
$\Phi$, $\Psi$ are all-H\"{o}lder conjugate. Once established, clearly
this stronger assertion proves the claim. Since even the stronger
assertion trivially is correct for $m=1$, henceforth assume
$m\ge 2$. Consider the map $h_a:\R^m \to \R^m$ given by
\begin{equation}\label{eqr2}
h_a(x)_j = \sum\nolimits_{k=0}^{m-j} \frac{ (\log |x_{m+1 - j
  - k }|)^k}{k! \, a^k} \,  x_{m+1 - j - k } \qquad \forall x \in \R^m, j \in \{1, \ldots , m \} \, ,
\end{equation}
with the convention that $0 (\log 0)^k = 0$ for every $k\in
\N_0$. Note that $h_a(0) = 0$, and since $u \mapsto u (\log |u|)^k$
satisfies a $\beta$-H\"{o}lder condition near $0$ for every $0<\beta <
1$, so does $h_a$. Moreover, it is readily seen that $h_a$ is a
homeomorphism, with the components of $h_a^{-1}$ determined
recursively by 
$$
h_a^{-1}(x)_j = x_{m+1 - j} - \sum\nolimits_{k=1}^{j-1} \frac{ (\log
  |h_a^{-1} (x)_k|)^{j-k} }{(j-k)!\, 
  a^{j-k}} \,  h_a^{-1} (x)_k \qquad \forall
x \in \R^m, j\in \{1, \ldots , m\} \, ;
$$
in particular, $h_a^{-1}$ also satisfies a $\beta$-H\"{o}lder
condition near $0$ for every $0<\beta < 1$. Note that
$$
h_a(re_1) = r \sum\nolimits_{j=0}^{m-1}\frac{(\log r)^{j}}{j!\, 
  a^{j}} \, e_{m-j} \qquad \forall r \in \R^+ \, ,
$$
so clearly $h_a$ does not satisfy a Lipschitz condition near $0$. In
summary, therefore, $h_a \in \cH_{1^-} \setminus \cH_1$ whenever $m\ge 2$.
Now, observe that for every $t\in \R$, $x\in \R^m$, and $j\in \{1, \ldots , m\}$,
\begin{align*}
h_a(e^{a t} x)_j & = e^{at} \sum\nolimits_{k=0}^{m-j} \frac{(at  +
                   \log |x_{m+1 - j - k}|)^k}{k! \, a^k} \, x_{m+1 - j
                  - k}   \\
  & = e^{at} \sum\nolimits_{k=0}^{m-j} \sum\nolimits_{\ell =0}^k
    \left( \!\!  \begin{array}{c} k \\ \ell \end{array}\!\!\right)
  \frac{a^{\ell}  t^{\ell} (\log |x_{m+1 - j - k}|)^{k-\ell} }{k!\,  a^k}
  \, x_{m+1 - j - k}  \\
  & = e^{at} \sum\nolimits_{\ell = 0}^{m-j} \frac{t^{\ell}}{\ell !}
    \sum\nolimits_{k= \ell}^{m-j} \frac{(\log |x_{m+1 - j -
    k}|)^{k-\ell} }{(k-\ell)! \, a^{k-\ell}} x_{m+1 - j - k}   = e^{at} \sum\nolimits_{\ell = 0}^{m-j} \frac{t^{\ell}}{\ell !}
    h_a(x)_{j+\ell}\\
  & = \bigl( e^{t J_m(a)} h_a(x)\bigr)_j  \, ,
\end{align*}
and consequently
$$
h_a (\Psi_t x) = h_a (e^{at}x) = e^{tJ_m(a)} h_a(x) = \Phi_t h_a(x)
\qquad \forall t \in \R ,  x\in \R^m \, .
$$
In other words, $\Psi \stackrel{h_a}{\cong} \Phi$, and so $\Phi
\stackrel{1^-}{\cong} \Psi$ as claimed.
\end{proof}

The next result is an analogue of Lemma \ref{lemr1} for the case of
$z=a+ib \in \C \setminus (\R \cup i\R)$.

\begin{lem}\label{lemr2}
Given $m\in \N$ and $a,b\in \R \setminus \{0\}$, let $\Phi$, $\Psi$ be
the flows on $\R^{2m}$ generated by $J_m(a+ib)$, $\mbox{\rm
  diag} \, [J_1(a+ib), \ldots , J_1(a+ib)]$ respectively. Then $\Phi
\stackrel{1^-}{\thicksim} \Psi$.
\end{lem}

\begin{proof}
Again, it turns out that in fact  $\Phi
\stackrel{1^-}{\cong} \Psi$, and it is this stronger assertion
that will be established here. Since there is nothing to prove for
$m=1$, henceforth assume $m\ge 
2$. To mimic the proof of Lemma \ref{lemr1}, for
every $j\in \{1, \ldots , m\}$ let $E_j = \mbox{\rm span} \{e_{2j-1} ,
e_{2j}\}$, and denote by $P_j$ the orthogonal projection of $\R^{2m}$ onto $E_j$. In
analogy to (\ref{eqr2}), consider $h_{a+ib} :\R^{2m} \to \R^{2m}$ given by
$$
\left[ \begin{array}{c}
         h_{a+ib} (x)_j \\
         h_{a+ib} (x)_{j+m}
       \end{array}
     \right]
     = \sum\nolimits_{k=0}^{m-j} \frac{(\log
       |P_{m+1 - j - k} x|)^{k}}{k! \, a^k}\,\left[
       \begin{array}{c}
         x_{2(m+1-j-k)-1} \\
         x_{2(m+1-j-k)}
       \end{array}
     \right]
     \qquad \forall x \in
\R^{2m} ,  j \in \{1, \ldots , m\} \, .
$$
As in the proof of Lemma \ref{lemr1}, it is readily seen that
$h_{a+ib} \in \cH_{1^-}\setminus \cH_1$, and an essentially identical
calculation yields, for every $t\in\R$, $x\in \R^{2m}$, and $j\in \{1, \ldots , m\}$,
\begin{align*}
h_{a+ib} \left(
  \mbox{\rm diag} \! \left[
e^{t J_1(a+ib)} , \ldots , e^{t J_1(a+ib)}
    \right] x
  \right)_j & = \left(  e^{t J_m(a+ib)} h_{a+ib} (x)\right)_j   \, ,
  \\
  h_{a+ib} \left(
  \mbox{\rm diag} \! \left[
e^{t J_1(a+ib)} , \ldots , e^{t J_1(a+ib)}
    \right] x
  \right)_{j+m} & = \left(  e^{t J_m(a+ib)} h_{a+ib} (x)\right)_{j+m}
                  \, .
\end{align*}
In other words, for every $t\in \R$ and $x\in \R^{2m}$,
$$
h_{a+ib} (\Psi_t x) = h_{a+ib} \left(
  \mbox{\rm diag} \! \left[
e^{t J_1(a+ib)} , \ldots , e^{t J_1(a+ib)}
    \right] x
  \right) = e^{t J_m (a+ib)} h_{a+ib} (x) = \Phi_t h_{a+ib} (x) \, ;
$$
that is, $\Psi \stackrel{h_{a+ib}}{\cong} \Phi$, and so $\Phi
\stackrel{1^-}{\cong} \Psi$ as claimed.
\end{proof}

Each individual block $J_1(a+ib)$ appearing in Lemma \ref{lemr2} can be
simplified further by means of an equivalence that is even more regular.

\begin{lem}\label{lemr3}
Given $a,b\in \R\setminus \{0\}$, let $\Phi$, $\Psi$ be the 
flows on $\R^2$ generated by $J_1(a+ib)$, $aI_2$ respectively. Then
$\Phi \stackrel{1}{\thicksim}\Psi$, i.e., $\Phi$, $\Psi$ are Lipschitz equivalent.
\end{lem}

\begin{proof}
  For convenience, let $R_s = e^{s J_1(i)} = \left[
    \begin{array}{cr}
\cos s & - \sin s \\ \sin s & \cos s
      \end{array}
    \right]\in \R^{2\times 2}$ for every $s\in \R$. The map $g:
    \R^2 \to \R^2$ given by $g(0) = 0$ and
    $$
g(x) = R_{-b\log |x|/a} x  \qquad \forall x \in \R^2 \setminus \{0\}
\, ,
    $$
is a bi-Lipschitz homeomorphism, with $g^{-1}(x) = R_{b\log
  |x|/a}x$ for $x\ne 0$. Thus $g\in \cH_1$; furthermore, for every $x\in \R^2
\setminus \{0\}$,
$$
g(\Phi_t x) = g(e^{at} R_{bt}x) = R_{-bt - b\log|x|/a} (e^{at} R_{bt}
x) = e^{at} R_{-b\log|x|/a}x= e^{at} g(x) = \Psi_t g(x) \qquad \forall
t \in \R \, ,
$$
and the two outermost expressions agree for $x=0$ also. Thus $\Phi
\stackrel{g}{\cong} \Psi$, and hence $\Phi \stackrel{1}{\cong} \Psi$,
i.e., $\Phi$, $\Psi$ are Lipschitz conjugate. Clearly, therefore,
$\Phi \stackrel{1}{\thicksim}\Psi$ as well. 
\end{proof}

Using Lemma \ref{lemr3}, it is straightforward to bring Lemma
\ref{lemr2} fully in line with Lemma \ref{lemr1}.

\begin{lem}\label{lemr4}
Given $m\in \N$ and $a,b\in \R \setminus\{0\}$, let $\Phi$, $\Psi$ be
the flows on $\R^{2m}$ generated by $J_m(a+ib)$, $a I_{2m}$
respectively. Then $\Phi \stackrel{1^-}{\thicksim} \Psi$.
\end{lem}

\begin{proof}
Denote by $\widetilde{\Phi}$ the flow on $\R^{2m}$ generated by
$\mbox{\rm diag}\, [J_1(a+ib), \ldots , J_1(a+ib)]$. By Lemma
\ref{lemr2}, $\Phi \stackrel{1^-}{\thicksim} \widetilde{\Phi}$, and by
Lemma \ref{lemr3}, $\widetilde{\Phi} \stackrel{1}{\thicksim} \Psi$, via
the $m$-fold product $g\times \ldots \times g$, the latter being Lipschitz
due to Proposition \ref{prop3add0}. Hence $\Phi
\stackrel{1^-}{\thicksim} \Psi$, by the transitivity of
$\stackrel{1^-}{\thicksim}$.
\end{proof}

From Lemmas \ref{lemr1} and \ref{lemr4}, it is readily deduced that
any two irreducible flows are all-H\"{o}lder equivalent, provided that
they are hyperbolic.

\begin{prop}\label{propr4A}
Let $\Phi$, $\Psi$ be irreducible flows on $X$ with $\sigma
(\Phi), \sigma(\Psi)\subset \C \setminus i\R$. Then $\Phi
\stackrel{1^-}{\thicksim}\Psi$. 
\end{prop}

Letting $\Phi$, $\Psi$ be the flows on $\R^2$
generated by $J_2(1)$, $J_1(1+i)$ respectively, shows that the conclusion
$\Phi \stackrel{1^-}{\thicksim}\Psi$ in Proposition \ref{propr4A}
cannot in general be strengthened to $\Phi 
\stackrel{1}{\thicksim}\Psi$ when $d\ge 2$; see also \cite{BW3}.

When extending Lemmas \ref{lemr1} and \ref{lemr4} to arbitrary (un)stable
flows, one may suspect that the presence of two or more irreducible
components for $\Phi$, $\Psi$ will decrease the maximal possible
regularity of $h$ in $\Phi \stackrel{h}{\thicksim}\Psi$, if indeed
$\Phi$, $\Psi$ are related at all. The remainder of the present
section confirms this suspicion by providing a detailed analysis of $\Phi
\stackrel{h}{\thicksim}\Psi$ with $h\in \cH_{\beta}$ and $0<\beta <
1$, assuming $\Phi$, $\Psi$ to both be (un)stable. In this analysis,
as well as in subsequent sections, the {\em topological invariance of
  dimension\/} is used in its following basic form; see, e.g.,
\cite[Sec.\ 2B]{hatch}.

\begin{prop}\label{prop4top}
Given $m,n\in \N$, let $U\subset \R^m$ be non-empty and open. There
exists a continuous one-to-one function $f:U\to \R^n$ if and only if $m\le n$.
\end{prop}

To extend Proposition \ref{propr4A}, let $\Phi$ be a stable
flow. (For an unstable flow, simply consider its time-reversal
instead.) Lemmas \ref{lemr1} and \ref{lemr4}, applied individually to
each irreducible component, together with Proposition \ref{prop3add0}, show
that $\Phi$ is all-H\"{o}lder equivalent to the flow generated by
$\Lambda^{\Phi}$. As far as $\beta$-H\"{o}lder relations between
stable flows are concerned, therefore, it suffices to study flows
generated by $\mbox{\rm diag} \,  [a_1, \ldots , a_m]$ with negative
$a_j$; for convenience, fix $a,b\in \R^m$ with
\begin{equation}\label{eqr9}
a_1 \le \ldots \le a_m < 0 \quad \mbox{\rm and} \quad b_1 \le
\ldots \le b_m < 0 \, .
\end{equation}
The following result characterizes $\Phi
\stackrel{\beta}{\thicksim} \Psi$ for any two flows $\Phi$, $\Psi$
thus generated.

\begin{theorem}\label{lemr9}
Given $m\in \N$ and $a,b\in \R^m$ as in
{\rm (\ref{eqr9})}, let $\Phi$, $\Psi$ be the flows on $\R^m$
generated by $\mbox{\rm
  diag}\, [a_1, \ldots , a_m]$, $\mbox{\rm
  diag}\, [b_1, \ldots , b_m]$ respectively. For every $0<\beta< 1$
the following are equivalent:
\begin{enumerate}
\item $\Phi \stackrel{\beta}{\simeq} \Psi$;
\item $\Phi \stackrel{\beta^-}{\thicksim} \Psi$;
\item $\beta^2 \le {\displaystyle \frac{\min_{j=1}^m (a_j/b_j)}{\max_{j=1}^m (a_j/b_j)}}$.
 \end{enumerate}
\end{theorem}

The proof of Theorem \ref{lemr9} makes use of the
elementary fact that, informally put, two
{\em different\/} $\Phi$-orbits cannot approach one another at a rate
faster than $e^{a_1 t}$ as $t\to \infty$. To state this precisely, as
usual let $\mbox{\rm dist} (x,W) = \inf_{w\in W}|x-w|$ for any $x\in
\R^m$ and $\varnothing \ne W \subset \R^m$.

\begin{lem}\label{lemr9A}
Given $m\in \N$, $x,y\in  \R^m$, $s\in \R$, and with
$\Phi$ as in Theorem \ref{lemr9},
$$
|\Phi_t x - \Phi_{s} y| \ge e^{a_1  t} \mbox{\rm dist} (x , \Phi_{\R}
y) \qquad \forall t \ge 0 \, .
$$
\end{lem}

\begin{proof}
Note that for every $t\ge 0$,
\begin{align*}
e^{-2a_1  t} |\Phi_t x - \Phi_{s} y|^2 & =
                                                \sum\nolimits_{j=1}^m
                                                \left( e^{(a_j -
                                                a_1 )t} x_j - e^{a_j
                                         s - a_1  t}
                                                y_j\right)^2 \\
  & = \sum\nolimits_{j=1}^m e^{2(a_j - a_1)t} \left( x_j - e^{a_j(s
    - t)}y_j\right)^2  \\
  & \ge \sum\nolimits_{j=1}^m  \left( x_j - e^{a_j(s
    - t)}y_j\right)^2 =   |x - \Phi_{s - t} y|^2 \ge \mbox{\rm dist} ( x,  \Phi_{\R}
    y )^2     \, ,
\end{align*}
where the first inequality is due to $(a_j - a_1)t \ge 0$ for every $j\in \{1,
\ldots , m\}$.
\end{proof}

\begin{proof}[Proof of Theorem \ref{lemr9}]
For $m=1$, all three statements are true for every $0<\beta < 1$,
as $\Phi \stackrel{{\sf lin}}{\simeq} \Psi$, and (iii) reads $\beta^2
\le 1$. Hence, assume $m\ge 2$ from now on. Obviously
(i)$\Rightarrow$(ii) by definition. 

To prove that (ii)$\Rightarrow$(iii),
fix any $0<\gamma < \beta$, and assume that $\Phi
\stackrel{h}{\thicksim}\Psi$ with some $h\in \cH_{\gamma}$. Note that
$\tau_x$ is increasing for every $x\in \R^m\setminus \{0\}$. Throughout the proof, it
will be useful to adopt the following classically-inspired notation
\cite{Hardy}: Given any two functions $f, g:\R \to \R^+$, write
$f(t)\prec g(t)$, or equivalently $g(t)\succ f(t)$ whenever
$\limsup_{t\to \infty} f(t)/g(t)<\infty$, and write $f(t)\asymp g(t)$
if both $f(t)\prec g(t)$ and $f(t)\succ g(t)$. Importantly, $\prec$ is
reflexive and transitive, and $\asymp$ is an equivalence relation.

Now, let $E_j = \mbox{\rm
  span}\{e_1 , \ldots , e_j\}$ for every $j\in \{1, \ldots , m\}$ and
fix $j\ge 2$. Since $h(E_j) \not \subset E_{j-1}$ by Proposition
\ref{prop4top}, and since $E_j \setminus
E_{j-1}$ is dense in $E_j$, there exists an $x\in E_j \setminus E_{j-1}$
so that $h(x)\not \in E_{j-1}$; in addition, it can be assumed that
$x_1 \cdot \ldots \cdot x_j \ne 0$. Then $|\Phi_t x|\asymp e^{a_j t}$,
and consequently $|h(\Phi_t x)|\prec e^{\gamma a_j t}$, but also
$$
|h(\Phi_t x)| = |\Psi_{\tau_x(t)} h(x)| \succ e^{b_j \tau_x (t)} \, .
$$
The transitivity of $\prec$ yields $e^{b_j \tau_x(t)} \prec e^{\gamma a_j
t}$, and hence
\begin{equation}\label{eqpp1}
\liminf\nolimits_{t\to  \infty} \left( \tau_x(t) - \frac{\gamma
    a_j}{b_j} \, t\right) > - \infty \, .
\end{equation}
Next, fix any $k\in \{1,\ldots, j-1\}$, and let
$$
y_w = h^{-1} (h(x) + w) \qquad \forall w \in E_k \, ;
$$
here usage of the subscript $w$ highlights the $w$-dependence of
$y_{w}$. Clearly $y_0 = x$. Moreover,
$$
|\Psi_{\tau_x(t)} h(y_w) - \Psi_{\tau_x(t)}h(x)| = |\Psi_{\tau_x(t)}
w| \prec e^{b_k \tau_x(t)} \qquad \forall w \in E_k \setminus \{0\} \, ,
$$
and consequently, with $\sigma_w := \tau_{y_w}^{-1} \circ \tau_x$,
\begin{equation}\label{eqpp1A}
|\Phi_{\sigma_w(t)} y_w - \Phi_t x| \prec e^{\gamma b_k \tau_x(t)}
\qquad \forall w \in E_k \setminus \{0\} \, .
\end{equation}
It will be shown below that
\begin{equation}\label{eqpp2}
\mbox{\rm dist} ( \Phi_t x ,  \Phi_{\R} y_w ) \succ
e^{a_k t} \quad \mbox{\rm for some } w\in E_k\setminus \{0\}  \, .
\end{equation}
Assuming (\ref{eqpp2}) for the time being, observe how (iii) follows
rather directly from it: Indeed, picking $w\in E_k \setminus \{0\}$ as
in (\ref{eqpp2}) implies, together with (\ref{eqpp1A}), that
$e^{a_k t} \prec e^{\gamma b_k \tau_x(t)}$, and hence
\begin{equation}\label{eqpp3}
\limsup\nolimits_{t\to \infty} \left( \tau_x(t) -\frac{a_k}{\gamma
    b_k} \, t\right) <  \infty \, .
\end{equation}
Combining (\ref{eqpp1}) and (\ref{eqpp3}) yields
\begin{equation}\label{eqpp3A}
\frac{\gamma a_j}{b_j} \le \frac{a_k}{\gamma b_k} \, .
\end{equation}
So far, it has been assumed that $1 \le k < j \le m$, but obviously
(\ref{eqpp3A}) is correct also for $j=k$. Since $\gamma <
\beta$ has been arbitrary,
\begin{equation}\label{eqpp4}
\beta^2 \le \frac{a_k/b_k}{a_j/b_j} \qquad \forall 1 \le k \le j \le m
\, .
\end{equation}
Identical reasoning with the roles of $\Phi$, $\Psi$ interchanged
yields (\ref{eqpp4}) with $j$, $k$ interchanged. In summary,
$$
\beta^2 \le \frac{a_k/b_k}{a_j/b_j} \qquad \forall j,k\in \{1,
\ldots , m\} \, ,
$$
which immediately implies (iii).

To complete the proof of (ii)$\Rightarrow$(iii), it remains to
establish (\ref{eqpp2}). For this, let $E = \mbox{\rm span} \{e_k,
\ldots , e_m\}$, and denote by $P$ the orthogonal projection of $\R^m$
onto $E$. The subspace $E$ is
$\Phi$- and $\Psi$-invariant. Denoting the
restriction $\Phi|_{\R \times E}$ by $\widetilde{\Phi}$ for convenience, observe that $\widetilde{\Phi}$ can
be identified with the flow on $\R^{m-k+1}$ generated by $\mbox{\rm
  diag}\, [a_k, \ldots , a_m]$. Clearly $\widetilde{\Phi}_t P = P
\Phi_t$ for all $t\in \R$. Moreover, recall that $x_1 \cdot \ldots \cdot
x_j \ne 0$, and let
$$
\cM_x = \left\{ z \in \R^m : \frac{z_k}{x_k}>0, \ldots ,\frac{
    z_j}{x_j} >0, z_{j+1} =
  \ldots = z_m = 0  \enspace
  \mbox{\rm and}\enspace \left( \frac{z_k}{x_k}\right)^{1/a_k} =
  \ldots = \left( \frac{z_j}{x_j}\right)^{1/a_j}
  \right\} \, .
  $$
The set $\cM_x\subset \R^m$ is $\Phi$-invariant, and $x\in
\cM_x$. Given any $z\in \R^m$, note that $z\in \cM_x$ precisely if $Pz \in P \Phi_{\R}
x$. To see that $(h(x) + E_k)\setminus h(\cM_x) \ne \varnothing$,
suppose by way of contradiction that
\begin{equation}\label{eq4pp1}
h(\cM_x) \supset h(x) + E_k \, .
\end{equation}
Then $y_w \in \cM_x$ for every $w\in E_k$, and hence, by the $\Phi$-invariance of $\cM_x$,
$$
\Psi_t (h(x) + w) = \Psi_t h(y_w) = h (\Phi_{\tau_{y_w}^{-1}(t)} y_w) \in
h(\cM_x) \qquad \forall (t,w) \in \R \times E_k \, .
$$
Thus, with $\cC_x:= \{ \Psi_t (h(x) + w): t\in \R, w \in E_k \}$ for
convenience, (\ref{eq4pp1}) implies that
\begin{equation}\label{eq4pp2}
h(\cM_x) \supset \cC_x \, .
\end{equation}
The map $f:\R \times E_k \to \cC_x$ given by $f(t,w) = \Psi_t (h(x) +
w)$ is continuous and onto; since $h(x) \not \in E_k$ it also
is one-to-one. Consequently, $\cC_x$ is homeomorphic to $\R \times
E_k$, and hence to $\R^{k+1}$. By contrast, $\cM_x$ is homeomorphic to
$\R^+\times E_{k-1}$, and hence to $\R^k$. Thus (\ref{eq4pp2}) and
indeed (\ref{eq4pp1}) are impossible by Proposition
\ref{prop4top}. In other words, $(h(x) +
E_k)\setminus h(\cM_x)\ne \varnothing$ as claimed. Pick any $w\in E_k$ with
$h(x)+ w \not \in h(\cM_x)$, that is, $y_w \not \in \cM_x$. Then
$Px \not \in \widetilde{\Phi}_{\R} Py_w = P \Phi_{\R}y_w$, and Lemma
\ref{lemr9A} applied to $Px$, $Py_w$, and $\widetilde{\Phi}$ yields
$$
\big|\widetilde{\Phi}_t Px - \widetilde{\Phi}_s P y_w\big| \ge e^{a_k t}
\mbox{\rm dist} \bigl( Px, \widetilde{\Phi}_{\R} P y_w \bigr) \qquad \forall t \ge
0, s\in \R  \, .
$$
With $c:= \mbox{\rm dist} (Px, P\Phi_{\R} y_w)>0$, therefore,
$$
|\Phi_t x - \Phi_s y_w| \ge \big|\widetilde{\Phi}_t Px -
\widetilde{\Phi}_s P y_w\big| \ge e^{a_k t} c \qquad \forall t \ge 0 ,
s \in \R \, ,
$$
which establishes (\ref{eqpp2}). As seen earlier,
this completes the proof of (ii)$\Rightarrow$(iii).

Finally, to prove that (iii)$\Rightarrow$(i), assume $0<\beta < 1$
satisfies (iii). Recalling that $a_j/b_j >0$ for every $j\in \{1,
\ldots, m\}$, let
$$
\alpha = \sqrt{ \min\nolimits_{j=1}^m (a_j/b_j) \max\nolimits_{j=1}^m
  (a_j/b_j) } > 0 \, ,
$$
and define $h:\R^m\to \R^m$ as
$$
h(x)_j = (\mbox{\rm sign} \, x_j) |x_j|^{\alpha b_j/a_j} = \left\{
  \begin{array}{rl}
    x_j^{\alpha b_j/a_j} & \mbox{\rm if } x_j \ge 0 \, ,  \\
    -|x_j|^{\alpha b_j/a_j} & \mbox{\rm if } x_j <  0 \, , 
    \end{array}
  \right.
  \qquad \forall x \in \R^m, j\in \{1, \ldots , m\} \, .
$$
Then $h$ is a homeomorphism, and in fact $h\in \cH_{\gamma}(\R^m)$, with
$$
\gamma = \min\nolimits_{j=1}^m \left\{
  \frac{\alpha b_j}{a_j}, \frac{a_j}{\alpha b_j} \right\}
=
\min \left\{ \frac{\alpha}{
      \max_{j=1}^m (a_j/b_j)} , \frac{\min_{j=1}^m (a_j/b_j)}{\alpha}
  \right\} = \sqrt{\frac{\min_{j=1}^m (a_j/b_j)}{\max_{j=1}^m
      (a_j/b_j)}} \ge \beta \, .
$$
Furthermore, observe that
$$
h(\Phi_t x)_j = e^{\alpha b_j t} h(x)_j = \bigl( \Psi_{\alpha t}
h(x)\bigr)_j \qquad \forall t\in \R , x \in \R^m , j \in \{1, \ldots ,
m\} \, ,
$$
that is, $\Phi \stackrel{h}{\simeq} \Psi$ with $\tau_x(t) =\alpha
t$ for all $x\in \R^m$, and hence (i) holds.
\end{proof}

To re-state Theorem \ref{lemr9} concisely, and to extend it slightly, the following
tailor-made terminology is useful: Given two hyperbolic flows $\Phi$,
$\Psi$ on $X$, let
$$
\rho_+ (\Phi, \Psi) = \frac{\min_{j=1}^d (\lambda_j^{\Phi}/ \lambda_j^{\Psi})}{|\max_{j=1}^d
  (\lambda_j^{\Phi}/\lambda_j^{\Psi})|} \le 1 \, ,
$$
and define the {\bf Lyapunov cross ratio} $\rho (\Phi,
\Psi)$ as
$$
  \rho (\Phi, \Psi) = \max \{ \rho_+ (\Phi, \Psi), \rho_+ (\Phi^*, \Psi)\}
  \, .
$$
Clearly, $\rho_+(\Phi,\Psi) \ne 0$, and $\rho_+(\Phi,\Psi)>0$ if and
only if $\lambda_j^{\Phi}/\lambda_j^{\Psi}>0$ for every $j\in \{1,
\ldots , m \}$. Thus, $\rho(\Phi,\Psi)>0$ if and only if $\{d_{\sf
  S}^{\Phi}, d_{\sf U}^{\Phi}\} = \{d_{\sf S}^{\Psi}, d_{\sf
  U}^{\Psi}\}$, and then also $\rho(\Phi,\Psi) =
\rho(\Psi,\Phi)$. Notice, however, that $\rho_+$ is not symmetric,
i.e.,  $\rho_+(\Phi,\Psi) \ne \rho_+ (\Psi,\Phi)$ in general, and
neither is $\rho$. 

When expressed using a Lyapunov cross ratio, Theorem
\ref{lemr9}(iii) simply reads $\beta^2 \le \rho (\Phi, 
\Psi)$, or equivalently $\beta^2 \le \rho (\Psi,\Phi)$. The following
corollary shows that this condition carries over to any two (un)stable
flows, as does the fact, implicit in Theorem \ref{lemr9} or an immediate
consequence thereof, that the relations
$\stackrel{\bigstar}{\simeq}$, $\stackrel{\bigstar}{\thicksim}$
coalesce for such flows for each $\bigstar \in \{\beta^-, \beta^+\}$
with $0<\beta < 1$.

\begin{cor}\label{cor48}
Let $\Phi$, $\Psi$ be stable or unstable flows on $X$. Then, for
every $0<\beta<1$:
\begin{enumerate}
\item $\Phi \stackrel{\beta^-}{\simeq} \Psi \enspace \Longleftrightarrow \enspace
  \Phi \stackrel{\beta^-}{\thicksim} \Psi \enspace \Longleftrightarrow \enspace
  \beta^2 \le \rho (\Phi,\Psi)$;
\item $\Phi \stackrel{\beta^+}{\simeq} \Psi \enspace\Longleftrightarrow \enspace
  \Phi \stackrel{\beta^+}{\thicksim} \Psi \enspace \Longleftrightarrow \enspace
  \beta^2 < \rho (\Phi,\Psi)$.
\end{enumerate}
\end{cor}

\begin{proof}
Since $\Phi \stackrel{{\sf lin}}{\simeq} \Phi^*$ and clearly $\rho(\Phi^*,
\Psi) = \rho(\Phi,\Psi) $, it can be assumed that $\Phi$, $\Psi$ 
either both are stable or else both are unstable, and in either case $\rho(\Phi, \Psi) = \displaystyle
\frac{\min_{j=1}^d (\lambda_j^{\Phi}/\lambda_j^{\Psi})}{\max_{j=1}^d
  (\lambda_j^{\Phi}/\lambda_j^{\Psi})}$.

To prove (i), let $\widetilde{\Phi}$, $\widetilde{\Psi}$ be generated by
$\Lambda^{\Phi}$, $\Lambda^{\Psi}$ respectively. Recall that $\Phi
\stackrel{1^-}{\simeq} \widetilde{\Phi}$ and $\Lambda^{\Phi} =
\Lambda^{\widetilde{\Phi}}$, and similarly for $\Psi$. With this, for every
$0<\beta < 1$,
$$
\Phi \stackrel{\beta^-}{\thicksim} \Psi \enspace \Longrightarrow \enspace 
\widetilde{\Phi} \stackrel{\beta^-}{\thicksim} \widetilde{\Psi}
\enspace \Longrightarrow \enspace \beta^2 \le \rho \bigl( \widetilde{\Phi},
\widetilde{\Psi} \bigr) = \rho (\Phi, \Psi) \enspace \Longrightarrow
\enspace \widetilde{\Phi} \stackrel{\beta^-}{\simeq} \widetilde{\Psi}
\enspace \Longrightarrow
\enspace  \Phi \stackrel{\beta^-}{\simeq} \Psi \, ;
$$
here the first and fourth implications are due to the fact that $h_1
\circ h_2 \circ h_3 \in \cH_{\beta^-}$ whenever $h_1,h_3\in \cH_{1^-}$
and $h_2 \in \cH_{\beta^-}$, while the second and third implications
are due to Theorem \ref{lemr9}.

To prove (ii), assume first that $\Phi \stackrel{\beta^+}{\thicksim}
\Psi $. Then $\Phi \stackrel{\gamma}{\thicksim}
\Psi $ for some $\beta < \gamma < 1$, so $\beta^2 < \gamma^2
\le \rho (\Phi,\Psi)$, by (i). Conversely, assume that $\beta^2 < \rho
(\Phi,\Psi)$, and pick $\beta <\gamma < 1$ so that $\gamma^2 \le \rho
(\Phi,\Psi)$. By (i), $\Phi \stackrel{\gamma^-}{\simeq} \Psi $, and
hence $\Phi \stackrel{\beta^+}{\simeq} \Psi$ as well. 
\end{proof}

\begin{rem}\label{remr13}
  (i) In the context of Corollary \ref{cor48}, notice that if $\Phi$,
  $\Psi$ are irreducible then $\rho (\Phi,\Psi)=1$. This explains why
  $\Phi \stackrel{1^-}{\simeq} \Psi$ is automatic for
  irreducible (un)stable flows, as seen in Proposition \ref{propr4A}.

  (ii) The same strategy as in the proof of Theorem
  \ref{lemr9} can be utilized to establish a characterization of
  $\Phi\stackrel{\beta}{\thicksim} \Psi$ for $0<\beta <  1$. Since
  Corollary \ref {cor48} suffices for the purpose of the present
article, this topic will not be pursued further here. Notice, however, that unlike
for $\stackrel{\bigstar}{\thicksim}$ with $\bigstar \in \{\beta^- ,
  \beta^+\}$, such a characterization must depend on finer geometric
  properties of $\Phi$, $\Psi$, not merely on their Lyapunov
  exponents. For a simple example illustrating this, consider the
  flows $\Phi$, $\Psi$ on $\R^5$ generated by $\mbox{\rm 
    diag} \, [1,1, J_2(2), 4]$, $I_5$ respectively, for which $\Phi
  \stackrel{0.5}{\simeq}  \Psi$ but $\Phi\cancel{
  \stackrel{0.5^+}{\thicksim}} \Psi$; with $\widetilde{\Phi}$ generated
  by $\mbox{\rm diag} \, [J_2(1), 2, 2, 4]$, by contrast,
  $\widetilde{\Phi} \! \stackrel{0.5^-}{\simeq} \! \Psi$ but
  $\widetilde{\Phi} \, \cancel{\stackrel{0.5}{\thicksim} } \, \Psi$,
  notwithstanding the fact that $\Lambda^{\Phi} =
  \Lambda^{\widetilde{\Phi}}$.
\end{rem}

The calculation proving (iii)$\Rightarrow$(i) in Theorem
\ref{lemr9} can be used to extend one part of 
Corollary \ref{cor48}(i) to hyperbolic flows; the straightforward
details are left to the interested reader.

\begin{cor}\label{cor48A}
Let $\Phi$, $\Psi$ be hyperbolic flows on $X$ and $0<\beta < 1$. If
$\beta^2 \le \rho (\Phi,\Psi)$ then $\Phi \stackrel{\beta^-}{\simeq} \Psi$.
\end{cor}

The authors conjecture that the converse of Corollary \ref{cor48A} is
true also, even in a stronger form. More precisely, they conjecture
that Corollary \ref{cor48} remains valid with {\em stable or
  unstable\/} replaced by {\em hyperbolic}, and hence in particular
that the relations $\stackrel{\bigstar}{\simeq}$, $\stackrel{\bigstar}{\thicksim}$
coalesce for hyperbolic flows and each $\bigstar \in \{\beta^-,
\beta^+\}$ with $0<\beta < 1$, just as they do in the case of (un)stable flows.

%%%%%%%%%%%%%%%%%%%%%%%%%%%%%%%%%%%%%%%%%%%%%%%%%%%%%%%%%%%%%%%%

\section{Preserving minimal periods}\label{sec4}

%%%%%%%%%%%%%%%%%%%%%%%%%%%%%%%%%%%%%%%%%%%%%%%%%%%%%%%%%%%%%%%%

This short section presents a simple observation regarding
minimal periods in linear flows. Though solely of an auxiliary nature in this
article, the result may also be of independent interest. In preparation
for the statement and proof, for any (not necessarily linear) flow $\varphi$ on
$X$, denote the minmal $\varphi$-period of $x\in X$ by $T_x^{\varphi} = \inf
\{t\in \R^+: \varphi_t(x) = x\}$, with the usual convention that $\inf
\varnothing = \infty$. Thus $x\in \mbox{\rm Fix}\, \varphi$ if and only if $T_x^{\varphi} =
0$. If $T_x^{\varphi} \in \R^+$ then $x$ is $T$-periodic with $T\in
\R^+$, i.e., $\varphi_T(x) = x$, precisely if $T = nT_x^{\varphi}$ for
some $n\in \N$. For convenience, let $\mbox{\rm Per}_T\varphi = \{x\in
X: \varphi_T(x) = x\}$ for every $T\in \R^+$, and let $\mbox{\rm
  Per}\, \varphi = \bigcup_{T\in \R^+} \mbox{\rm Per}_T \varphi$.

The following is a characterization of a certain rigidity for minimal
periods in linear flows where all orbits are bounded.

\begin{lem}\label{lem2a}
Given $m_0,n_0\in \N_0$ and $m, n \in \N$ with $m_0+2m = n_0 + 2n=d$, as
well as $b\in  (\R^+)^m$, $c \in (\R^+)^n$, let $\Phi, \Psi$ be
the flows on $\R^{d}$ generated by
$$
\mbox{\rm diag}\,  [O_{m_0},  J_1 (ib_1), \ldots , J_1(ib_m) ] \, , \qquad
\mbox{\rm diag}\,  [O_{n_0}, J_1 (ic_1), \ldots , J_1(ic_n) ] \, ,
$$
with $O_{m_0}$, $O_{n_0}$ understood to be present
only if $m_0\ge 1$, $n_0\ge 1$, respectively. Then the following are equivalent:
\begin{enumerate}
\item there exists an open set $U\subset \R^{d}$ with $0\in U$
  and a continuous one-to-one function $f:U\to \R^{d}$ with
  $T_x^{\Phi} = T_{f(x)}^{\Psi}$ for every $x\in U$;
\item $(m_0,m)=(n_0,n)$, and there exists a permutation $p$ of $\{1, \ldots , m\}$ so that
  $b_j = c_{p(j)}$ for every $j\in \{1,\ldots , m\}$;
\item $\Phi \stackrel{{\sf lin}}{\cong}\Psi$, i.e., $\Phi$, $\Psi$ are
  linearly conjugate.
\end{enumerate}
\end{lem}

\begin{proof}
Note first that $\mbox{\rm Fix}\, \Phi = \mbox{\rm span} \{e_{1},
\ldots , e_{m_0}\}$ and $\mbox{\rm Fix}\, \Psi = \mbox{\rm span} \{e_{1},
\ldots , e_{n_0}\}$, with $\mbox{\rm span}\, \varnothing = \{0\}$ as
usual. If (i) holds then $f(0)\in \mbox{\rm Fix}\, \Psi$, and
replacing $f$ by $f - f(0)$ otherwise, it can be assumed that
$f(0)=0$. Moreover, $f(U\cap \mbox{\rm Fix}\, \Phi) = f(U)\cap
\mbox{\rm Fix}\, \Psi$, so Proposition \ref{prop4top} yields $m_0 =
\dim \mbox{\rm Fix}\, \Phi \le \dim \mbox{\rm Fix}\, \Psi =
n_0$. Since $f(0)=0$, one may interchange the roles of $\Phi$, $\Psi$,
which yields $m_0 = n_0$. Trivially, $m_0=n_0$ also if either (ii) or
(iii) holds. In other words, if $m_0\ne n_0$ then (i), (ii), and (iii)
all are false, so to prove the lemma it suffices to consider the case
of $m_0=n_0$. Thus, assume $(m_0,m)=(n_0,n)$ from now on. 
Furthermore, for every $j\in \{1,\ldots , m\}$ let
$E_j = \mbox{\rm span} \{e_{ m_0+ 2j -1}, e_{m_0+ 2j}\}$ and denote by
$P_j$ the orthogonal projection of $\R^d$ onto $E_j$; also let $E_0=
\mbox{\rm span} \{ e_{1}, \ldots , e_{m_0}\} = \mbox{\rm Fix}\, \Phi =
\mbox{\rm Fix}\, \Psi$ and denote by $P_0= I_d - \sum_{j=1}^m P_j$ the orthogonal
projection of $\R^d$ onto $E_0$.

Notice first that for $m=1$ the asserted equivalence of (i), (ii), and
(iii) trivially is correct if
$m_0=0$, and if $m_0 \ge 1$ then, for every $T\in \R^+$,
$$
U \cap \{T^{\Phi} = T\} := \{x\in U : T_x^{\Phi} = T\} = \left\{
  \begin{array}{ll} U \setminus E_0 & \mbox{\rm if } T = 2\pi/b_1 \, , \\
    \varnothing & \mbox{\rm otherwise}\, ,
  \end{array}                
  \right.
$$
and similarly with $\Phi$, $b_1$ replaced by $\Psi$, $c_1$
respectively. Thus for $m=1$ and $m_0 \ge 1$, (i) implies that $f(U
\cap \{T^{\Phi} = T\}) = f(U) \cap \{T^{\Psi} = T\}\ne \varnothing$
precisely if $2\pi/b_1 = T = 2\pi/c_1$, and hence $b_1 =
c_1$. Also, if $b_1 = c_1$ then $\Phi \stackrel{h}{\cong}\Psi$ with
$h=I_d$, so $\Phi \stackrel{{\sf
    lin}}{\cong}\Psi$. Finally, if $h \Phi_t = \Psi_t h$ for some
$h\in \cH_{\sf lin}(\R^d)$ and all $t\in \R$, then $U= \R^d$ and 
$f= I_d$ obviously satisfy (i). In 
summary, (i), (ii), and (iii) are equivalent
whenever $m=1$, so henceforth assume $m\ge 2$. 

Given $m_0 \in \N_0$ and $m\in \N \setminus \{1\}$, assume (i), and w.l.o.g.\
let $b_1 \le \ldots \le b_m$ and $c_1 \le \ldots \le c_m$. To prove
that (i)$\Rightarrow$(ii), it suffices to show that 
\begin{equation}\label{eqp1}
b_j = c_j \qquad \forall j \in \{1, \ldots , m \} \, .
\end{equation}
To prove (\ref{eqp1}) fix any $T\in \R^+$, and note that $\mbox{\rm Per}_T
\Phi$, $\mbox{\rm Per}_T \Psi$ are subspaces of $\R^d$, in fact
$$
\mbox{\rm Per}_T \Phi = E_0 \oplus \bigoplus_{j: b_j T \in 2\pi \N}
E_j \, , \qquad
\mbox{\rm Per}_T \Psi = E_0 \oplus \bigoplus_{j: c_j T \in 2\pi \N}
E_j \, .
$$
Moreover, observe that
$$
\{T^{\Phi} = T\} = \bigcap_{\ell \in \N \setminus \{1\}} \bigl( \mbox{\rm
  Per}_T \Phi \setminus   \mbox{\rm Per}_{T/\ell} \Phi \bigr) \, .
$$
Since $\mbox{\rm Per}_{T/\ell} \Phi$ is a subspace of $\mbox{\rm
  Per}_T \Phi$, and since $\mbox{\rm Per}_{T/\ell} \Phi = E_0$ for all
sufficiently large $\ell$, if the set $\{ T^{\Phi} =
T\}$ is non-empty then it is open and dense in $\mbox{\rm Per}_T
\Phi$. Similarly, $\{T^{\Psi} = T\}$ is open and dense in $\mbox{\rm
  Per}_T \Psi$ whenever non-empty. By assumption, $f(U\cap \{T^{\Phi}
= T\}) = f(U) \cap \{T^{\Psi} = T\}$. Since $f:U \to f(U)$ is continuous
and one-to-one, $f(U\cap \mbox{\rm Per}_T \Phi) = f(U) \cap \mbox{\rm Per}_T
\Psi$, and Proposition \ref{prop4top} yields
$$
 m_0 + 2\# \{1\le j \le m : b_j T \in 2\pi \N\}  = \dim \mbox{\rm Per}_T
\Phi \le \dim \mbox{\rm Per}_T \Psi = m_0 + 
 2\# \{1\le j \le m : c_j T \in 2\pi \N\}  \, .
$$
Interchanging the roles of $\Phi$, $\Psi$ yields, 
since $T\in \R^+$ has been arbitrary,
\begin{equation}\label{eqp2}
 \# \{1\le j \le m : b_j \in  s \N\} = \# \{1\le j \le m : c_j \in  s
 \N\} \qquad \forall s \in \R^+ \, .
\end{equation}
Utilizing (\ref{eqp2}), the desired conclusion (\ref{eqp1}) is now easily
obtained as follows: First observe that if $b_m<c_m$ then the integer
on the left in (\ref{eqp2}) for $s=c_m$ would be zero, whereas the
integer on the right would be positive, an obvious 
contradiction. Similarly, $b_m>c_m$ is impossible, and hence $b_m =
c_m$. Taking $s=b_m=c_m$ in (\ref{eqp2}) yields
$$
\# \{1\le j \le m : b_j = b_m \} = \# \{1\le j \le m : c_j = c_m\} \, . 
$$
Let $j_1 = \min \{ 1\le j \le m : b_j = b_m\}-1$. If $j_1 = 0$ then
$b_1 = b_m$ and $c_1 = c_m$, and hence (\ref{eqp1}) holds; otherwise, clearly $b_{j_1}< b_{j_1 +
  1} = \ldots = b_m$ and $c_{j_1} < c_{j_1 + 1} = \ldots c_m$. By
interchanging the roles of $\Phi, \Psi$ if necessary, assume w.l.o.g.\
that $b_{j_1}\ge c_{j_1}$. Notice that if $c_j = b_{j_1} \ell$ for
some $\ell \in \N \setminus \{1\}$ then $j \ge j_1 + 1$ and hence $b_j
= c_j$. This, together with (\ref{eqp2}) for $s= b_{j_1}$, yields
\begin{align*}
\# \{1\le j \le m : b_j = b_{j_1}\} & = \# \{1\le j \le m : b_j \in
                                      b_{j_1}\N \}  -
                                      \sum\nolimits_{\ell =2}^{\infty}
                                      \# \{1\le j \le m : b_j = \ell
                                      b_{j_1}\} \\
  & = \# \{1\le j \le m : c_j \in
                                      b_{j_1}\N \}  -
                                      \sum\nolimits_{\ell =2}^{\infty}
                                      \# \{1\le j \le m : c_j = \ell
    b_{j_1}\} \\
  & = \# \{1\le j \le m : c_j = b_{j_1}\} \, ,
\end{align*}
and so in particular $b_{j_1} = c_{j_1}$. Repeating this argument with
$$
j_{k+1} = \min \{ 1\le j \le m : b_j = b_{j_k}\} - 1 
$$
instead of $j_1$ yields
$$
\# \{1\le j \le m : b_j = b_{j_{k+1}}\} = \# \{1\le j \le m : c_j = b_{j_{k+1}}\}
$$
for every $k$, and hence also $b_{j_{k+1}} = c_{j_{k+1}}$. Since
$j_{k+1} \le j_k - 1\le m - (k+1)$, this procedure terminates after at
most $m$ steps, i.e., $j_k=0$ for some $1\le k \le m$, which in turn
establishes (\ref{eqp1}). As seen earlier, this proves that
(i)$\Rightarrow$(ii).

Next, given $m_0\in \N_0$ and $m\ge 2$, to prove that
(ii)$\Rightarrow$(iii), assume $b_j = c_{p(j)}$ for some permutation $p$
of $\{1, \ldots , m\}$ and every $j\in \{1, \ldots , m\}$. Defining a
linear map as
$$
f(x) = P_0 x + \sum\nolimits_{j=1}^m (x_{m_0 + 2j - 1} e_{m_0 + 2p(j)
  - 1} + x_{m_0 + 2j} e_{m_0 + 2p(j)}) \qquad \forall x \in \R^d \, ,
$$
note that $f:\R^d \to \R^d$ is an isomorphism with $\Psi_t f = f
\Phi_t$ for all $t\in \R$. Thus $\Phi \stackrel{\sf lin}{\cong} \Psi$,
i.e., (iii) holds.

Finally, to prove that (iii)$\Rightarrow$(i), assume $\Phi
\stackrel{h}{\cong}\Psi$ for some $h\in \cH_{\sf lin}$. Then $h\Phi_t
x = \Psi_t h x$ for all $t\in \R$, $x\in \R^d$, so clearly 
$T_x^{\Phi} = T_{hx}^{\Psi}$ for every $x\in \R^d$. In other words,
(i) holds with $U=\R^d$ and $f=h$.
\end{proof}

In an appropriately adjusted form, Lemma \ref{lem2a} extends to
all linear flows. 

\begin{theorem}\label{thm2b} 
Let $\Phi, \Psi$ be linear flows on $X$. Assume there exists an open
set $U\subset X$ with $0\in U$ and a continuous one-to-one
function $f:U\to X$ with $T_x^{\Phi} = T_{f(x)}^{\Psi}$ for every
$x\in U$. Then $\sigma (\Phi)\cap i\R = \sigma (\Psi) \cap i\R$.
\end{theorem}

\begin{proof}
Notice first that $\sigma (\Phi) \cap i \R = \varnothing$ if and only
if $T_x^{\Phi} = \infty$ for every $x\in U \setminus \{0\}$. By
assumption, $f(U\cap \{T^{\Phi} <\infty \}) = f(U) \cap \{T^{\Psi}
<\infty \}$,
and so $\sigma (\Phi) \cap i\R = \varnothing$ precisely if $\sigma
(\Psi) \cap i\R  =\varnothing$, in which case the assertion trivially
is correct. Thus, assume henceforth that $\sigma (\Phi) \cap
i\R\ne \varnothing$ and $\sigma (\Psi) \cap
i\R\ne \varnothing$.

For convenience, write $A^{\Phi}$, $A^{\Psi}$ as $A$, $B$
respectively. From
$$
U\cap \{T^{\Phi} =0\} = U \cap \mbox{\rm Fix}\, \Phi = U\cap \ker A \,
, \qquad
f(U)\cap \{T^{\Psi} =0\} = f(U) \cap \mbox{\rm Fix}\, \Psi = f(U) \cap \ker B \, , 
$$
and since $f(U \cap \{T^{\Phi} = 0\}) = f(U)\cap \{T^{\Psi} = 0\}$ by
assumption, it is clear that $0\in \sigma (\Phi)$ if and only if $0\in
\sigma (\Psi)$, and it follows from Proposition \ref{prop4top} that $\dim \ker
A = \dim \ker B$ always; for convenience, denote the latter number by
$m_0 \in \N_0$. It will be shown below that
\begin{equation}\label{eqp3}
\sigma (\Phi) \cap is \Q = \sigma (\Psi) \cap is \Q \qquad \forall s
\in \R^+ \, .
\end{equation}
Notice that (\ref{eqp3}) immediately proves the assertion of the
theorem, since
$$
\sigma (\Phi)\cap i\R = \bigcup\nolimits_{s\in \R^+} (\sigma (\Phi)
\cap i s \Q) =
\bigcup\nolimits_{s\in \R^+} (\sigma (\Psi)
\cap i s \Q) =
\sigma (\Psi)\cap i\R \, .
$$
As a first step towards establishing (\ref{eqp3}), fix any $T\in \R^+$
and consider the set
$$
X_T:= \{T^{\Phi} \in T \Q^+\} \cup \{T^{\Phi} = 0\} = \{x\in X :
T_x^{\Phi} = 0 \: \mbox{\rm or $T_x^{\Phi} = rT$ for some} \: r\in
\Q^+\} \, .
$$
Deduce from
$$
X_T = \ker A \oplus \bigoplus_{r\in \Q^+} \ker \left ( A^2 +
  \frac{4\pi^2}{r^2 T^2} I_d\right)
$$
that $X_T$ is a $\Phi$-invariant subspace with $X_T \supset \ker A$,
and the spectrum of the restricted flow $\Phi|_{\R \times X_T}$ is given by
$$
\sigma (\Phi|_{\R \times X_T}) = \sigma (\Phi) \cap \frac{2\pi i}{T} \Q \, .
$$
Similarly,
$$
Y_T:= \{T^{\Psi} \in T \Q^+\} \cup \{T^{\Psi} = 0\} = \ker B \oplus \bigoplus_{r\in \Q^+} \ker \left ( B^2 +
  \frac{4\pi^2}{r^2 T^2} I_d \right)
$$
is a $\Psi$-invariant subspace, and
$$
\sigma (\Psi|_{\R \times Y_T}) = \sigma (\Psi) \cap \frac{2\pi i}{T} \Q \, .
$$
By assumption, $f(U \cap X_T) = f(U)\cap Y_T$, and hence in particular
$\dim X_T = \dim Y_T\ge m_0$.

Now, consider first the case where $\dim X_T = \dim Y_T = m_0$, or
equivalently where $X_T = \ker A$ and $Y_T = \ker B$. In this case,
$$
\sigma (\Phi|_{\R \times X_T}) = \left\{
  \begin{array}{cl}
    \varnothing & \mbox{\rm if } m_0 = 0 \\
    \{0\} & \mbox{\rm if } m_0 \ge 1 
  \end{array}
    \right\} = \sigma (\Psi|_{\R \times Y_T}) \, ,
 $$
 and consequently
 \begin{equation}\label{eqp4}
\sigma (\Phi) \cap \frac{2\pi i }{T} \Q = \sigma (\Psi) \cap
\frac{2\pi i }{T} \Q \, .
 \end{equation}

 Next, consider the case where $\dim X_T = \dim Y_T > m_0$, and hence
 $\dim X_T = m_0 + 2m$ for some $m\in \N$. With appropriate
 $b_j,c_j \in 2\pi T^{-1}\Q^+$ for all $j\in \{1, \ldots , m\}$,
 $$
\sigma (\Phi|_{\R \times X_T}) \setminus \{0\} = \{\pm i b_j : 1\le j \le m\} \,
, \qquad
\sigma (\Psi|_{\R \times Y_T}) \setminus \{0\} = \{\pm i c_j : 1\le j \le m\} \, .
$$
For convenience, let $\widetilde{\Phi}$, $\widetilde{\Psi}$ be
the flows on $\R^{m_0+2m}$ generated by 
$$
\mbox{\rm diag}\,  [O_{m_0},  J_1 (ib_1), \ldots , J_1(ib_m) ] \, , \qquad
\mbox{\rm diag}\,  [O_{m_0}, J_1 (ic_1), \ldots , J_1(ic_m) ] \, ,
$$
respectively, and pick isomorphisms $H_X: X_T \to \R^{m_0+2m}$ and
$H_Y:Y_T \to \R^{m_0+2m}$ so that
$$
H_X \Phi_t x = \widetilde{\Phi}_t H_X x \, , \quad H_Y \Psi_t y =
\widetilde{\Psi}_t H_Y y \qquad \forall t\in \R, x \in X_T, y \in Y_T
\, .
$$
The set $\widetilde{U}:= H_X U \subset \R^{m_0+2m}$ is open
with $0\in \widetilde{U}$, and the function $\widetilde{f}:
\widetilde{U}\to \R^{m_0+2m}$ given by
$$
\widetilde{f} (z) = H_Y f (H_X^{-1} z) \qquad \forall z \in
\widetilde{U} \, ,
$$
is continuous and one-to-one. Moreover,
$$
T_z^{\widetilde{\Phi}} = T_{H_X^{-1}z}^{\Phi} =
T_{f(H_X^{-1}z)}^{\Psi} = T_{\widetilde{f}(z)}^{\widetilde{\Psi}}
\qquad \forall z \in \widetilde{U} \, ,
$$
so by Lemma \ref{lem2a}, $b_j = c_{p(j)}$ for
some permutation $p$ of $\{1, \ldots , m\}$ and every $j\in \{1,
\ldots , m\}$. In particular,
therefore,
$$
\sigma (\Phi) \cap \left( \frac{2\pi i}{T} \Q \setminus \{0\}\right)  = \{\pm i b_j:
1\le j \le m \} =  \{\pm i c_j:
1\le j \le m \} = \sigma (\Psi) \cap \left( \frac{2\pi i}{T} \Q \setminus
\{0\} \right) \, .
$$
Moreover, as seen above,
$$
\sigma (\Phi )\cap \{0\} = \left\{
  \begin{array}{cl}
    \varnothing & \mbox{\rm if } m_0 = 0 \\
    \{0\} & \mbox{\rm if } m_0 \ge 1 
  \end{array}
    \right\} = \sigma (\Psi)\cap \{0\} \, ,
$$
and hence (\ref{eqp4}) holds in this case also. In summary,
(\ref{eqp4}) is valid for {\em every\/} $T\in \R^+$. This clearly
establishes (\ref{eqp3}), and the latter in turn proves the theorem,
as discussed previously. 
\end{proof}

\begin{rem}
The non-trivial implications (i)$\Rightarrow \! \! \!
\{$(ii),(iii)$\}$ in Lemma \ref{lem2a}, and hence the 
conclusion in Theorem \ref{thm2b} as well, may fail if $U$
is not open or $0\not \in U$, and also if $f$ is not continuous or not
one-to-one. In other words, every single assumption on $U$ and $f$ in
these results is indispensable. For a simple example illustrating this,
using the same symbols as in Lemma \ref{lem2a} and its proof, take
$(m_0,m) = (n_0,n) =(0,3)$ as well as $b=(2,2,3)\in (\R^+)^3$ and $c=(1,2,3)\in
(\R^+)^3$, so clearly
$$
\sigma (\Phi) \cap i\R = \sigma (\Phi) = \{\pm 2i , \pm 3i\} \ne \{\pm
i , \pm 2i , \pm 3i\} = \sigma (\Psi) = \sigma (\Psi) \cap i\R \, .
$$
Here Lemma \ref{lem2a}(ii),(iii) and the conclusion in Theorem \ref{thm2b}
all fail. 
However, take $f:U\to \R^6$ as $f(x)=x$ for all $x\in U$, so $f$ is
continuous and one-to-one; moreover, $T_x^{\Phi}= T_{f(x)}^{\Psi}$ for every
$x\in U$ if, for instance, $U=E_3$ (containing $0$, but not open) or
$U= B_{1/2} (e_1+e_3+e_5)$ (open, but $0\not \in U$). Similarly, take
$f:\R^6 \to \R^6$ to be any (necessarily discontinuous) bijection with
$$
f(0) = 0 \, , \quad f(E_1 \oplus E_2 \setminus \{0\}) = E_2 \setminus \{0\} \, , \qquad
f(E_3 \setminus \{0\}) = E_3 \setminus \{0\} \, .
$$
Then $T_x^{\Phi} = T_{f(x)}^{\Psi}$ for every $x\in U=\R^6$. Finally,
take $f:\R^6 \to \R^6$ as $f(x) = |P_1 x + P_2 x|e_3 + P_3 x$, so $f$
is continuous but not one-to-one, and again $T_x^{\Phi} =
T_{f(x)}^{\Psi}$ for every $x\in U=\R^6$.  
\end{rem}

%%%%%%%%%%%%%%%%%%%%%%%%%%%%%%%%%%%%%%%%%%%%%%%%%%%%%%%%%%%%%%%%

\section{H\"{o}lder classifications}\label{sec5}

%%%%%%%%%%%%%%%%%%%%%%%%%%%%%%%%%%%%%%%%%%%%%%%%%%%%%%%%%%%%%%%%

Utilizing the tools developed above, this section proves
the main results previewed already in the Introduction, namely the some-H\"{o}lder
and all-H\"{o}lder classifications of linear flows on $X$.
As it turns out, {\em some}-H\"{o}lder equivalence is easily
characterized. The following is a mildly extended form of Theorem \ref{thma}.

\begin{theorem}\label{thm6some}
Let $\Phi$, $\Psi$ be linear flows on $X$. Then the following
statements are equivalent:
\begin{enumerate}
\item $\Phi \stackrel{0^+}{\simeq}\Psi$;
\item $\Phi \stackrel{0^+}{\thicksim}\Psi$;
\item $\Phi \stackrel{0}{\simeq}\Psi$;
\item $\Phi \stackrel{0}{\thicksim} \Psi$;
  \item $\{d_{\sf S}^{\Phi}, d_{\sf U}^{\Phi}\} = \{d_{\sf S}^{\Psi},
    d_{\sf U}^{\Psi}\} $, and there exists an $\alpha \in \R \setminus
    \{0\}$ so that $A^{\Phi_{\sf C}}$, $\alpha A^{\Psi_{\sf C}}$ are
    similar. 
\end{enumerate}
\end{theorem}

A simple preparatory observation is helpful for proving
(v)$\Rightarrow$(i) in Theorem \ref{thm6some}.

\begin{lem}\label{lem61}
Let $\Phi$, $\Psi$ be linear flows on $X$ and $0<\beta < 1$. Assume that $A^{\Phi_{\sf
        C}}$, $\alpha A^{\Psi_{\sf C}}$ are similar for some $\alpha \in \R ^+$, and 
    \begin{equation}\label{eq61}
\beta \le \min \left\{ \frac{\alpha \lambda_j^{\Psi_{\sf
        H}}}{\lambda_j^{\Phi_{\sf H}}},  \frac{  \lambda_j^{\Phi_{\sf
        H}}}{\alpha \lambda_j^{\Psi_{\sf H}}}\right\} \qquad \forall j
\in \{1, \ldots , d_{\sf H}^{\Phi}\} \, .
\end{equation}
Then $\Phi \stackrel{\beta^-}{\simeq} \Psi$. 
\end{lem}

\begin{proof}
 By the similarity of $A^{\Phi_{\sf
        C}}$, $\alpha A^{\Psi_{\sf C}}$, clearly $d_{\sf C}^{\Phi} =
    d_{\sf C}^{\Psi} =:\ell$. If $\ell = 0$ then $\Phi$, $\Psi$ are
    hyperbolic, $\alpha \in \R^+$ is
    arbitrary, and (\ref{eq61}) yields $\beta \le \sqrt{\rho_+
      (\Phi,\Psi)}$. Thus the conclusion follows directly from
    Corollary \ref{cor48A}. Henceforth assume $\ell \ge 1$, and let
    $Q\in \R^{\ell \times \ell}$ be invertible with $Q A^{\Phi_{\sf
        C}} = \alpha A^{\Psi_{\sf C}}Q$. If $\ell = d$ then
    (\ref{eq61}) is void, and the assertion is correct, since in fact $\Phi
    = \Phi_{\sf C} \stackrel{{\sf lin}}{\simeq} \Psi_{\sf C} =
    \Psi$. Thus assume $1\le \ell \le d-1$ from now on, and
    consequently $d_{\sf
      H}^{\Phi} = d_{\sf S}^{\Phi} + d_{\sf U}^{\Phi} = d_{\sf
      S}^{\Psi} + d_{\sf U}^{\Psi} = d -
    \ell \in \{1, \ldots , d-1\}$. Note that if $d_{\sf S}^{\Phi}<
    d_{\sf S}^{\Psi}$
    then $\lambda_j^{\Psi_{\sf H}}/ \lambda_j^{\Phi_{\sf H}} <0$ for
    $j = d_{\sf S}^{\Phi} + 1 \le d_{\sf S}^{\Psi}$, contradicting
    (\ref{eq61}). Similarly, $d_{\sf S}^{\Phi}>d_{\sf S}^{\Psi}$ is impossible, and
    hence $(d_{\sf S}^{\Phi}, d_{\sf U}^{\Phi}) = (d_{\sf S}^{\Psi},
    d_{\sf U}^{\Psi})=: (k,m)$. As in
    Section \ref{sec3}, Lemmas \ref{lemr1} and \ref{lemr4}, applied
    individually to each irreducible component associated with an
    eigenvalue (pair) not in $i\R$, together with Proposition
    \ref{prop3add0}, show that $\Phi \stackrel{1^-}{\simeq}
    \widetilde{\Phi}$, with $\widetilde{\Phi}$ generated by
    $\mbox{\rm diag}\, [\Lambda^{\Phi_{\sf S}}, A^{\Phi_{\sf C}},
    \Lambda^{\Phi_{\sf U}}]$; similarly $\Psi \stackrel{1^-}{\simeq}
    \widetilde{\Psi}$, with $\widetilde{\Psi}$ generated by
    $\mbox{\rm diag}\, [\Lambda^{\Psi_{\sf S}}, A^{\Psi_{\sf C}},
    \Lambda^{\Psi_{\sf U}}]$. In analogy to the proof of Theorem \ref{lemr9},
    define $h:\R^d\to \R^d$ as
    $$
h(x)_j = (\mbox{\rm sign}\,  x_j) |x_j|^{\alpha
  \lambda^{\Psi}_j/\lambda^{\Phi}_j} \qquad \forall x \in \R^d,  j \in
\{1, \ldots , k\} \cup \{k+\ell + 1, \ldots d\} \, , 
$$
whereas, with $P_{\sf C}$ denoting the orthogonal projection of $\R^d$
onto $\mbox{\rm span}  \{e_{k+1}, \ldots , e_{k+\ell}\}= X_{\sf
  C}^{\widetilde{\Phi}} = X_{\sf C}^{\widetilde{\Psi}}$,
$$
P_{\sf C} h(x) = Q P_{\sf C} x \qquad \forall x \in \R^d \, .
$$
Then $h\in \cH_{\beta}(\R^d)$ by
Proposition \ref{prop3add0} and (\ref{eq61}); furthermore, for every
$t\in \R$, $x\in \R^d$, 
$$
h\bigl( \widetilde{\Phi}_t x\bigr) = h \left(
  \mbox{\rm diag} \left[
e^{t \Lambda^{\Phi_{\sf S}}}, e^{tA^{\Phi_{\sf C}}}, e^{t
  \Lambda^{\Phi_{\sf U}}}
    \right] x
  \right) =\mbox{\rm diag} \left[
e^{\alpha t \Lambda^{\Psi_{\sf S}}}, Q e^{tA^{\Phi_{\sf C}}}Q^{-1},
e^{\alpha t \Lambda^{\Psi_{\sf U}}}
    \right] h(x) = \widetilde{\Psi}_{\alpha t} h(x) \, , 
$$
where the last equality is due to $Q  A^{\Phi_{\sf C}} Q^{-1} =
\alpha A^{\Psi_{\sf C}}$. In other words, $\widetilde{\Phi}
\stackrel{h}{\simeq} \widetilde{\Psi}$ with $h\in \cH_{\beta}$,
and hence $\Phi \stackrel{\beta^-}{\simeq} \Psi$ as claimed.    
\end{proof}

Given any $\alpha \in \R^+$, and assuming $d_{\sf H}^{\Phi} = d_{\sf
  H}^{\Psi}$, note that the RHS of the inequality in
(\ref{eq61}) is positive for all $j$, and hence (\ref{eq61}) holds for
{\em some\/} $\beta >0$, precisely if $(d_{\sf S}^{\Phi}, d_{\sf
  U}^{\Phi}) =  ( d_{\sf S}^{\Psi} , d_{\sf  U}^{\Psi})$. If so, and
if in addition $\sigma (\Phi_{\sf C}) = \sigma (\Psi_{\sf C}) \subset 
\{0\}$, i.e., if the flows $\Phi$, $\Psi$ either are hyperbolic or
else have $0$ as their only eigenvalue on the imaginary axis, then
$\alpha \in \R ^+$ is arbitrary, and hence (\ref{eq61}) can be
optimized over $\alpha$; this results in it taking the form $\beta \le
\sqrt{\rho_+ (\Phi_{\sf H}, \Psi_{\sf H})}$, which is consistent with
Corollary \ref{cor48A}. 

\begin{proof}[Proof of Theorem \ref{thm6some}]
  Obviously (i)$\Rightarrow\! \! \{$(ii),(iii)$\}\!\! \Rightarrow$(iv) by definition.

  To prove that (iv)$\Rightarrow$(v), assume $\Phi
  \stackrel{0}{\thicksim} \Psi$, and so $\Phi
  \stackrel{0}{\thickapprox} \Psi$ by Proposition \ref{lemH1}. On the
  one hand, if $\tau_x$ is increasing for every $x\in X\setminus \{ 0 \}$, then
  \cite[Thm.1.1]{BW} yields $(d_{\sf S}^{\Phi}, d_{\sf U}^{\Phi})
  = (d_{\sf S}^{\Psi}, d_{\sf U}^{\Psi})$, and $A^{\Phi_{\sf C}}$, $\gamma A^{\Psi_{\sf C}}$ are
  similar for some $\gamma \in \R^+$. On the other hand, if
  $\tau_x$ is decreasing for every $x$, then the same argument
  with $\Phi$ replaced by $\Phi^*$ yields $(d_{\sf U}^{\Phi}, d_{\sf S}^{\Phi}) =
  (d_{\sf S}^{\Phi^*}, d_{\sf U}^{\Phi^*}) = (d_{\sf S}^{\Psi}, d_{\sf
  U}^{\Psi})$, and $A^{\Phi_{\sf
      C}^*} = - A^{\Phi_{\sf C}}$,  $\gamma A^{\Psi_{\sf C}}$ are
  similar for some $\gamma \in \R^+$. In either case, therefore, (v)
  holds.

Finally, to prove that (v)$\Rightarrow$(i), note that this implication
clearly is correct if $d_{\sf
  H}^{\Phi} = d_{\sf S}^{\Phi} + d_{\sf U}^{\Phi}=0$. If $d_{\sf
  H}^{\Phi}\ge 1$ then no generality is
lost by assuming $(d_{\sf S}^{\Phi}, d_{\sf U}^{\Phi})
  = (d_{\sf S}^{\Psi}, d_{\sf U}^{\Psi})$ (otherwise replace $\Phi$ by $\Phi^*$) and $\alpha > 0$ (since
$A^{\Phi_{\sf C}}$, $-A^{\Phi_{\sf C}}$ are similar). Then $\lambda_j^{\Psi_{\sf
    H}}/ \lambda_j^{\Phi_{\sf H}}>0$ for every $j\in \{1, \ldots ,
d_{\sf H}^{\Phi} \}$, and Lemma \ref{lem61} shows that $\Phi
\stackrel{\beta}{\simeq} \Psi$ for all sufficiently small $\beta > 0$.
\end{proof}

To motivate one concise reformulation of Theorem \ref{thm6some}, note
that if $\Phi \stackrel{0^+}{\thicksim} \Psi$ then 
automatically $\{h(X_{\sf S}^{\Phi}), h(X_{\sf U}^{\Phi})\} =\{ 
X_{\sf S}^{\Psi}, X_{\sf U}^{\Psi} \} $, whereas simple examples show
that $h(X_{\sf H}^{\Phi}) \ne X_{\sf H}^{\Psi}$ and $h(X_{\sf
  C}^{\Phi}) \ne X_{\sf C}^{\Psi}$ in general. In other words,
some-H\"{o}lder equivalence between linear flows does not typically preserve
hyperbolic and central {\em spaces}. It does, however, preserve the
{\em flows\/} induced on these spaces, in the following sense.

\begin{cor}\label{cor62A}
Let $\Phi$, $\Psi$ be linear flows on $X$. Then $\Phi
\stackrel{0^+}{\thicksim} \Psi$ if and only if $\Phi_{\sf H}
\stackrel{0^+}{\simeq} \Psi_{\sf H}$ and $\Phi_{\sf C}
\stackrel{{\sf lin}}{\simeq} \Psi_{\sf C}$.
\end{cor}

\begin{proof}
To prove the ``if'' part, assume that $\Phi_{\sf H}
\stackrel{\beta}{\simeq} \Psi_{\sf H}$ for some $\beta > 0$, and
furthermore assume that $QA^{\Phi_{\sf C}} = \alpha A^{\Psi_{\sf C}}Q$ for some linear
isomorphism $Q: X_{\sf C}^{\Phi} \to X_{\sf C}^{\Psi}$ and $\alpha \in
\R \setminus \{0\}$. Now $\{d_{\sf S}^{\Phi},
d_{\sf U}^{\Phi} \} = \{  d_{\sf S}^{\Psi}, d_{\sf U}^{\Psi}\}$ by
Theorem \ref{thm6some}, and again it can be assumed that
$(d_{\sf S}^{\Phi} , d_{\sf U}^{\Phi})=(d_{\sf S}^{\Psi} , d_{\sf
  U}^{\Psi})=:(k,m)$ and $\alpha > 0$. If $(k,m)=0$ then there is
nothing to prove since $\Phi = \Phi_{\sf C} \stackrel{{\sf
    lin}}{\simeq} \Psi_{\sf C} = \Psi$. Assuming $k\ge 1$, recall from
Section \ref{sec3} that there exists a homeomorphism $f:X_{\sf
  S}^{\Phi} \to \R^k$ with $f(0)=0$ so that $f$, $f^{-1}$ both
satisfy, for instance, a $\frac12$-H\"{o}lder condition near $0$, and
\begin{equation}\label{eq6abbo1}
f(\Phi_{{\sf S}, t} x) = \widetilde{\Phi}_t f(x) \qquad \forall (t,x)
\in \R \times X_{\sf S}^{\Phi} \, ,
\end{equation}
with $\widetilde{\Phi}$ generated by $\mbox{\rm diag}\,
[\lambda_1^{\Phi_{\sf S}}, \ldots , \lambda_k^{\Phi_{\sf S}}]$;
similarly for $g:X_{\sf S}^{\Psi} \to \R^k$, where (\ref{eq6abbo1})
holds with $f$, $\Phi$ replaced by $g$, $\Psi$ respectively. Since
$\lambda_j^{\Psi_{\sf S}}/\lambda_j^{\Phi_{\sf S}}>0$ for every $j\in
\{1, \ldots , k\}$, it is possible to pick $0<\beta_{\sf S} < 1$ so
that
$$
\beta_{\sf S} < \min\nolimits_{j=1}^k  \min \left\{ \frac{\alpha \lambda_j^{\Psi_{\sf
        S}}}{\lambda_j^{\Phi_{\sf S}}},  \frac{  \lambda_j^{\Phi_{\sf
        S}}}{\alpha \lambda_j^{\Psi_{\sf S}}}\right\} \, .
$$
As seen in the proof of Lemma \ref{lem61}, there exists an $h \in
\cH_{\beta_{\sf S}}(\R^k)$ with $h(\widetilde{\Phi}_t x) =
\widetilde{\Psi}_{\alpha t} h(x)$ for all $t\in \R$, $x\in \R^k$. Letting
$h_{\sf S} = g^{-1} \circ h \circ f : X_{\sf S}^{\Phi} \to X_{\sf
  S}^{\Psi}$ yields a homeomorphism with $h_{\sf S}(0) =0$ so that
$h_{\sf S}$, $h^{-1}_{\sf S}$ both satisfy a $\frac14 \beta_{\sf
  S}$-H\"{o}lder condition near $0$. Assuming $m\ge 1$, a completely analogous argument
yields a homeomorphism $h_{\sf U}: X_{\sf U}^{\Phi} \to X_{\sf
  U}^{\Psi}$ with $h_{\sf U}(0) =0$ so that
$h_{\sf U}$, $h^{-1}_{\sf U}$ both satisfy a $\frac14 \beta_{\sf
  U}$-H\"{o}lder condition near $0$ for some $0<\beta_{\sf U}<1$. With
this, define $h:X\to X$ as
\begin{equation}\label{eq62a}
h(x) = h_{\sf S} (P_{\sf S}^{\Phi} x) + Q P_{\sf C}^{\Phi} x + h_{\sf U} (P_{\sf
  U}^{\Phi} x) \qquad \forall x \in X \, .
\end{equation}
Then $h \in \cH_{\gamma} (X)$ with $\gamma = \frac14 \min \{\beta_{\sf S},
\beta_{\sf U}\} >0$ by Proposition \ref{prop3add0}, and $h(\Phi_t x) = \Psi_{\alpha t}
h(x)$ for all $t\in \R$, $x\in X$. If $k=0$ or $m=0$, then the same
conclusion holds, though with the $h_{\sf S}$- or the $h_{\sf U}$-term deleted
from (\ref{eq62a}) and $\beta_{\sf S}:= 1$ or $\beta_{\sf U}:=1$,
respectively. In all cases, therefore, $\Phi \stackrel{0^+}{\simeq} \Psi$.

To prove the ``only if'' part, note that $\{d_{\sf S}^{\Phi}, d_{\sf U}^{\Phi}\} = \{d_{\sf
  S}^{\Psi}, d_{\sf U}^{\Psi}\}$ and $\Phi_{\sf C}
\stackrel{{\sf lin}}{\simeq} \Psi_{\sf C}$ by Theorem \ref{thm6some}. Thus, if
$(d_{\sf S}^{\Phi}, d_{\sf U}^{\Phi}) = (0,0)$ then there is nothing
else to prove. If $(d_{\sf S}^{\Phi}, d_{\sf U}^{\Phi}) \ne  (0,0)$
then  $\rho (\Phi_{\sf H},
\Psi_{\sf H})>0$, as seen in Section \ref{sec3}, and Corollary
\ref{cor48A} shows that $\Phi_{\sf H} \stackrel{\beta^-}{\simeq} \Psi_{\sf H}$
for every $0<\beta \le \sqrt{\rho(\Phi_{\sf H}, \Psi_{\sf H})}$.
\end{proof}

Unlike for its some-H\"{older} counter-part, a
characterization of {\em all}-H\"{o}lder equivalence does involve an
additional property of linear flows beyond the dimensions of their
(un)stable spaces, namely Lyapunov similarity. The ultimate result is
the following, mildly extended form of Theorem \ref{thmb}. 

\begin{theorem}\label{thm6all}
Let $\Phi$, $\Psi$ be linear flows on $X$. Then the following statements are
equivalent:
\begin{enumerate}
\item $\Phi \stackrel{1^-}{\simeq}\Psi$;
\item $\Phi \stackrel{1^-}{\thicksim} \Psi$;
  \item there exists an $\alpha \in \R\setminus \{0\}$ so that $A^{\Phi}$,
  $\alpha A^{\Psi}$ are Lyapunov similar and $A^{\Phi_{\sf
        C}}$, $\alpha A^{\Psi_{\sf C}}$ are similar;
  \item $\{d_{\sf S}^{\Phi}, d_{\sf U}^{\Phi}\} = \{d_{\sf S}^{\Psi},
    d_{\sf U}^{\Psi}\} $, and there exists an $\alpha \in \R^+$ so that
    $A^{\Phi_{\sf C}}$, $\alpha A^{\Psi_{\sf C}}$ are similar and 
    $$
\lambda_j^{\Phi_{\sf H}} = \alpha \lambda_j^{\Psi_{\sf H}} \quad
\mbox{or} \quad \lambda_j^{\Phi_{\sf H}} = -\alpha
\lambda_{d-j+1}^{\Psi_{\sf H}} \qquad \forall j \in \{1, \ldots , d_{\sf H}^{\Phi}\} \, .
$$
\end{enumerate}
\end{theorem}

The proof of Theorem
\ref{thm6all} is facilitated by two preparatory observations, Lemmas
\ref{lemH2} and \ref{lem65b} below. To motivate
the first of these lemmas, recall from Section \ref{sec3} that $\Phi_{\bullet} \stackrel{1^-}{\simeq}
\widetilde{\Phi}_{\bullet}$ for $\bullet \in \{{\sf S}, {\sf U}\}$,
where $\widetilde{\Phi}_{\bullet}$ is generated by
$\Lambda^{\Phi_{\bullet}}$. Thus $\Phi  \stackrel{1^-}{\simeq}
\widetilde{\Phi} $, with $\widetilde{\Phi}$ generated by $\mbox{\rm
  diag}\, [\Lambda^{\Phi_{\sf S}}, A^{\Phi_{\sf C}},
\Lambda^{\Phi_{\sf U}}]$. For convenience, let $( d_{\sf
  S}^{\Phi}, d_{\sf C}^{\Phi}, d_{\sf
  U}^{\Phi})=(k,\ell , m)$, so $k,\ell, m\in \N_0$ and $k+\ell+m =d $. Furthermore, if
$k\ge 1$ let
$$
S = \Lambda^{\Phi_{\sf S}} =  \mbox{\rm diag}\, [s_1, \ldots , s_k ] \in \R^{k\times k} \, ,
\quad \mbox{\rm with } s_1 \le \ldots \le s_k <0 \, ;
$$
if $\ell \ge 1$ let $C= A^{\Phi_{\sf C}}\in \R^{\ell \times \ell}$ with $\sigma (C)
\subset i\R$; and if $m\ge 1$ let
$$
U = \Lambda^{\Phi_{\sf U}} = \mbox{\rm diag}\, [u_1, \ldots , u_m] \in \R^{m\times m} \, ,
\quad \mbox{\rm with } 0 < u_1 \le \ldots \le u_m \, .
$$
With these ingredients, and for any $\alpha_{\sf S}, \alpha_{\sf
  U}\in \R^+$ and $\alpha_{\sf C} \in \R\setminus \{0\}$, consider the two $d\times d$-matrices
\begin{equation}\label{eqH1}
A = \mbox{\rm diag} \, [S,C,U] \, , \qquad B = \mbox{\rm diag} \,
[\alpha_{\sf S}S,\alpha_{\sf C}C, \alpha_{\sf U}U] \, ,
\end{equation}
where the $S$-, $C$-, and $U$-part in either matrix is
understood to be present only if $k\ge 1$, $\ell \ge 1$, and $m\ge 1$
respectively. For convenience, let $E_{\sf S} = \mbox{\rm span}
\{e_1, \ldots, e_k\}$, $E_{\sf C} = \mbox{\rm
  span}\{e_{k+1}, \ldots , e_{k+\ell}\}$, and $E_{\sf U} = \mbox{\rm
  span}\{e_{k+\ell+1}, \ldots , e_d\}$.  As presented below, the proof
of (ii)$\Rightarrow$(iii) in Theorem \ref{thm6all} 
crucially depends on the following auxiliary result.

\begin{lem}\label{lemH2}
Given $k, \ell, m\in \N_0$ with $k+\ell + m = d$, as well as $\alpha_{\sf S}, \alpha_{\sf
  U}\in \R^+$ and $\alpha_{\sf C}\in \R\setminus \{0\}$, let $\Phi$, $\Psi$ be the flows on $\R^d$ generated by $A$, $B$
in {\rm (\ref{eqH1})} respectively. Assume that $\Phi
\stackrel {h} {\thicksim} \Psi$ for some $h\in \cH_{1^-}(\R^d)$ with
$h(E_{\sf S})= E_{\sf S}$.
\begin{enumerate}
\item If $\min \{k,m\}\ge 1$ then $\alpha_{\sf S} = \alpha_{\sf U}$.
\item If $\min \{k,\ell\}\ge 1$ and $\sigma (C) \ne \{0\}$, then
  $\alpha_{\sf S} = |\alpha_{\sf C}|$.
\item If $\min \{\ell, m\}\ge 1$ and $\sigma (C) \ne \{0\}$, then
  $|\alpha_{\sf C}| = \alpha_{\sf U}$.
\end{enumerate}
\end{lem}

\begin{proof}
By Proposition \ref{lemH1}, it can be assumed that
$\Phi\stackrel {{\scriptstyle 1^-}} {\thickapprox} \Psi$, and since
$h(E_{\sf S})= E_{\sf S}$, clearly $\tau_x$ is increasing for every
$x\in \R^d\setminus \{0\}$; moreover,  $X_{\bullet}^{\Phi} =  E_{\bullet} = X_{\bullet}^{\Psi}$ for $\bullet
\in \{{\sf S}, {\sf C}, {\sf U}\}$. 
Denoting by $P_{\bullet}$ the orthogonal projection
of $\R^d$ onto $E_{\bullet}$, note that $P_{\bullet}$ commutes with
$\Phi_t$, $\Psi_t$ for every $t\in \R$. Moreover, if $k\ge 1$ then
\begin{equation}\label{eqH3}
e^{s_1 t}|x| \le |\Phi_t x| \le e^{s_k t} |x| \qquad \forall t \ge 0  ,  x \in E_{\sf S}\, ;
\end{equation}
if $\ell \ge 1$, and with an appropriate $\mu \in \R^+$,
\begin{equation}\label{eqH4}
|\Phi_t x| \le \mu \sqrt{1 + t^{2\ell -2}} |x| \qquad \forall t \in \R
,  x \in E_{\sf C} \, ;
\end{equation}
and if $m\ge 1$ then
\begin{equation}\label{eqH5}
e^{u_1 t}|x| \le |\Phi_t x| \le e^{u_m t} |x| \qquad \forall t\ge 0 ,  x \in E_{\sf U}\, .
\end{equation}
As a consequence, if $m\ge 1$, and with an appropriate $\nu \in
\R^+$,
\begin{equation}\label{eqH6}
|\Phi_t x| \le \nu e^{u_m t} |x|\qquad \forall t \ge 0 , x \in \R^d \, .
\end{equation}
Similar universal bounds are valid with $\Psi$ instead of $\Phi$,
provided that $s_1, s_k$ and $u_1, u_m$ are replaced by $\alpha_{\sf
  S} s_1, \alpha_{\sf S} s_k$ and $\alpha_{\sf U} u_1, \alpha_{\sf U}
u_m$ respectively. By assumption, $\Phi
\stackrel {h} {\thicksim} \Psi$ with $h\in \cH_{1^-}$ and $h(E_{\sf
  S}) = E_{\sf S}$. Recall that this implies $h(E_{\sf U}) = E_{\sf U}$,
whereas it is possible that $h(E_{\sf C}) \ne E_{\sf C}$. To prepare
for the elementary but somewhat intricate arguments below, fix
$0<\beta < 1$. By means of an appropriate rescaling, it can be assumed
that $h,h^{-1}$ both satisfy a $\beta$-H\"{o}lder condition on $V:=
B_2(0) \cup h\bigl( B_2(0) \bigr)$, i.e., with some $\kappa\in \R^+$,
$$
|h(x) - h(y)| +  |h^{-1}(x) - h^{-1} (y)| \le \kappa |x -
y|^{\beta} \qquad \forall x,y \in V \, .
$$

\medskip

\noindent
{\bf Proof of (i):} Assume that $\min \{k,m\}\ge 1$. Notice that $\Psi$
can, and henceforth will be assumed to be generated by $\mbox{\rm
  diag}\, [S,\alpha_{\sf C} C, \alpha_{\sf U} U]$, i.e., assume
w.l.o.g.\ that $\alpha_{\sf S}=1$. Establishing (i) therefore amounts
to proving that $\alpha_{\sf U} = 1$. To this end, for every
$r\in \R^+$ let $x_r = e_k + e^{-u_m r} e_d$, and
hence
$$
\Phi_t x_r = e^{s_k t} e_k + e^{u_m (t - r)} e_d \qquad
\forall t \in \R \, ,
$$
from which it is clear that $|\Phi_t x_r|< \sqrt{2}$ for every
$t\in [0, r ]$. Also, $x_r \to e_k$ as $r\to \infty$, and hence
$h(x_r)\to h(e_k)\in E_{\sf S}\setminus \{0\}$, whereas $\Phi_{r} x_r =
e^{s_k r } e_k + e_d \to e_d$ and 
$h(\Phi_{r} x_r) \to h(e_d)\in E_{\sf U} \setminus
\{0\}$. This shows that $P_{\sf U} h(x_r)\ne 0$ for all $r>0$
sufficiently large since otherwise $h (\Phi_{r}
x_r) \in \Psi_{\R} h(x_r) \subset E_{\sf S}
\oplus E_{\sf C}$, which clearly is not the case for large $r$. With an appropriate
$r_1\ge 0$, therefore, $P_{\sf U} h(x_r) \ne 0$ for every
$r > r_1$, and it makes sense to consider
$$
y_r:= h^{-1} \left(
h(x_r) + \frac{|P_{\sf C} h(x_r)|^{1/\beta^3}}{ |P_{\sf U}
  h(x_r)|} P_{\sf U} h(x_r)
  \right) \qquad \forall r > r_1 \, .
$$
Recalling that $h(e_k) \in E_{\sf S}$ and hence $P_{\sf C} h(e_k)=0$,
deduce from
\begin{align*}
|h(y_r) - h(x_r)|  = |P_{\sf C} h(x_r)|^{1/\beta^3} =
                              \big|
P_{\sf C}  \bigl( h(x_r) - h(e_k) \bigr)
                              \big|^{1/\beta^3} \le |h(x_r) - h
                              (e_k)|^{1/\beta^3} & \le \kappa^{1/\beta^3}|x_r
                              - e_k|^{1/\beta^2} \\
  & = \kappa^{1/\beta^3} e^{-u_m r/\beta^2} 
\end{align*}
that $h(y_r) \in V$ for every $r>r_2$, with an
appropriate $r_2 \ge r_1$. Consequently,
\begin{equation}\label{eqH7}
|y_r - x_r| \le \kappa |h(y_r) - h(x_r)|^{\beta}
\le \kappa^{1+1/\beta^2} e^{-u_m r/\beta} \qquad \forall r > r_2 
\, .
\end{equation}
Together with (\ref{eqH6}), this yields for every $r> r_2$,
\begin{equation}\label{eqH8}
|\Phi_t y_r - \Phi_t x_r| \le \nu  e^{u_m t} |y_r -
x_r| \le \nu \kappa^{1+1/\beta^2} e^{-u_m r(1/\beta - 1)} \qquad
\forall t \in [0,  r] \, .
\end{equation}
Since $\beta< 1$, picking $r_3 \ge r_2$ sufficiently large
guarantees that $\Phi_t y_r \in V$ for all $t\in [0, r]$ and $r >
r_3$. For convenience, henceforth denote $\tau_{y_r}
(r )>0$ simply by $T_r$.

First, rough (lower and upper) bounds on $T_r$ are going to be
established; these bounds will show in particular that, as the reader
no doubt suspects already, $T_r \to \infty$ as $r\to \infty$. To obtain a {\em lower\/}
bound, recall that $h(\Phi_{r}x_r) \to
h(e_d) \in E_{\sf U} \setminus \{0\}$ as $r\to \infty$, and hence by
(\ref{eqH8}) also $P_{\sf U}h(\Phi_{r }y_r ) \to P_{\sf U}
h(e_d) = h(e_d) \ne 0$. Thus, using the notation $f(r) \prec g(r)$
exactly as in the proof of Theorem \ref{lemr9},
\begin{align*}
{\textstyle \frac12} |h(e_d)|  \prec |P_{\sf U} h (\Phi_{r} y_r)| =
                   |P_{\sf U} \Psi_{T_r} h(y_r)| & = \left|
                   \Psi_{T_r} \left(
P_{\sf U} h(x_r) + \frac{|P_{\sf C}
                   h(x_r)|^{1/\beta^3}}{|P_{\sf U} h(x_r)|}
                   P_{\sf U} h(x_r)
                   \right)
                   \right| \\
  & \le  e^{\alpha_{\sf U} u_m T_r}\left(|P_{\sf U} h(x_r)| +
    |P_{\sf C} h(x_r)|^{1/\beta^3}\right) \, ,
\end{align*}
where the last inequality is due to the $\Psi$-version of
(\ref{eqH5}). This yields
\begin{equation}\label{eqH9}
e^{T_r} \succ  \left( |P_{\sf U} h(x_r)| +
    |P_{\sf C} h(x_r)|^{1/\beta^3}\right) ^{-1/(\alpha_{\sf U} u_m)} \, .
\end{equation}
Note that $|P_{\sf U} h(x_r)|, |P_{\sf C} h(x_r)|\to 0$ as $r\to
\infty$, so (\ref{eqH9}) shows that indeed $T_r \to \infty$, as suspected.

Similarly, to obtain an {\em
  upper\/} bound on $T_r$, deduce from
\begin{align*}
2 |h(e_d)|  \succ |P_{\sf U} h (\Phi_{r } y_r)| 
 & = \left|  \Psi_{T_r} \left(
P_{\sf U} h(x_r) + \frac{|P_{\sf C}
                   h(x_r)|^{1/\beta^3}}{|P_{\sf U} h(x_r)|}
                   P_{\sf U} h(x_r)
                   \right)
                   \right| \\
  & \ge e^{\alpha_{\sf U} u_1 T_r} \left(|P_{\sf U} h(x_r)| +
    |P_{\sf C} h(x_r)|^{1/\beta^3} \right) \, ,
\end{align*}
with the last inequality again due to the $\Psi$-version of
(\ref{eqH5}), that
$$
e^{T_r} \prec  \left( |P_{\sf U} h(x_r)| +
    |P_{\sf C} h(x_r)|^{1/\beta^3} \right)^{-1/(\alpha_{\sf U} u_1)} \, .
    $$
It is possible, therefore, to pick $r_4 \ge r_3$ so large that
\begin{equation}\label{eqH10}
e^{T_r} \le  \left(|P_{\sf U} h(x_r)| +
    |P_{\sf C} h(x_r)|^{1/\beta^3}\right)^{-1/(\alpha_{\sf U} \beta u_1)}
    \qquad \forall r > r_4 \, ;
  \end{equation}
notice the additional factor $\beta$ in the exponent on the right
allowing for the unspecified upper bound implied by $\prec$ to be
chosen as $1$, or indeed any positive constant.

Building on these preparations, the overall strategy now is to estimate
$|P_{\sf S} \Phi_{r} y_r|$ in two different ways: First directly,
utilizing (\ref{eqH3}) and (\ref{eqH7}), which leads to a lower bound,
and then by considering $\big|P_{\sf S} h^{-1} \bigl(  \Psi_{T_r}
h(y_r)\bigr)\big|$ instead, utilizing the bounds on $T_r$ just established,
which leads to an upper bound; see also Figure \ref{fig3}.
Concretely, deduce from (\ref{eqH3}) that on the one hand
\begin{align*}
|P_{\sf S} \Phi_{r } y_r| & = |P_{\sf S} \Phi_{r} x_r  + P_{\sf S}  \Phi_{r} (y_r - x_r)|
                                         = |e^{s_k r} e_k + P_{\sf S} \Phi_{r} (y_r -
                                         x_r)| \\
& \ge e^{s_k r} - e^{s_k r} |y_r - x_r| \ge e^{s_k r} \left( 1  -  \kappa^{1+1/\beta^2}
                                                         e ^{-u_m r
                                                    /\beta}\right)  \succ
                                                    e^{s_k r}
                                                         \, ;
\end{align*}
here the first and second $\ge$ are due to the reverse triangle
inequality and (\ref{eqH7}), respectively. Thus clearly
\begin{equation}\label{eqH11}
|P_{\sf S} \Phi_{r} y_r| \succ   e^{-|s_k| r} \, .
\end{equation}
On the other hand, let
$z_r = h^{-1} \bigl(\Psi_{T_r} P_{\sf U} h(y_r) \bigr)$
for convenience, so $z_r \in E_{\sf U}$. Note that $z_r \to e_d$ as
$r\to \infty$, and hence $z_r \in E_{\sf U} \cap V$ for
every $r>r_5$ with an appropriate $r_5\ge r_4$. Since $P_{\sf S} z_r =0$,
\begin{align*}
  |P_{\sf S} \Phi_{r } y_r| & = \big|  P_{\sf S} \bigl( h^{-1} \circ h (\Phi_{r } y_r) -
                                           h^{-1} \circ h (z_r)
                               \bigr)   \big| \le
                               \big| h^{-1} \bigl( \Psi_{T_r} h(y_r) \bigr) - h^{-1} \circ h (z_r)
                                           \big| \\
                                         & \le \kappa | \Psi_{T_r} h(y_r) - h(z_r)|^{\beta} =
                                           \kappa ( |\Psi_{T_r}
                                           P_{\sf S} h(x_r)|^2 + |\Psi_{T_r}
                                           P_{\sf C}
                                           h(x_r)|^2)^{\beta/2}
  \\
                                         & \le \kappa \left(
e^{2s_k  T_r} |P_{\sf S} h(x_r)|^2 + \mu^2 (1+
                                           T_r^{2\max \{\ell,1\} - 2})
                                           |P_{\sf C} h(x_r)|^2
                                           \right)^{\beta/2} \, ,
\end{align*}
with the last inequality due to the $\Psi$-versions of (\ref{eqH3})
and (\ref{eqH4}). (Recall that $\alpha_{\sf S}=1$.) Moreover, note
that for every $r> 0$, 
$$
|P_{\sf U} h(x_r)| = \big|
P_{\sf U} \bigl( h(x_r) - h (e_k) \bigr)
\big| \le |h(x_r) - h(e_k)| \le \kappa |x_r - e_k|^{\beta} =
\kappa e^{-\beta u_m r} \, ,
$$
and similarly $|P_{\sf C} h(x_r)| \le \kappa e^{-\beta u_m r}$. This yields
\begin{equation}\label{eqH111}
|P_{\sf U} h(x_r)| + |P_{\sf C} h(x_r)|^{1/\beta^3} \le
\kappa e^{-\beta u_m r} + \kappa^{1/\beta^3} e^{-u_mr/\beta^2} =
\kappa e^{-\beta u_m r} \left( 1 + \kappa^{1/\beta^3 - 1} e^{- u_m r
  (1/\beta^2 - \beta)} \right) \, ,
\end{equation}
and hence (\ref{eqH9}) implies that
$$
e^{2s_k T_r}  \prec  \left(
  \kappa e^{-\beta u_m r } \left(
1 + \kappa^{1/\beta^3 - 1} e^{-u_m r  (1/\beta^2 - \beta)}
    \right)
               \right)^{-2s_k /(\alpha_{\sf U} u_m) }
\prec e^{2s_k \beta r/\alpha_{\sf U}} \, .
$$
Since $P_{\sf S} h(x_r) \to P_{\sf S} h(e_k) = h(e_k) \ne 0$ as
$r\to \infty$, it is clear that
\begin{equation}\label{eqH9A}
e^{2 s_k T_r} |P_{\sf S} h(x_r)|^2  \prec  e^{-2|s_k| \beta r
  /\alpha_{\sf U}} \, .
\end{equation}
By contrast, deduce from (\ref{eqH10}) that
$$
T_r \le - \frac1{\alpha_{\sf U} \beta u_1 } \log \left( |P_{\sf U}
h(x_r)| + |P_{\sf C} h(x_r)|^{1/\beta^3}\right) \qquad \forall r > r_4 \, ,
$$
and hence
\begin{align*}
  (1 & +T_r^{2\max \{\ell, 1\} - 2})  |P_{\sf C} h(x_r)|^2  \le \\
  & \qquad \le \left(
1 + \frac1{(\alpha_{\sf U} \beta u_1 )^{2\max\{\ell,1\} - 2}} \Big| \log \left( |P_{\sf U}
h(x_r)| + |P_{\sf C} h(x_r)|^{1/\beta^3} \right) \Big|^{2\max\{\ell,1\} - 2}
  \right) |P_{\sf C} h(x_r)|^2 \\
 & \qquad  \prec \Big| \log \left( |P_{\sf U}
h(x_r)| + |P_{\sf C} h(x_r)|^{1/\beta^3} \right) \Big|^{2\ell}  \left(|P_{\sf U}
  h(x_r)| + |P_{\sf C} h(x_r)|^{1/\beta^3}\right)^{2\beta^3}  \\
 &  \qquad \prec  \left(|P_{\sf U}
   h(x_r)| + |P_{\sf C} h(x_r)|^{1/\beta^3}\right)^{2\beta^4}  \prec e^{-2\beta^5 u_m r} \, ;  
\end{align*}
here the second $\prec$ is due to the fact that 
$\sup_{0<u<1} |\log u|^a u^b <\infty$ for every $a,b\in \R^+$, whereas
the last $\prec$ is due to (\ref{eqH111}). Using this and (\ref{eqH9A}), therefore, 
\begin{equation}\label{eqH10A}
|P_{\sf S} \Phi_{r} y_r| \prec \left( e^{-2|s_k| \beta r/\alpha_{\sf
      U}} + e^{-2\beta^5 u_m r}\right)^{\beta /2} \prec e^{-r \min \{
  |s_k| \beta^2 /\alpha_{\sf U}, \beta^6
  u_m\}}  \, .
\end{equation}

\begin{figure}[ht]
\psfrag{thh}[]{$h$}
\psfrag{thhm1}[]{$h^{-1}$} 
\psfrag{tes}[]{$E_{\sf S}$}
\psfrag{tesu}[l]{$E_{\sf H}$}
\psfrag{thesu}[l]{$h(E_{\sf H})$}
\psfrag{tec}[]{$E_{\sf C}$}
\psfrag{teu}[]{$E_{\sf U}$}
\psfrag{t0}[]{$0$}
\psfrag{ted}[]{$e_d$}
\psfrag{thed}[l]{$h(e_d)$}
\psfrag{ptrxr}[r]{$\Phi_{r}x_r$}
\psfrag{hptrxr}[r]{$h(\Phi_{r}x_r)$}
\psfrag{ptryr}[l]{$\Phi_{r}y_r$}
\psfrag{hptryr}[l]{$\Psi_{T_r} h(y_r) = h(\Phi_{-r}y_r)$}
\psfrag{hptryra}[l]{$h(\Phi_{r}y_r)$}
\psfrag{hptryrb}[l]{$= \! \Psi_{T_r} h(y_r)$}
\psfrag{tptryr}[]{$P_{\sf S}\Phi_{r}y_r$}
\psfrag{tyr}[]{$y_r$}
\psfrag{txr}[]{$x_r$}
\psfrag{thyr}[]{$h(y_r)$}
\psfrag{thxr}[]{$h(x_r)$}
\psfrag{tzr}[r]{$z_r$}%=h^{-1} \bigl( P_{\sf U} \Psi_{T_r}
% h(y_r)\bigr)$}
\psfrag{thzr}[r]{$h(z_r)\! =\! P_{\sf U} \Psi_{T_r} h(y_r)$}
\psfrag{tek}[]{$e_k$}
\psfrag{thek}[]{$h(e_k)$}
\psfrag{tpyr}[]{$P_{\sf S} y_r$}
%
% scale 0.85
%
\begin{center}
\includegraphics{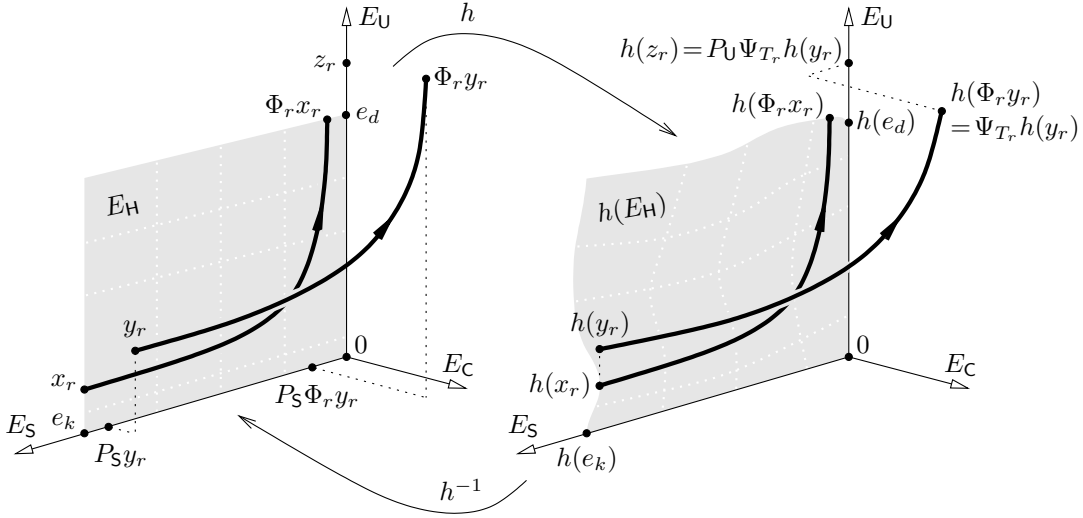}
\end{center}
\caption{Proving Lemma \ref{lemH2}(i) by estimating $|P_{\sf S}
  \Phi_{r} y_r| = \big|P_{\sf S} h^{-1} \bigl(\Psi_{T_r} h(y_r)
  \bigr) \big|$ in two different ways which lead to (\ref{eqH11}) and
(\ref{eqH10A}), respectively.}\label{fig3}
\end{figure}

Thus both desired estimates for $|P_{\sf S} \Phi_{r} y_r|$ alluded to
earlier have been obtained, in the form of the lower bound
(\ref{eqH11}) and the upper bound (\ref{eqH10A}). Combining these yields
$$
e^{-|s_k|r} \prec e^{-r \min \{ |s_k| \beta^2 /\alpha_{\sf U}, \beta^6
  u_m\}} \, .
$$
It follows that $|s_k| \ge \min \{|s_k|\beta^2/\alpha_{\sf U},
\beta^6 u_m\}$. Recall that $0<\beta < 1$ has been arbitrary, so
letting $\beta \uparrow 1$ yields $|s_k| \ge
\min \{|s_k|/\alpha_{\sf U}, u_m\}$. In other words, $1 \ge \min \{1/\alpha_{\sf
  U}, u_m/|s_k|\}$, or equivalently
\begin{equation}\label{eqH11A}
  1 \le \max \left \{
\alpha_{\sf U}, \frac{|s_k|}{u_m}
    \right\} \, .
\end{equation}
Identical reasoning with the roles of $\Phi$, $\Psi$ interchanged
yields (\ref{eqH11A}) with $\alpha_{\sf U}$, $u_m$ replaced by
$1/\alpha_{\sf U}$, $\alpha_{\sf U} u_m$ respectively, that is,
$\alpha_{\sf U} \le \max \{1, |s_k|/u_m\}$. In total, therefore,
\begin{equation}\label{eqH12}
  1 \le \max \left\{
\alpha_{\sf U} , \frac{|s_k|}{u_m}
\right\} \quad \mbox{\rm and} \quad
\alpha_{\sf U} \le \max \left\{
1, \frac{|s_k|}{u_m}
  \right\} \, .
\end{equation}
Note that $\Phi^* \stackrel{1^-}{\thicksim} \Psi^*$ also, with 
$\Phi^*$, $\Psi^*$, generated by $\mbox{\rm diag}\, [-U, -C, -S]$, 
$\mbox{\rm diag}\, [-\alpha_{\sf U}U, - \alpha_{\sf C}C, - \alpha_{\sf
S} S]$ respectively. Identical reasoning shows that (\ref{eqH12})
remains valid in this situation as well, provided that $\alpha_{\sf U}, |s_k|,
u_m$ are replaced by $1/\alpha_{\sf U}, u_1, |s_1|$
respectively. This yields
\begin{equation}\label{eqH13}
  1 \ge \min \left\{
\alpha_{\sf U}, \frac{|s_1|}{u_1}
    \right\} \quad \mbox{\rm and} \quad  \alpha_{\sf U} \ge \min
    \left\{
1, \frac{|s_1|}{u_1}
      \right\} \, .
 \end{equation}
    
 Now, in order to conclude the argument by combining (\ref{eqH12}) and (\ref{eqH13}), it is
 helpful to distinguish three disjoint cases corresponding to the
 possible value of $|s_k|/u_m \in \R^+$: First, if
 $|s_k|/u_m < 1$ then the left inequality in (\ref{eqH12}) yields
 $\alpha_{\sf U} \ge 1$, whereas the right inequality reads
 $\alpha_{\sf U}\le 1$. Thus, $\alpha_{\sf U} = 1$. Second, if
 $|s_k|/u_m = 1$ then the left inequality in (\ref{eqH12}) automatically
 holds, whereas the right inequality reads $\alpha_{\sf U} \le
 1$. Moreover, $|s_1|/u_1 \ge |s_k|/u_1 \ge |s_k|/u_m = 1$, and so the right
 inequality in (\ref{eqH13}) reads $\alpha_{\sf U} \ge 1$. Again,
 therefore, $\alpha_{\sf U} =1$. Finally, if $|s_k|/u_m >1$ then
 $|s_1|/u_1>1$, so the left and right inequalities in (\ref{eqH13}) yield
 $\alpha_{\sf U} \le 1$ and $\alpha_{\sf U} \ge 1$ respectively,
 hence $\alpha_{\sf U} = 1$ once again. Thus $\alpha_{\sf U} = 1$ in
 all three cases. This completes the proof of (i).

\medskip
   
\noindent
{\bf Proof of (ii):} Assume that $\min \{k,\ell\}\ge 1$ and $\sigma(C)\ne
 \{0\}$. As in the above proof of (i), it can be assumed w.l.o.g.\ that
 $\alpha_{\sf S} = 1$, so establishing (ii) amounts to proving that
 $|\alpha_{\sf C}| = 1$. To this end, note that $\mbox{\rm Per}\, \Phi
 \setminus \mbox{\rm Fix}\, \Phi \ne \varnothing$. Pick $r_1>0$ so
 small that $\Phi_{\R}p \subset V$ and $h ( \Phi _{\R} p)  =
 \Psi _{\R} h(p) \subset V$ for every $p\in B_{r_1}(0) \cap
 \mbox{\rm Per}\, \Phi $. For the following argument, pick any $p\in
 B_{r_1}(0) \cap (\mbox{\rm Per}\, \Phi \setminus \mbox{\rm Fix}\, \Phi )$.
Clearly $h(p)\in \mbox{\rm Per}\, \Psi \setminus \mbox{\rm Fix}\, \Psi $. For
 convenience, write $T_p^{\Phi}$, $T_{h(p)}^{\Psi}$ as $T$,
 $S$ respectively, so $T, S\in \R^+$. For all that follows, it will
 be useful to consider the set
$$
K_p:= \bigl\{ x\in \R^d : P_{\sf U} x = 0 , P_{\sf C} x \in \Phi_{\R} p \bigr\}
\subset E_{\sf S} \oplus E_{\sf C} \, .
$$
It is readily seen that $x\in K_p$ precisely if $\Phi_t x$ approaches
the compact set or ``loop'' $\Phi_{\R} p$ as $t\to \infty$, or
equivalently if, given any sequence $(t_n)$ in $\R$ with $t_n \to \infty$,
there exists a sequence $(s_n)$ with $0\le s_n < T$ so that $\lim_{n\to
  \infty}|\Phi_{t_n} x - \Phi_{s_n} p|=0$. The set $K_p$ clearly is
$\Phi$-invariant; it may be thought of as a $(k+1)$-dimensional ``cylinder'' over the closed orbit $\Phi_{\R}
p$. Note that $p+E_{\sf S} \subset K_p$, and $\Phi_T (p+ E_{\sf S}) =
p + E_{\sf S}$. Thus, with $\iota_p$ denoting the isometry
$$
\iota_p : \left\{
  \begin{array}{ccl}
    \R^k & \to & p + E_{\sf S} \, , \\
    x & \mapsto & p + \sum_{j=1}^k x_j e_j \, , 
    \end{array}
  \right.
$$
the map $\Phi_T$ induces the linear (Poincar\'{e}) map $F^{\Phi}: \R^k \to \R^k$
given by
$$
F^{\Phi} = \iota_p^{-1} \circ \Phi_T|_{p+E_{\sf S}} \circ \iota_p =
\mbox{\rm diag}\left[ e^{s_1 T}, \ldots , e^{s_k T} \right] \, ;
$$
see Figure \ref{fig4}.
\begin{figure}[ht]
  \psfrag{tx}[r]{$x$}
  \psfrag{tgx}[l]{$G(x)$}
  \psfrag{thh}[]{$h$}
   \psfrag{thhi}[]{$h^{-1}$}
\psfrag{tip}[]{$\iota_p$}
\psfrag{tihp}[]{$\iota_{h(p)}$}
\psfrag{tfpx}[r]{$F^{\Phi} x$}
\psfrag{tfpgx}[l]{$F^{\Psi} G( x)\! = \! G(F^{\Phi} x)$} 
\psfrag{trk}[l]{$\R^k$}
\psfrag{t0}[r]{$0$}
\psfrag{t0r}[l]{$0$}
\psfrag{tkp}[]{$K_p$}
\psfrag{tkhp}[]{$K_{h(p)}\! =\! h(K_p)$}
\psfrag{tpx}[l]{$p+x$}
\psfrag{tgpx}[r]{$g(p+x)$}
\psfrag{thpx}[l]{$h(p+x)$}
\psfrag{tptpx}[l]{$\Phi_{T}(p+x)$}
\psfrag{tppes}[]{$p+E_{\sf S}$}
\psfrag{thppes}[]{$h(p+E_{\sf S})$}
\psfrag{tphpes}[]{$h(p)+E_{\sf S}$}
\psfrag{tp}[l]{$p$}
\psfrag{torbp}[]{$\Phi_{\R} p$}
\psfrag{torbhp}[]{$\Psi_{\R} h(p)$}
\psfrag{thp}[r]{$h(p)$}
\psfrag{tes}[l]{$E_{\sf S}$}
%
% scale 0.70
%
\begin{center}
\includegraphics{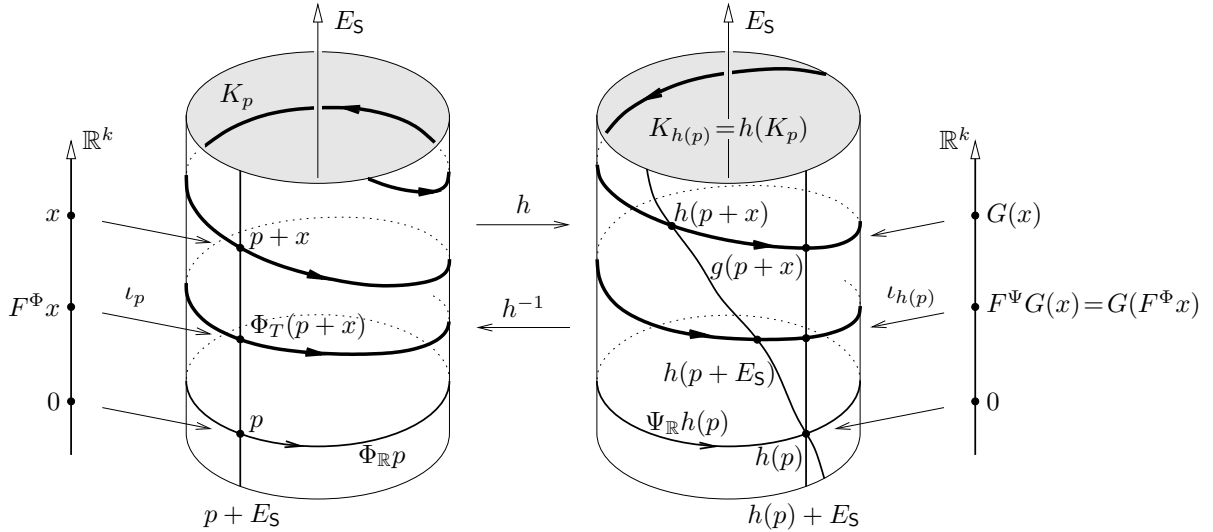}
\end{center}
\caption{Proving Lemma \ref{lemH2}(ii) by linking the 
  Poincar\'{e} maps $F^{\Phi}$, $F^{\Psi}$ induced on $\R^k$, via the
  homeomorphism $G = \iota_{h(p)}^{-1}
  \circ g \circ \iota_p$.}\label{fig4}
\end{figure}
A completely analogous construction, utilizing the $\Psi$-invariant
``cylinder'' $K_{h(p)}$ over $\Psi_{ \R} h(p)$, yields the
linear map $F^{\Psi}: \R^k \to \R^k$ given by
$$
F^{\Psi} = \iota_{h(p)}^{-1} \circ \Psi_S|_{h(p)+E_{\sf S}} \circ \iota_{h(p)} =
\mbox{\rm diag}\left [e^{s_1 S}, \ldots , e^{s_k S} \right] \, ,
$$
induced in this case by $\Psi_S$ as $\Psi_S( h(p) + E_{\sf
  S}) = h(p)+E_{\sf S}$. Now, to link $F^{\Phi}$, $F^{\Psi}$ by means of the
homeomorphism $h$, notice that $h(K_p) = K_{h(p)}$. While a
point $h(p+x)$ with $x\in E_{\sf S}$ will not in general be an element
of $h(p) + E_{\sf S}$ if $x\ne 0$, letting the point flow with $\Psi$
for a small amount of time will bring it into $h(p)+E_{\sf
  S}$. Formally, there exists $0<r_2\le r_1$ and a smooth function $\theta:
K_{h(p)}\cap B_{r_2} \bigl( h(p) \bigr)\to \R$ with $\theta \bigl(
h(p)\bigr)= 0$ so that $\Psi_{\theta (y)}y \in h(p)+ E_{\sf S}$ for all
$y \in K_{h(p)}\cap B_{r_2} \bigl( h(p) \bigr)$. Pick $0<r_3 \le r_2$
small enough to ensure $B_{r_3} (p)\subset V$ as well as $h \bigl(
B_{r_3} (p)\bigr)\subset B_{r_2} \bigl( h(p) \bigr)$, and define
$$
g :\left\{
  \begin{array}{ccl}
    (p + E_{\sf S})\cap B_{r_3} (p) & \to & h(p) + E_{\sf S} \, , \\
    p + x & \mapsto & \Psi_{\theta ( h(p+x) )} h(p+x) \, .
    \end{array}
  \right.
$$
Clearly $g(p) = h(p)$. Since $\theta$ is smooth and $h$ satisfies a
$\beta$-H\"{o}lder condition on $(p + E_{\sf S}) \cap B_{r_3}(p) = p +
E_{\sf S} \cap B_{r_3}(0)$, so does $g$. Furthermore, $g$ is invertible,
and $g^{-1}$ satisfies a $\beta$-H\"{o}lder condition on $h(p) + E_{\sf
  S}\cap B_{r_4}(0)$ for any sufficiently small $0<r_4 \le r_3$.
Since $T$, $S$ are the minimal periods of $p$, $h(p)$
respectively,
$$
g \bigl( \Phi_T (p+x) \bigr) = \Psi_S g (p+x) \qquad \forall x \in
E_{\sf S} \cap B_{r_4} (0) \, .
$$
With $G: B_{r_4}(0) \to \R^k$ given by $G = \iota_{h(p)}^{-1} \circ g
\circ \iota_p$, this means that
$$
G ( F^{\Phi} x) = F^{\Psi}  G (x) \qquad \forall x \in
B_{r_4} (0) \, .
$$
Note that $G(0)=0$, and $G$ is a local homeomorphism, as is $G^{-1} =
\iota_p^{-1} \circ g^{-1} \circ \iota_{h(p)}$; moreover, $G$, $G^{-1}$
both satisfy a $\beta$-H\"{o}lder condition near $0$. Since $F^{\Phi}$, $F^{\Psi}$ are
contractions on $\R^k$, it is clear that with an appropriate $0<r_5\le r_4$,
\begin{equation}\label{eqxx1}
(F^{\Phi})^n (x) = G^{-1} \circ (F^{\Psi})^n \circ G (x) \qquad
\forall n \in \N , x \in B_{r_5}(0) \, .
\end{equation}
With $\kappa$ denoting a $\beta$-H\"{o}lder constant for $G$, $G^{-1}$
on $B_{r_5}(0)$, take $x= \frac12 r_5 e_k$ and deduce from
(\ref{eqxx1}) that
$$
{\textstyle \frac12} e^{n s_k  T} r_5 \le \kappa \left| (F^{\Psi})^n
  \circ G \left( {\textstyle \frac12} r_5
e_k\right)\right|^ \beta \le \kappa e^{n s_k  S \beta} \left| G \left(
  {\textstyle \frac12} r_5 e_k\right) \right|^{\beta}
\qquad \forall n \in \N \, .
$$
Letting $n\to \infty$ yields $T - S \beta \ge 0$, i.e., $T/S\ge
\beta$. Since $0<\beta < 1$ has been arbitrary, it follows that $T\ge
S$, and interchanging the roles of $\Phi$, $\Psi$ yields $T=S$.

Recall that $p \in B_{r_1}(0) \cap (\mbox{\rm Per}\, \Phi \setminus
\mbox{\rm Fix}\, \Phi)$ has been arbitrary. In summary, therefore, it has
been shown that $T_p^{\Phi} = T_{h(p)}^{\Psi}$ for 
every $p \in B_{r_5} (0) \cap (\mbox{\rm Per}\, \Phi \setminus
\mbox{\rm Fix}\, \Phi)$. Clearly, $T_x^{\Phi} = 0 = T_{h(x)}^{\Psi}$
for every $x\in  B_{r_5} (0) \cap  \mbox{\rm Fix}\, \Phi$, whereas if
$x\in B_{r_5}(0)\setminus  \mbox{\rm Per}\, \Phi$ then $T_x^{\Phi} =
\infty = T_{h(x)}^{\Psi}$. In other words, $T_x^{\Phi} =
T_{h(x)}^{\Psi}$ for every $x\in B_{r_5}(0)$, and so Theorem \ref{thm2b} yields
$$
\sigma (C) = \sigma (\Phi) \cap i \R = \sigma (\Psi) \cap i\R = \sigma
(\alpha_{\sf C} C)\, .
$$
Since $\sigma(C) \ne \{0\}$, necessarily
$\alpha_{\sf C} \in \{-1,1\}$. Thus $|\alpha_{\sf C}| = 1$, and the proof of (ii) is complete.

\medskip 

\noindent
{\bf Proof of (iii):} Assume that $\min \{\ell, m\}\ge 1$ and
$\sigma (C)\ne \{0\}$. In this case, it can be assumed w.l.o.g.\ that
$\alpha_{\sf U} = 1$, i.e., $\Psi$ is generated by $\mbox{\rm diag}\,
[\alpha_{\sf S}S, \alpha_{\sf C} C, U]$, and it only needs to be shown
that $|\alpha_{\sf C}| = 1$. To this end, recall that
$\Phi^*\stackrel{1^-}{\thicksim}\Psi^*$, with $\Phi^*$, $\Psi^*$ generated by $\mbox{\rm diag}\,
[- U, -C, -S]$, $\mbox{\rm diag}\,
[- U, - \alpha_{\sf C}C, -\alpha_{\sf S} S]$ respectively. Applying
(ii), with $m,k,-U,-C,-S$, and $\alpha_{\sf S}$ instead of $k,m,S,C,U$,
and $\alpha_{\sf U}$ respectively, yields $|\alpha_{\sf C}| = 1$ and
hence completes the proof overall.
\end{proof}

In Lemma \ref{lemH2}, note that if $k\ne m$ then the additional condition
$h(E_{\sf S}) = E_{\sf S}$ automatically is satisfied. By contrast, if
$k=m\ge 1$ then that condition is essential, as can be seen, for instance,
from the flows $\Phi$, $\Psi$ on $\R^4$ generated by $\mbox{\rm
  diag}\, [-1, J_1(i), 2]$, $\mbox{\rm diag}\, [-4, J_1(2i), 2]$
respectively, for which $\Phi \stackrel{{\sf lin}}{\simeq} \Psi$, and
yet $(\alpha_{\sf S}, \alpha_{\sf C}, \alpha_{\sf U}) = (4,2,1)$, so
all three conclusions in Lemma \ref{lemH2} fail.

The second preparatory observation for the proof of Theorem
\ref{thm6all} is a strengthening of Theorem \ref{lemr9} as $\beta
\uparrow 1$ which may also be of independent interest. Informally put,
it asserts that for stable flows the Lyapunov exponents, and in fact
even the Lyapunov {\em spaces} as defined in
Section \ref{sec2a}, behave naturally under all-H\"{o}lder equivalence.

\begin{lem}\label{lem65b}
Let $\Phi$, $\Psi$ be stable flows on $X$. Assume that $\Phi
\stackrel{h}{\thicksim} \Psi$ for some $h\in \cH_{1^-}(X)$. Then there
exists a unique $\alpha \in \R^+$ so that $\Lambda^{\Phi} = \alpha
\Lambda^{\Psi}$,
\begin{equation}\label{eq6xx1}
h \bigl( L^{\Phi} (\alpha s)\bigr) = L^{\Psi} (s) \qquad \forall s \in
\R \, ,
\end{equation}
as well as
\begin{equation}\label{eq620a}
\lim\nolimits_{t\to \infty} \frac{\tau_x(t)}{t} = \alpha \qquad
\forall x \in X \setminus \{0\} \, .
\end{equation}
\end{lem}

\begin{proof}
The first two assertions clearly are correct for $d=1$. To prove
(\ref{eq620a}) for $d=1$ as well, pick any $x\ne 0$, so $h(x)\ne 0$,
and let $\Lambda^{\Phi} = [\alpha \lambda]$, $\Lambda^{\Psi} =
[\lambda]$ with $\lambda<0$ and a unique $\alpha \in
\R^+$. Then
$$
|\Phi_t x| = e^{\alpha \lambda t}|x| \, , \quad |\Psi_{\tau_x(t)} h(x)| =
e^{\lambda \tau_x (t)} |h(x)| \qquad \forall t \in \R \, .
$$
Fix $0<\beta < 1$. With the symbols $\prec$, $\succ$, and $\asymp$
used exactly as in earlier proofs, observe that
$$
e^{\lambda \tau_x(t)} \asymp |h(\Phi_t x)| \prec |\Phi_t x|^{\beta}
\asymp e^{\alpha \beta \lambda t} \, .
$$
Thus $e^{\tau_x(t)}\succ e^{\alpha \beta t}$ because $\lambda < 0$,
and consequently
\begin{equation}\label{eql66p1}
\liminf\nolimits_{t\to \infty} (\tau_x (t) - \alpha \beta t) > -
\infty \, .
\end{equation}
Similarly,
$$
e^{\alpha \lambda t} \asymp \big|
h^{-1} \bigl( \Psi_{\tau_x(t)} h(x)\bigr)
\big| \prec |\Psi_{\tau_x(t)}h(x)|^{\beta} \asymp e^{\beta \lambda
  \tau_x(t)}  \, ,
$$
thus $e^{\beta\tau_x(t)} \prec e^{\alpha t}$, and hence
\begin{equation}\label{eql66p2}
\limsup\nolimits_{t\to \infty} (\beta \tau_x (t) - \alpha t) <
\infty \, .
\end{equation}
Combining (\ref{eql66p1}) and (\ref{eql66p2}) yields
$$
\alpha \beta \le \liminf\nolimits_{t\to \infty} \frac{\tau_x(t)}{t}
\le
\limsup\nolimits_{t\to \infty} \frac{\tau_x(t)}{t} \le
\frac{\alpha}{\beta} \, ,
$$
and since $0<\beta < 1$ has been arbitrary, $\lim_{t\to
  \infty}\tau_x(t)/t = \alpha$. Thus all assertions of the Lemma are
correct for $d=1$. Assume
$d\ge 2$ from now on. Obviously, at most one $\alpha \in \R^+$ can have the desired
properties. Notice that by first applying a linear change of coordinates to obtain
the Jordan normal form of $A^{\Phi}$, and by then applying Lemmas \ref{lemr1}
and \ref{lemr4} individually to each irreducible component, together with Proposition \ref{prop3add0}, it is
straightforward to construct
$h_1 \in \cH_{1^-} (X)$ so that $\Phi \stackrel{h_1}{\cong}
\widetilde{\Phi}$, with $\widetilde{\Phi}$ generated by
$\Lambda^{\Phi}$, as well as $h_1\bigl( L^{\Phi} (s)\bigr) =
L^{\widetilde{\Phi}}(s)$ for all $s\in \R$. Similarly, there exists an 
$h_2 \in \cH_{1^-} (X)$ so that $\Psi \stackrel{h_2}{\cong}
\widetilde{\Psi}$ and $h_2\bigl( L^{\Psi} (s)\bigr) =
L^{\widetilde{\Psi}}(s)$ for all $s\in \R$, with $\widetilde{\Psi}$ generated by
$\Lambda^{\Psi}$. As a consequence,
$\widetilde{\Phi} \stackrel{\widetilde{h}}{\thicksim}
\widetilde{\Psi}$ with $\widetilde{h} = h_2 \circ h \circ h_1^{-1} \in
\cH_{1^-} (X)$, and if the assertions of the lemma can be proved for
$\widetilde{\Phi}$, $\widetilde{\Psi}$, $\widetilde{h}$ instead of
$\Phi$, $\Psi$, $h$, then also $\Lambda^{\Phi} =
\Lambda^{\widetilde{\Phi}} = \alpha \Lambda^{\widetilde{\Psi}} =
\alpha \Lambda^{\Psi}$,
$$
h\bigl( L^{\Phi} (\alpha s)\bigr)  = h \circ h_1^{-1} \bigl(
L^{\widetilde{\Phi}} (\alpha s)\bigr) = h_2^{-1} \circ \widetilde{h} \bigl(
L^{\widetilde{\Phi}} (\alpha s)\bigr) = h_2^{-1} \bigl( L^{\widetilde{\Psi}}
(s)\bigr) = L^{\Psi} (s ) \qquad \forall s \in \R \, ,
$$
and $\lim_{t\to \infty}\tau_x(t)/t=\alpha$ since
$\widetilde{h} \bigl( \widetilde{\Phi}_t h_1(x) \bigr) =
\widetilde{\Psi}_{\tau_x(t)} \widetilde{h}\bigl( h_1(x)\bigr)$ for all $t\in \R$,
$x\in X\setminus \{0\}$.
In other words, no generality is lost by assuming that $\Phi$, $\Psi$
are generated by $\Lambda^{\Phi}$, $\Lambda^{\Psi}$
respectively. With this, letting $\beta\uparrow 1$ in Theorem
\ref{lemr9} yields an $\alpha\in \R^+$ so that
$\lambda_j^{\Phi}/\lambda_j^{\Psi} = \alpha$ for all $j\in \{1,
\ldots, d\}$, that is, $\Lambda^{\Phi} = \alpha
\Lambda^{\Psi}$. Otherwise replacing $\lambda_j^{\Phi}$ by $\alpha
\lambda_j^{\Phi}$ for each $j$, assume that $\alpha = 1$ and write $\lambda_j^{\Phi}$ simply
as $\lambda_j$. Thus, it remains to show that
\begin{equation}\label{eq6yy1}
h\bigl( L^{\Phi} (s)\bigr) = L^{\Psi} (s) \qquad \forall s \in \R \, ,
\end{equation}
as well as
\begin{equation}\label{eq621a}
\lim\nolimits_{t\to \infty} \frac{\tau_x (t)}{t} = 1 \qquad \forall x
\in X \setminus \{0\} \, .
\end{equation}

To prove (\ref{eq6yy1}), notice first that $L^{\Phi} (s) = L^{\Psi}
(s) = \mbox{\rm span}  \{ e_j: \lambda_j \le s\}$ for every $s\in
\R$, and hence $h
\bigl( L^{\Phi} (s)\bigr) = \{0\} = L^{\Psi} (s)$ whenever
$s<\lambda_1$. To establish equality in (\ref{eq6yy1}) for
$s=\lambda_1$, pick any $x\in L^{\Phi} (\lambda_1) \setminus \{0\}$ 
and $0<\beta < 1$. Then $|\Phi_t x|\asymp e^{\lambda_1 t}$ and
$|h(\Phi_t x)|\prec e^{\beta \lambda_1 t}$. Also, $h(x) \in L^{\Psi}
(s)\setminus L^{\Psi} (s^-)$ for some $s<0$, and hence
$$
e^{\beta \lambda_1 t} \succ |h(\Phi_t x)| = |\Psi_{\tau_x(t)} h(x)|
\asymp e^{s \tau_x(t)} \, ,
$$
from which it follows that
\begin{equation}\label{eq6yy2}
\liminf\nolimits_{t\to \infty} \left( \tau_x(t) - \frac{\beta
    \lambda_1}{s} t\right) > - \infty \, .
\end{equation}
Now, suppose that $s>\lambda_1$. Then $h(x) + e_1 \not
\in \Psi_{\R} h(x)$, and hence $y:= h^{-1} (h(x) + e_1) \not \in
\Phi_{\R} x$. This yields
$$
e^{\lambda_1 t} \prec \mbox{\rm dist} (\Phi_t x, \Phi_{\R} y) \le
|\Phi_t x - \Phi_{\tau_y^{-1} \circ \tau_x (t)} y| \prec
|\Psi_{\tau_x(t)} h(x) - \Psi_{\tau_x(t)} h(y)|^{\beta} =
|\Psi_{\tau_x(t) }e_1|^{\beta} = e^{\beta\lambda_1 \tau_x(t)} \, ,
$$
where the left-most $\prec$ is due to Lemma \ref{lemr9A}. Thus
$e^{\lambda_1 t} \prec e^{\beta \lambda_1 \tau_x(t)}$, and
consequently
\begin{equation}\label{eq6yy3}
\limsup\nolimits_{t\to \infty}  \left( \tau_x(t) - \frac1{\beta}
  t\right) <  \infty \, .
\end{equation}
Combining (\ref{eq6yy2}) and (\ref{eq6yy3}) yields $\beta^2 \lambda_1 /s
\le 1$, that is, $s\le \beta^2 \lambda_1$ because $s<0$. Since $0<\beta <
1$ has been arbitrary, $s\le \lambda_1$, and this obviously
contradicts $s>\lambda_1$. Thus
$h(x) \in L^{\Psi} (\lambda_1)$, and since $x\in L^{\Phi}(\lambda_1)$
has been arbitrary as well, $h\bigl( L^{\Phi}
(\lambda_1) \bigr)\subset L^{\Psi} (\lambda_1)$. Interchanging the
roles of $\Phi$, $\Psi$ yields $h\bigl( L^{\Phi}
(\lambda_1) \bigr)= L^{\Psi} (\lambda_1)$. Thus equality in
(\ref{eq6yy1}) holds for all $s\le \lambda_1$.

To prepare for an induction argument, assume that $h\bigl( L^{\Phi}
(s) \bigr)= L^{\Psi} (s)$ for some $j\in \{1, \ldots , d-1\}$ and all
$s\le \lambda_j$. Similarly to before, pick $x \in L^{\Phi}
(\lambda_{j+1})$ and $0<\beta < 1$. If $x\in L^{\Phi} (\lambda_j)$
then $h(x) \in L^{\Psi} (\lambda_j) \subset L^{\Psi} (\lambda_{j+1})$.
Otherwise, $x\in L^{\Phi} (\lambda_{j+1})\setminus L^{\Phi}
(\lambda_{j+1}^-)$, and (\ref{eq6yy2}) remains valid with $\lambda_{j+1}$
instead of $\lambda_1$. In this situation, and in analogy to before, suppose that
$\lambda_{j+1} < s < 0$. Since clearly (\ref{eq6yy3}) remains valid in this situation also, the same argument as
before leads to the contradiction $s\le \lambda_{j+1}$. In summary,
$h(x) \in L^{\Psi} (\lambda_{j+1})$ for every $x\in L^{\Phi}
(\lambda_{j+1})$, and interchanging the roles of $\Phi$, $\Psi$ yields $h\bigl( L^{\Phi}
(\lambda_{j+1}) \bigr)= L^{\Psi} (\lambda_{j+1})$. In other words, $h\bigl( L^{\Phi}
(s) \bigr)= L^{\Psi} (s)$ for all $s\le \lambda_{j+1}$. Induction therefore
establishes (\ref{eq6yy1}).

To prove (\ref{eq621a}), denote by $s_1 <
\ldots < s_k<0$ all $k\le d$ {\em different\/} Lyapunov exponents of
$\Phi$, and recall that $X\setminus \{0\} = \bigcup_{j=1}^k L^{\Phi}
(s_j)\setminus L^{\Phi} (s_j^-)$. Thus, given any $x\ne 0$, there
exists a unique $j\in \{1, \ldots , k\}$ so that $x \in  L^{\Phi}
(s_j)\setminus L^{\Phi} (s_j^-)$, and hence also $h(x) \in  L^{\Psi}
(s_j)\setminus L^{\Psi} (s_j^-)$ by (\ref{eq6yy1}). It follows that
$|\Phi_t x| \asymp e^{s_j t}$ and $|\Psi_{\tau_x(t)}h(x)| \asymp e^{s_j \tau_x(t)}$.
Fix $0<\beta < 1$. Recalling that $s_j<0$, deduce from
$$
e^{s_j \tau_x(t)}  \asymp |h(\Phi_t x)| \prec |\Phi_t
x|^{\beta} \asymp e^{\beta s_j t} 
$$
that $e^{\tau_x(t)} \succ e^{\beta t}$, and hence
$\liminf_{t\to \infty} (\tau_x(t) - \beta t) > -\infty$.
Similarly, deduce from
$$
e^{s_j t} \asymp \big|
h^{-1} \bigl(
\Psi_{\tau_x(t)} h(x)
\bigr)
\big| \prec |\Psi_{\tau_x(t)} h(x)|^{\beta} \asymp e^{\beta s_j
  \tau_x(t)} 
$$
that $e^{\beta \tau_x(t)} \prec e^t$, and consequently
$\limsup_{t\to \infty} (\beta \tau_x(t) - t) < \infty$.
In summary, therefore,
$$
\beta \le \liminf\nolimits_{t\to \infty} \frac{\tau_x(t)}{t} \le
\limsup\nolimits_{t\to \infty} \frac{\tau_x(t)}{t} \le \frac1{\beta}
\, , 
$$
and since $0<\beta < 1$ has been arbitrary, $\lim_{t\to
  \infty}\tau_x(t)/t=1$. This establishes (\ref{eq621a}) and hence
completes the proof.
\end{proof}

\begin{rem}\label{rm6xx1}
(i) Lemma \ref{lem65b} carries over to unstable flows $\Phi$, $\Psi$ in an
obvious way: If  $\Phi \stackrel{h}{\thicksim} \Psi$ with $h\in
\cH_{1^-}$, then $\Lambda^{\Phi} = \alpha
\Lambda^{\Psi}$, $h \bigl( L^{\Phi^*} (\alpha s)\bigr) = L^{\Psi^*}
( s) $, and $\lim_{t\to -\infty} \tau_x(t)/t = \alpha$ for the appropriate $\alpha \in \R^+$ and all $s\in
\R$, $x\in X\setminus \{0\}$. Beyond (un)stable flows, the conclusion that $\Lambda^{\Phi} = \alpha
\Lambda^{\Psi}$ for some $\alpha \in \R\setminus \{0\}$ whenever $\Phi
\stackrel{1^-}{\thicksim} \Psi$ remains valid for {\em all\/} linear
flows $\Phi$, $\Psi$, as a consequence of Lemma
\ref{lemH2}. By contrast, the much stronger properties
(\ref{eq6xx1}) and (\ref{eq620a}) do not even carry over to {\em
  hyperbolic\/} flows.

(ii) For {\em Lipschitz\/} equivalences, the conclusions in Lemma
\ref{lem65b} take a significantly stronger form. For instance, whereas
the convergence in (\ref{eq620a}) can in general be arbitrarily slow,
it turns out that actually $\sup_{t\ge 0}|\tau_x (t) - \alpha
t|<\infty$ for every $x\in X\setminus \{0\}$ whenever $h\in \cH_1(X)$;
see \cite{BW3} for details.
\end{rem}

At long last, the scene is now set for a short

\begin{proof}[Proof of Theorem \ref{thm6all}]
Obviously (i)$\Rightarrow$(ii) by definition.

To show that (ii)$\Rightarrow$(iii), assume $\Phi
\stackrel{h}{\thicksim} \Psi$ for some $h\in \cH_{1^-}(X)$. Theorem \ref{thm6some} yields
$\{d_{\sf S}^{\Phi}, d_{\sf U}^{\Phi}\} = \{d_{\sf S}^{\Psi}, d_{\sf U}^{\Psi}\}$, and $A^{\Phi_{\sf
    C}}$, $\alpha_{\sf C} A^{\Psi_{\sf C}}$ are similar for some $\alpha_{\sf C} \in
\R \setminus \{0\}$; otherwise replacing $\Phi$ by $\Phi^*$, it again can be assumed that $(d_{\sf S}^{\Phi},
d_{\sf U}^{\Phi}) = (d_{\sf S}^{\Psi}, d_{\sf U}^{\Psi}) = :(k, m)$
and $h(X^{\Phi}_{\sf S}) = X^{\Psi}_{\sf S}$. Consider first the case
where $km\ge 1$ and $k+m\le d-1$, or equivalently $d_{\bullet}^{\Phi}=d_{\bullet}^{\Psi}>0$ for
each $\bullet\in \{{\sf S}, {\sf C}, {\sf U}\}$. Recall from the proof of
Lemma \ref{lem61} that $\Phi\stackrel{1^-}{\simeq} \widetilde{\Phi}$,
where $\widetilde{\Phi}$ is generated by $\mbox{\rm diag}\,
[\Lambda^{\Phi_{\sf S}}, A^{\Phi_{\sf C}}, \Lambda^{\Phi_{\sf U}}]$,
and similarly $\Psi\stackrel{1^-}{\simeq} \widetilde{\Psi}$, with
$\widetilde{\Psi}$ generated by $\mbox{\rm diag}\,
[\Lambda^{\Psi_{\sf S}}, A^{\Psi_{\sf C}}, \Lambda^{\Psi_{\sf
    U}}]$. Moreover, $\widetilde{\Phi}
\stackrel{1^-}{\thicksim}\widetilde{\Psi}$ and $h(E_{\sf S}) = E_{\sf
  S}$, so $\widetilde{\Phi}_{\sf S} \stackrel{1^-}{\thicksim}
\widetilde{\Psi}_{\sf S}$ as well. Since $k\ge 1$, Lemma
\ref{lem65b} yields
$\Lambda^{\Phi_{\sf S}} = \alpha_{\sf S} \Lambda^{\Psi_{\sf 
    S}}$ with the appropriate $\alpha_{\sf S}\in \R^+$. Similarly, since $m\ge 1$ also $\Lambda^{\Phi_{\sf U}} =
\alpha_{\sf U} \Lambda^{\Psi_{\sf 
    U}}$ with the appropriate $\alpha_{\sf U}\in \R^+$. Thus, $\widetilde{\Phi}$ is generated by $\mbox{\rm diag}\,
[\alpha_{\sf S} \Lambda^{\Psi_{\sf S}}, \alpha_{\sf C} A^{\Psi_{\sf
    C}}, \alpha_{\sf U } \Lambda^{\Psi_{\sf U}}]$. Note that if
$\sigma (\Phi_{\sf C}) = \sigma (\Psi_{\sf C}) = \{0\}$ then
$\alpha_{\sf S} A^{\Psi_{\sf C}}$, $\alpha_{\sf C} A^{\Psi_{\sf C}}$
are similar, so it can be assumed that $\alpha_{\sf S} = \alpha_{\sf
  C}$. Lemma \ref{lemH2}
now shows that the set $\{\alpha_{\sf S}, |\alpha_{\sf C}|, \alpha_{\sf
  U}\}$ actually is the singleton $\{\alpha\}$ for some $\alpha
\in \R^+$. Consequently,
$$
\Lambda^{\Phi} = \Lambda^{\widetilde{\Phi}} = \mbox{\rm diag}\,
[\Lambda^{\Phi_{\sf S}}, O_{\ell}, \Lambda^{\Phi_{\sf U}}] = \alpha \,
\mbox{\rm diag}\,
[\Lambda^{\Psi_{\sf S}}, O_{\ell}, \Lambda^{\Psi_{\sf U}}] =\alpha
 \Lambda^{\widetilde{\Psi}} = \alpha \Lambda^{\Psi} \, , 
 $$
 that is, $A^{\Phi}$, $\alpha A^{\Psi}$ are Lyapunov similar, and
 clearly $A^{\Phi_{\sf C}}$, $\alpha A^{\Psi_{\sf C}}$ are
 similar. Via completely analogous arguments, the same conclusions
 remain valid for the other, simpler cases where
 $d_{\bullet}^{\Phi}=0$ for some $\bullet \in \{{\sf S}, {\sf C}, {\sf
 U}\}$. Thus (iii) holds.
 
That (iii)$\Leftrightarrow$(iv) is immediate from the definition of Lyapunov similarity.

Finally, to prove that (iv)$\Rightarrow$(i), it can once again be assumed
that $(d_{\sf S}^{\Phi} , d_{\sf U}^{\Phi})=(d_{\sf S}^{\Psi}, d_{\sf U}^{\Psi})$, and hence
$\lambda_j^{\Phi_{\sf H}} = \alpha \lambda_{j}^{\Psi_{\sf H}}$ for
every $j\in \{1, \ldots , d_{\sf H}^{\Phi}\}$. Now (\ref{eq61})
simply reads $\beta \le 1$, and Lemma \ref{lem61} yields $\Phi
\stackrel{1^-}{\simeq} \Psi$. 
\end{proof}

%%%%%%%%%%%%%%%%%%%%%%%%%%%%%%%%%%%%%%%%%%%%%%%%%%%%%%%%%%%%%%%%

\section{Linear flows on complex spaces}\label{sec7}

The analysis of $\Phi \stackrel{\bigstar}{\thicksim} \Psi$ thus far has
focussed entirely on {\em real\/} flows. It is worthwhile and
straightforward to extend this analysis to linear flows on arbitrary
finite-dimensional normed spaces. In doing so, this brief section
brings the discussion of the main results to a natural conclusion.

Let $(X, \|\cdot\|)$ be a finite-dimensional normed space over $\K = \R$
or $\K = \C$. Denote by $X^{\R}$ the {\bf realification} of $X$, i.e.,
$X^{\R}$ equals $X$ as a set but is a linear space with the field of
scalars restricted to $\R$, and define $\iota_X : X \to X^{\R}$ as $\iota_X(x)
= x$. Thus, if $\K = \C$ then $\iota_X$ is an $\R$-linear bijection,
and $\dim X^{\R} = 2\dim X$; moreover, $\|\cdot\|_{X^{\R}} :=
\|\cdot\|\circ \iota_X^{-1}$ is a norm on $X^{\R}$, and $\iota_X$ is
an isometry. (Trivially, if $\K = \R$ then $X^{\R} = X$ as linear spaces, and $\iota_X = I_X$.)
Every map $h:X\to X$ induces a map $h^{\R} = \iota_X \circ h \circ
\iota_X^{-1}: X^{\R} \to X^{\R}$, and clearly
\begin{equation}\label{eq71}
h \in \cH_{\bigstar} (X) \quad \Longleftrightarrow \quad h^{\R} \in
\cH_{\bigstar} (X^{\R}) \quad \forall \bigstar \in \{0,0^+, \beta^-
\beta , \beta^+, 1^- , \, 1\}, \, 0< \beta < 1 \, ,
\end{equation}
whereas, with $J_X:= (iI_X)^{\R}$,
\begin{equation}\label{eq72}
h \in \cH_{\bigstar} (X) \quad \Longleftrightarrow \quad h^{\R} \in
\cH_{\bigstar} (X^{\R}) \enspace \mbox{\rm and} \enspace
J_X D_0h^{\R} = D_0h^{\R} J_X
\quad \forall \bigstar \in \{ {\sf diff}, {\sf lin}\} \, .
\end{equation}
Given any (not necessarily linear) flow $\varphi$ on $X$, its
realification $\varphi^{\R}$ is the flow on $X^{\R}$ with
$(\varphi^{\R})_t = (\varphi_t)^{\R}$ for all $t\in \R$. By (\ref{eq71}),
given two flows $\varphi$, $\psi$ on $X$,
\begin{equation}\label{eq73}
\varphi \stackrel{\bigstar}{\thicksim} \psi \quad \Longleftrightarrow
\quad \varphi^{\R} \stackrel{\bigstar}{\thicksim} \psi^{\R} \quad \forall \bigstar \in \{0,0^+, \beta^-
\beta , \beta^+, 1^- , 1\}, \, 0< \beta < 1 \, ,
\end{equation}
and similarly for $\varphi \stackrel{\bigstar}{\simeq} \psi $, $
\varphi \stackrel{\bigstar}{\thickapprox} \psi$ etc. For a $\K$-linear flow
$\Phi$ on $X=\K^d$, it is readily seen that all the dynamical objects
associated with $\Phi$ that have been studied in earlier sections behave
naturally under realification: For instance, $A^{\Phi^{\R}} =
(A^{\Phi})^{\R}$, and $X_{\bullet}^{\Phi^{\R}} =
(X_{\bullet}^{\Phi})^{\R} = \iota_X (X_{\bullet}^{\Phi})$ for $\bullet
\in \{ {\sf S}, {\sf C}, {\sf U}, {\sf H}\}$, as well as
$(\Phi_{\bullet})^{\R} = (\Phi^{\R})_{\bullet}$. Also, if $\K = \C$
then $\lambda_{2j-1}^{\Phi^{\R}} = \lambda_{2j}^{\Phi^{\R}} =
\lambda_j^{\Phi}$ for every $j\in \{1, \ldots, d\}$. With this, the
topological and H\"{o}lder classifications of $\K$-linear flows follow
immediately from Theorem \ref{thm6some} and \ref{thm6all}
and may be seen as the ultimate versions of
Theorems \ref{thma} and \ref{thmb}, respectively. They reveal themselves as
being {\em real\/} results, in the sense that whether or not $\Phi
\stackrel{\bigstar}{\thicksim} \Psi$ for $\bigstar \in \{0^+, 1^-\}$
is determined solely by the associated realifications $\Phi^{\R}$,
$\Psi^{\R}$. In both statements, let $X\ne \{0\}$ be a
finite-dimensional normed space over $\K$.

\begin{theorem}\label{thm71}
Let $\Phi$, $\Psi$ be $\K$-linear flows on $X$. Then the following
statements are equivalent:
\begin{enumerate}
\item $\Phi \stackrel{0^+}{\simeq}\Psi$;
\item $\Phi \stackrel{0}{\thicksim} \Psi$;
\item $\Phi^{\R} \stackrel{0^+}{\simeq}\Psi^{\R}$;
\item $\Phi^{\R} \stackrel{0}{\thicksim} \Psi^{\R}$;
  
  \item $\{\dim X_{\sf S}^{\Phi}, \dim X_{\sf U}^{\Phi}\} = \{\dim
    X_{\sf S}^{\Psi}, \dim X_{\sf U}^{\Psi}\} $, and there exists an
    $\alpha \in \R\setminus \{0\}$ so that $A^{\Phi^{\R}_{\sf
        C}}$, $\alpha A^{\Psi^{\R}_{\sf C}}$ are similar.
\end{enumerate}
\end{theorem}

\begin{proof}
Obviously (i)$\Rightarrow$(ii) and (iii)$\Rightarrow$(iv) by
definition, but also (i)$\Leftrightarrow$(iii) and
(ii)$\Leftrightarrow$(iv) by (\ref{eq73}). Furthermore, it follows
from Theorem \ref{thm6some} that
(iii)$\Leftrightarrow$(iv)$\Leftrightarrow$(v), and so all five
statements are equivalent.
\end{proof}

\begin{theorem}\label{thm72}
Let $\Phi$, $\Psi$ be $\K$-linear flows on $X$. Then the following 
statements are equivalent:
\begin{enumerate}
\item $\Phi \stackrel{1^-}{\simeq}\Psi$;
\item $\Phi \stackrel{1^-}{\thicksim} \Psi$;
\item $\Phi^{\R} \stackrel{1^-}{\simeq}\Psi^{\R}$;
\item $\Phi^{\R} \stackrel{1^-}{\thicksim} \Psi^{\R}$;
  
  \item there exists an $\alpha \in \R \setminus \{0\}$ so that
    $A^{\Phi^{\R}}$, $\alpha A^{\Psi^{\R}}$ are Lyapunov similar and $A^{\Phi^{\R}_{\sf
        C}}$, $\alpha A^{\Psi^{\R}_{\sf C}}$ are similar.
\end{enumerate}
\end{theorem}

\begin{proof}
Instead of Theorem \ref{thm6some}, simply invoke Theorem \ref{thm6all}
in the above proof of Theorem \ref{thm71}.
\end{proof}

\begin{rem}\label{rem73}
With a view on (\ref{eq73}), the reader may find it unsurprising that
the Lipschitz counterpart of Theorems \ref{thm71} and \ref{thm72},
i.e., the extension of Proposition \ref{thmc} to any
finite-dimensional normed space,
also turns out to be a real theorem in the above sense; see \cite{BW3}
for details. By contrast, the corresponding extension of Proposition
\ref{thmd} is {\em not\/}
a real theorem; see, e.g., \cite[Sec.\ 6]{BW}. In light of (\ref{eq72}), this fact may not surprise the
reader either, and it is readily
illustrated by a very simple example: Let the flows $\Phi$, $\Psi$ on
$\C$ be generated by $A = [1+i]$,
$B= [1-i]$ respectively. Then
$A^{\R}$, $B^{\R}$ are similar, so $\Phi^{\R}
\stackrel{{\sf lin}}{\simeq} \Psi^{\R}$, and hence also $\Phi
\stackrel{1^-}{\simeq}\Psi$, indeed even $\Phi
\stackrel{1}{\simeq}\Psi$, by (the complex version of) Lemma \ref{lemr3}. However, $A$, $\alpha B$ are not
similar for any $\alpha \in \R \setminus \{0\}$, so
$\Phi\cancel{\stackrel{{\sf diff}}{\thicksim}} \Psi$. Thus, the {\em
  smooth\/} equivalence of linear flows $\Phi$, $\Psi$ on $X$ is not
the same as the smooth equivalence of $\Phi^{\R}$, $\Psi^{\R}$, with
the latter being necessary for the former, but not in general sufficient.
\end{rem}

In the Introduction, all four classifications of linear flows on
$\R^d$ for $d\in \{1,2\}$ have been described. It is illuminating to
compare these to their complex counterparts. For the latter, already
the case $d=1$ hints at the peculiarity of the smooth classification
alluded to in Remark \ref{rem73}: Whereas every linear flow on
$X=\C^1$ is smoothly (in fact, holomorphically) equivalent to the flow
generated by precisely one of
$$
[0], \, [i], \, [1+ib] \qquad b\in \R \, ,
$$
it is (Lipschitz, H\"{o}lder, or topologically) equivalent to the
flow generated by $[0]$, $[i]$, or $[1]$. For $d=2$, naturally the
classification is quite a bit richer: Every linear flow on $X=\C^2$ is {\em smoothly\/} equivalent to the flow
generated by precisely one of either
\begin{equation}\label{eq7fo}
  \left[
    \begin{array}{cc}
0 & 0 \\ 0 & 0 
    \end{array}
  \right] ,
\left[
    \begin{array}{cc}
0 & 1 \\ 0 & 0 
    \end{array}
  \right] ,
  \left[
    \begin{array}{cc}
i & 1 \\ 0 & i 
    \end{array}
  \right] ,
\end{equation}
or a (necessarily unique) matrix from
$$
\left[
    \begin{array}{cc}
ia & 0 \\ 0 & i 
    \end{array}
  \right] ,
  \left[
    \begin{array}{cc}
1+ib & 1 \\ 0 & 1+ib 
    \end{array}
  \right] ,
  \left[
    \begin{array}{cc}
a +ib  & 0 \\ 0 & 1+ic 
    \end{array}
  \right] \qquad a\in [-1,1], b, c\in \R;
$$
it is {\em Lipschitz\/} equivalent to the flow
generated by precisely one of either (\ref{eq7fo}) or
$$
\left[
    \begin{array}{cc}
i|a| & 0 \\ 0 & i 
    \end{array}
  \right] ,
  \left[
    \begin{array}{cc}
1 & 1 \\ 0 & 1 
    \end{array}
  \right], 
\left[
    \begin{array}{cc}
1+ib & 1 \\ 0 & 1+ib 
    \end{array}
  \right] ,
 \left[
    \begin{array}{cc}
a  & 0 \\ 0 & 1
    \end{array}
  \right] ,
   \left[
    \begin{array}{cc}
ib  & 0 \\ 0 & 1
    \end{array}
  \right] 
  \qquad a\in [-1,1], b \in \R^+;
$$
it is {\em H\"{o}lder\/} equivalent to the flow 
generated by precisely one of either (\ref{eq7fo}) or
$$
\left[
    \begin{array}{cc}
i |a| & 0 \\ 0 & i 
    \end{array}
  \right] ,
\left[
    \begin{array}{cc}
a & 0 \\ 0 & 1 
    \end{array}
  \right] ,
   \left[
    \begin{array}{cc}
ib  & 0 \\ 0 & 1
    \end{array}
  \right]
  \qquad a\in [-1,1] , b\in \R^+;
$$
and it is {\em topologically\/} equivalent  to the flow
generated by precisely one of either (\ref{eq7fo}) or
$$
\left[
    \begin{array}{cc}
i a & 0 \\ 0 & i
    \end{array}
  \right] ,
\left[
    \begin{array}{rc}
-1 & 0 \\ 0 & 1
    \end{array}
  \right] ,
  \left[
    \begin{array}{cc}
0 & 0 \\ 0 & 1
    \end{array}
  \right] ,
  \left[
    \begin{array}{cc}
1 & 0 \\ 0 & 1
    \end{array}
  \right] ,
\left[
    \begin{array}{cc}
i & 0 \\ 0 & 1
    \end{array}
  \right]
  \qquad a\in [0,1] .
  $$
Notice in particular that while the topological classification on
$\R^2$ yields precisely six discrete classes (as seen in the
Introduction and indicated in Figure \ref{fig0a}), the corresponding
classification on $\C^2$ leads to precisely seven discrete classes,
together with the infinite family $\{i \, \mbox{\rm diag}\,[a,1]: a\in
[0,1] \}$.

\subsubsection*{Acknowledgements}

The first author was partially supported by an {\sc Nserc} Discovery  
Grant.
% The authors gratefully acknowledge helpful comments and
% suggestions by

%%%%%%%%%%%%%%%%%%%%%%%%%%%%%%%%%%%%%%%%%%%%%%%%%%%%%%%%%%%%%%%%%%%%%


\begin{thebibliography}{99}

\bibitem{Amann} H.\ Amann, {\em Ordinary differential equations: an
    introduction to non-linear analysis}, deGruyter, 1990.
  
\bibitem{ACK1} V.\ Ayala, F.\ Colonius, and W.\ Kliemann, 
  Dynamical characterization of the Lyapunov form of matrices, {\em
    Linear Algebra Appl.} {\bf 402}(2005), 272--290.
  
\bibitem{AK} V.\ Ayala and C.\ Kawan, Topological conjugacy of real
  projective flows, {\em J.\ Lond.\ Math.\ Soc.\ (2)\/} {\bf 90}(2014), 49--66.

\bibitem{BaPe} L.\ Bareira and Y.\ Pesin, {\em Introduction to Smooth
    Ergodic Theory}, Graduate Studies in Mathematics {\bf 148}, American Mathematical Society, 2013.

\bibitem{BW} A.\ Berger and A.\ Wynne, On the classification of
  finite-dimensional linear flows, {\em J.\ Dynam.\ Differential
    Equations} {\bf 32}(2020), 23--59.

\bibitem{BW3} A.\ Berger and A.\ Wynne, On Lipschitz equivalence
  of finite-dimensional linear flows, in preparation.

\bibitem{CK} F.\ Colonius and W.\ Kliemann, {\em Dynamical systems and linear algebra},
   Graduate Studies in Mathematics {\bf 158}, American Mathematical Society, 2014.

\bibitem{cop} W.A.\ Coppel, {\it Dichotomies in stability theory},
  Lecture Notes in Mathematics {\bf 629}, Springer, 1978.
   
\bibitem{DSS} A.\ Da Silva, A.J.\ Santana, and S.N.\ Stelmastchuk,
  Topological conjugacy of linear systems on Lie groups, {\em Discrete
    Contin.\ Dyn.\ Syst.} {\bf 37}(2017), 3411--3421. 

\bibitem{F} R.\ Fiorenza, {\em H\"{o}lder and locally H\"{o}lder
    Continuous Functions, and Open Sets of Class $C^k$,
    $C^{k,\lambda}$}, Birkh\"{a}user, 2016.
  
\bibitem{Hardy} G.H.\ Hardy, {\em Orders of infinity}, second edition reprinted,
  Cambridge University Press, 1954.

\bibitem{hatch} A.\ Hatcher, {\em Algebraic Topology}, Cambridge University Press, 2002.

\bibitem{He} T.\ He, Topological Conjugacy of Non-hyperbolic Linear
  Flows, preprint (2017), arXiv:1703.4413. 

\bibitem{Hei} J.\ Heinonen, {\em Lectures on Lipschitz Analysis},
  report nr.\ 100, University of Jyv\"{a}skyl\"{a}, 2005.
 
\bibitem{Irwin} M.C.\ Irwin, {\em Smooth dynamical systems}, Advanced
  Series in Nonlinear Dynamics {\bf 17}, World Scientific, 2001.
  
\bibitem{Jam} R.C.\ James, Bases in Banach Spaces, {\em Amer.\ Math.\
    Monthly\/} {\bf 89}(1982), 625--640.

\bibitem{KS} C.\ Kawan and T.\ Stender, Lipschitz conjugacy of linear flows,
{\em J. Lond. Math. Soc. (2)\/} {\bf 80}(2009),  699--715.

\bibitem{Kuiper} N.H.\ Kuiper, The topology of the solutions of a
  linear differential equation on {$R^{n}$}, p.\ 195--203 in: {\em
    Manifolds--Tokyo 1973 (Proc.\ Internat.\ Conf., Tokyo, 1973)},
  Univ.\ Tokyo Press, 1975. 

\bibitem{Ladis} N.N.\ Ladis, Topological equivalence of linear flows,
  {\em Differencialnye Uravnenija\/} {\bf 9} (1973),
  1222--1235.

\bibitem{LZ} J.\ Li, and Z.\ Zhang, Topological classification of
  linear control systems---an elementary analytic approach, {\em
    J.\ Math.\ Anal.\ Appl.} {\bf 402}(2013), 84--102.

\bibitem{MM} P.D.\ McSwiggen and K.R.\ Meyer, Conjugate phase portraits
  of linear systems, {\em Amer. Math. Monthly\/} {\bf 115}(2008), 596--614.

\bibitem{R} C.\ Robinson, {\em Dynamical Systems. Stability, Symbolic
    Dynamics, and Chaos}, CRC Press, 1995.

\bibitem{Willems} J.C.\ Willems, Topological classification and
  structural stability of linear systems, {\em J.\ Differential
    Equations\/} {\bf 35}(1980), 306--318.

\end{thebibliography}
\end{document}